\documentclass{article}

\usepackage{amsmath,amssymb,amsthm,amsfonts}
\newtheorem{theorem}{Theorem}[section]
\newtheorem{lemma}[theorem]{Lemma}
\newtheorem{proposition}[theorem]{Proposition}
\newtheorem{corollary}[theorem]{Corollary}

\theoremstyle{definition}
\newtheorem{definition}[theorem]{Definition}
\newtheorem{construction}[theorem]{Construction}
\newtheorem{remark}[theorem]{Remark}
\newtheorem{example}[theorem]{Example}

\newcommand{\Qed}{\hbox{ }\hbox{ }\hfill\rule{2.5mm}{3mm}\medskip}

\newcommand{\PGammaL}{\mathop{\mathrm{P}\Gamma\mathrm{L}}}
\newcommand{\AGammaL}{\mathop{\mathrm{A}\Gamma\mathrm{L}}}
\newcommand{\GammaL}{\mathop{\Gamma\mathrm{L}}}
\newcommand{\PGL}{\mathop{\mathrm{PGL}}}
\newcommand{\PSp}{\mathop{\mathrm{PSp}}}
\newcommand{\PSigmaL}{\mathop{\mathrm{P}\Sigma\mathrm{L}}}
\newcommand{\AGL}{\mathop{\mathrm{AGL}}}
\newcommand{\AG}{\mathop{\mathrm{AG}}}
\newcommand{\PSL}{\mathop{\mathrm{PSL}}}

\newcommand{\GF}{\mathop{\mathrm{GF}}}
\newcommand{\GL}{\mathop{\mathrm{GL}}}
\newcommand{\Cos}{\mathop{\mathrm{Cos}}}

\newcommand{\lcm}{\mathop{\mathrm{lcm}}}
\newcommand{\norml}{\vartriangleleft}
\newcommand{\Wr}{\mathop{\mathrm{wr}}}
\renewcommand{\P}{\mathcal{P}}

\newcommand{\Aut}{\mathop{\mathrm{Aut}}}

\newcommand{\Sym}{\mathop{\mathrm{Sym}}}
\newcommand{\soc}{\mathop{\mathrm{soc}}}
\newcommand{\la}{\langle}
\newcommand{\ra}{\rangle}

\title{Primitive decompositions of Johnson graphs\footnote{The first author is a Postdoctoral Researcher of the Fonds National de la Recherche Scientifique (Belgium). This research was supported under the Australian Research Council's Discovery Projects funding scheme (project number DP0449429). The second author is a recipient of an ARC Postdoctoral Fellowship while the third author holds an ARC Queen Elizabeth II Fellowship.}}

\author{Alice Devillers\\ 
 Universit\'e Libre de Bruxelles\\
D\'epartement de math\'ematiques\\
G\'eom\'etrie- CP 216\\
Boulevard du Triomphe\\
B-1050 Bruxelles
Belgique\\ \\
 Michael Giudici, Cai Heng Li and Cheryl E.~Praeger \\ 
School of Mathematics and Statistics\\
The University of Western Australia\\
35 Stirling Highway\\
Crawley WA 6009\\
Australia}
\date{}

\begin{document}
\maketitle
\begin{abstract}
A \emph{transitive decomposition} of a graph is a partition of the
edge set together with a group of automorphisms which transitively
permutes the parts. In this paper we determine all transitive
decompositions of the Johnson graphs such that the group preserving
the partition is arc-transitive and acts primitively on the parts.
\end{abstract}

\section{Introduction}
A \emph{decomposition} of a graph is a partition of the edge set with at least two parts, which we interpret as subgraphs and call the \emph{divisors} of the decomposition. If
each divisor is a spanning subgraph we call the decomposition a \emph{factorisation} and the divisors \emph{factors}. Graph decompositions and factorisations have received much attention, see for example \cite{bosak,heinrich}. Of particular interest \cite{iso1,iso10} are decompositions where the divisors are pairwise isomorphic. These are known as \emph{isomorphic decompositions}.  

A \emph{transitive decomposition} is a decomposition $\mathcal{P}$ of
a graph $\Gamma$ together with a group of automorphisms $G$ which
preserves the partition and acts transitively on the set of divisors. We
refer to $(\Gamma,\P)$ as a $G$-transitive decomposition. This is a
special class of isomorphic decompositions and a general theory has
been outlined in \cite{genpaper}. Sibley \cite{sibley} has
described all $G$-transitive decompositions of the complete graph
$K_n$ where $G$ is 2-transitive on vertices. This generalised the Cameron-Korchmaros classification in \cite{CamKorch} of the
$G$-transitive 1-factorisations of $K_n$ (that is, the factors have
valency 1) with $G$ acting 2-transitively on vertices. Note that a subgroup of
$S_n$ is arc-transitive  on $K_n$ if and only if it is 2-transitive.
Also all $G$-transitive decompositions of graphs with $G$ inducing a 
rank three product action on vertices have been determined in 
\cite{BPP}. A special class of transitive decompositions called 
\emph{homogeneous factorisations}, are the $G$-transitive decompositions $(\Gamma,\P)$ such that the kernel $M$
of the action of $G$ on  $\P$ is vertex-transitive. This implies that
each divisor is a spanning subgraph and so $\P$ is indeed a
factorisation. Homogeneous factorisations were first introduced in
\cite{LP03} for complete graphs and extended to arbitrary graphs and
digraphs in \cite{GLPP1}.

The \emph{Johnson graph} $J(n,k)$ is the graph with vertices the $k$-element subsets of an $n$-set $X$, two sets being adjacent if they have $k-1$ points in common. Note that $J(n,1)\cong K_n$ and $J(n,k)\cong J(n,n-k)$ so we always assume that $2\leq k\leq \frac{n}{2}$.  Note that $J(4,2)\cong K_{2,2,2}$ while the complement of $J(5,2)$ is the Petersen graph. All homogeneous factorisations of $J(n,k)$ were determined in \cite{Cuaresma,CGP}. Examples only exist for $J(q+1,2)$ for prime powers $q\equiv 1\pmod 4$, $J(q,2)$ and $J(q+1,3)$ for $q=2^{r^f}$ with $r$ an odd prime, and for $J(8,3)$. However, examples of transitive decompositions exist for all values of $n$ and $k$ (see Construction \ref{con:classic}). Constructions \ref{con:classic}(1) and (2) were drawn to our attention by Michael Orrison. Both constructions were used in \cite{KK} to help determine maximal subgroups of symmetric groups while Construction \ref{con:classic}(1) was used in \cite{orrisonetal} for the statistical analysis of unranked data.  

In this paper we determine all $G$-transitive decompositions of the
Johnson graphs subject to two conditions on $G$. The first is that $G$
is arc-transitive while the second is that $G$ acts primitively on the set of divisors of the
decomposition.  We call $G$-transitive decompositions for which $G$
acts primitively on the set of divisors, \emph{$G$-primitive decompositions.}
We see in Lemma \ref{lem:primitive} that any
$G$-transitive decomposition is the refinement of some $G$-primitive
decomposition. By Theorem \ref{Sn}, a subgroup $G\leqslant S_n$ acts
transitively on the set of arcs of $J(n,k)$ if and only if $G$ is
$(k+1)$-transitive, or $(n,k)=(9,3)$ and $G=\PGammaL(2,8)$. Using this, we analyse the appropriate groups to determine all primitive
decompositions.   In particular we prove the following
theorem.

\begin{theorem}
\label{thm:main}
Let $G$ be an arc-transitive group of automorphisms of
$\Gamma=J(n,k)$ where $2\leq k\leq n/2$. If $(\Gamma,\P)$ is a $G$-primitive
decomposition then one of the following holds:
\begin{enumerate}
\item the divisors are matchings or unions of cycles,
\item the divisors are unions of $K_{n-k+1}, K_{k+1}$ or $K_3$, or
\item $(\Gamma,\P)$ is given by one of the rows of Table \ref{tab:interesting}.
\end{enumerate}
\end{theorem}

\begin{table}
\begin{center}
\caption{$G$-primitive decompositions of $J(n,k)$ for Theorem \ref{thm:main}}
\label{tab:interesting}
\begin{tabular}{|llll|}
\hline
$\Gamma$ & $G$ & Divisor & Comments \\
\hline\hline 
$J(6,3)$ & $A_6$ or $\la A_6, (1,2)\tau\ra$& Petersen
graph & Example \ref{eg:An}(2) \\ 
$J(12,4)$ & $M_{12}$ & $2J(6,4)$ & Construction
\ref{con:design} and \ref{con:primitive}\\
$J(12,4)$ & $M_{12}$ & $\Sigma$ & Construction
\ref{con:M11graph}\\
$J(24,4)$ & $M_{24}$ & $J(8,4)$ & Construction
\ref{con:design}\\
$J(23,3)$ & $M_{23}$ &  $J(7,3)$ & Construction
\ref{con:design}\\
$J(11,3)$ & $M_{11}$  & $J(5,2)$ & Construction
\ref{con:design}\\
$J(11,3)$ & $M_{11}$ &  2 Petersen graphs & Construction
\ref{con:petersen}\\
$J(11,3)$& $M_{11}$ & 11 Petersen graphs &Construction
\ref{con:PSL211}(2) \\
$J(11,3)$ & $M_{11}$ &$\Pi$ & Construction \ref{con:PSL211}(1) \\
$J(9,3)$& $\PGammaL(2,8)$ & $\PSL(2,8)$-orbits & Construction
\ref{con:pgammal}(1) \\
$J(9,3)$ & $\PGammaL(2,8)$ &  Heawood graph &
Construction \ref{con:pgammal}(4)\\
$J(22,2)$ & $M_{22}$ or $\Aut(M_{22})$ &
$J(6,2)$ & Construction \ref{con:design}\\
$J(2^d,2)$, $d\geq 3$ & $\AGL(d,2)$  & $2^{d-2}K_{2,2,2}$ & Construction
\ref{con:design} and \ref{con:primitive} \\
$J(16,2)$ &$C_2^4\rtimes A_7$ & $4K_{2,2,2}$& Construction
\ref{con:design} and \ref{con:primitive}\\
$J(q+1,2)$ & 3-transitive subgroup & $J(q_0+1,2)$ & Construction \ref{con:design}\\
           &of $\PGammaL(2,q)$  &$q=q_0^r$, $r$ prime&\\
$J(q+1,2)$ & 3-transitive subgroup & $\PSL(2,q)$-orbits & Construction \ref{con:PGL1} \\
$q\equiv 1\pmod 4$&of $\PGammaL(2,q)$ &&\\
\hline
\end{tabular}
\end{center}
\end{table}

The divisor graphs $\Sigma$ and $\Pi$ of Table \ref{tab:interesting} are
investigated further in \cite{DGLP}.
Construction \ref{con:design} allows us to construct  transitive
decompositions of $J(n,k)$ with divisors isomorphic to $J(l,k)$ for any
Steiner system $S(k+1,l,n)$ and this accounts for many of the examples
in Table \ref{tab:interesting}. Further constructions of transitive
decompositions from Steiner systems are given in Section
\ref{sec:cons} and these have divisors isomorphic to unions of cliques or
matchings.

\section{General constructions}
\label{sec:cons}
First we show that the study of transitive decompositions can be reduced
to the study of primitive decompositions. We denote by $V\Gamma$, $E\Gamma$ and  $A\Gamma$, the sets of vertices, edges and arcs respectively, of the graph $\Gamma$.

\begin{construction}
{\rm 
\label{con:primitive}
Let $(\Gamma,\P)$ be a $G$-transitive decomposition and let $\mathcal{B}$
be a system of imprimitivity for $G$ on $\P$. For each
$B\in\mathcal{B}$, let $Q_B=\cup_{P\in B}P$ and let 
$\mathcal{Q}=\{Q_B\mid B\in\mathcal{B}\}$. Then $(\Gamma,\mathcal{Q})$ is a $G$-transitive decomposition.
}%
\end{construction}

\begin{lemma}
\label{lem:primitive} 
Any $G$-transitive decomposition $(\Gamma,\P)$ with $|\P|$ finite is the refinement of a $G$-primitive decomposition  $(\Gamma,\mathcal{Q})$.
\end{lemma}
\proof
If $G^{\P}$ is primitive then we are done. If not, let $\mathcal{B}$
be a nontrivial system of imprimitivity for $G$ on $\P$ with maximal block size. Then
$G^{\mathcal{B}}$ is primitive and $\P$ is a refinement of the
partition $\mathcal{Q}$ yielded by Construction \ref{con:primitive}.
Thus $(\Gamma,\mathcal{Q})$ is a $G$-primitive decomposition.
\Qed

We have the following general construction of transitive
decompositions.

\begin{construction}
{\rm 
\label{con:general}
Let $\Gamma$ be a graph with an  arc-transitive group $G$ of
automorphisms.  Let $e$ be 
an edge of $\Gamma$ and suppose that there exists a subgroup $H$ of
$G$ such that $G_e<H<G$. Let $P=e^H$ and $\P=\{P^g\mid g\in G\}$. 
}%
\end{construction}
\begin{lemma}
\label{lem:general}
Let $(\Gamma,\P)$ be obtained as in Construction
$\ref{con:general}$. Then $(\Gamma,\P)$ is a $G$-transitive
decomposition. Conversely, every $G$-transitive decomposition with $G$
arc-transitive arises in such a manner. Moreover, if the subgroup $H$ is maximal in $G$, then  $(\Gamma,\P)$ is a $G$-primitive decomposition.
\end{lemma}
\proof
Since $G$ is arc-transitive and $G_e<H<G$, then $\P$ is a partition of
$E\Gamma$ which is preserved by $G$ and such that $G^{\P}$ is
transitive. Thus $(\Gamma,\P)$ is a $G$-transitive
decomposition. Conversely, let $(\Gamma,\P)$ be a $G$-transitive
decomposition such that $G$ is arc-transitive. Let $e$ be an edge of
$\Gamma$ and $P$ the divisor containing $e$. Since $\P$ is a
system of imprimitivity for $G$ on $E\Gamma$ it follows that for
$H=G_P$ we have $G_e<H<G$ and $P=e^H$. Moreover, $\P=\{P^g\mid g\in
G\}$ and so $(\Gamma,\P)$ arises from Construction \ref{con:general}.
The last statement follows from the fact that $H$ is the stabiliser in $G$ of the divisor $P$.
\Qed

\begin{remark}
\label{rem:2ways}
{\rm Lemma \ref{lem:general} implies that there are two possible ways to
determine all  $G$-transitive decompositions such that the
divisor stabilisers are in a given conjugacy class $H^G$ of subgroups of
$G$. One is to fix an edge $e$ and run 
over all subgroups conjugate to $H$ which contain the stabiliser of
$e$. Note that different conjugates may give different partitions. The
second is to run over all edges whose stabiliser is contained in
$H$. Again, different edges may give different partitions.
}%
\end{remark}

We say that two decompositions $(\Gamma,\P_1)$ and $(\Gamma,\P_2)$ are \emph{isomorphic}
if there exists $g\in\Aut(\Gamma)$ such that $\P_1^g=\P_2$. If both are $G$-transitive decompositions, then they are \emph{isomorphic $G$-transitive decompositions} if there is such an element $g\in N_{\Aut(\Gamma)}(G)$. The
following lemma gives us a condition for determining when different
conjugates give the same decomposition.  
\begin{lemma}
\label{lem:iso}
Let $(\Gamma,\P_1)$, $(\Gamma,\P_2)$ be two $G$-transitive
decompositions with $G$ arc-transitive. 
\begin{enumerate}
\item Let $e$ be an edge of $\Gamma$
and $P_1,P_2$ be the divisors of $\P_1$, $\P_2$ respectively that contain $e$. 
If there exists an automorphism $g\in N_{\Aut(\Gamma)}(G)$ fixing $e$
such that $G_{P_1}^g=G_{P_2}$ then $(\Gamma,\P_1)$ and $(\Gamma,\P_2)$ are isomorphic.
\item  Let $e_1,e_2$ be two edges of $\Gamma$ with divisors $P_1=e_1^H$ and $P_2=e_2^H$ of $\P_1$, $\P_2$ respectively. If there exists an automorphism $g\in N_{\Aut(\Gamma)}(G)$ mapping $e_1$ onto $e_2$
such that $H^g=H$ then $(\Gamma,\P_1)$ and $(\Gamma,\P_2)$ are isomorphic.
\end{enumerate}
\end{lemma}
\proof
\begin{enumerate} \item By Lemma \ref{lem:general},
$P_1=e^{G_{P_1}}$ and $P_2=e^{G_{P_2}}$. Thus
$P_2=e^{g^{-1}G_{P_1}g}=e^{G_{P_1}g}=P_1^g$. Moreover,
$\P_2=P_2^G=(P_1^g)^G=(P_1^G)^g=\P_1^g$ and so $(\Gamma,\P_1)$ and $(\Gamma,\P_2)$ are isomorphic.
\item We have $P_2=e_2^H=(e_1^g)^H=(e_1^H)^g=P_1^g$. Hence we get the same conclusion.\Qed
\end{enumerate}

We also have the following useful lemma.

\begin{lemma}\label{lem:restrict}
Let $(\Gamma,\P)$ be a $G$-primitive decomposition, with $H$ the stabiliser of a divisor $P$. If $L\leqslant G$ is such that $L\not\leqslant H$, $L$ is arc-transitive on $\Gamma$ and $L\cap H$ is maximal in $L$, then  $(\Gamma,\P)$ is a $L$-primitive decomposition.
\end{lemma}
\proof Since $L$ is arc-transitive and contained in $G$, it follows that $L$ acts transitively on $\P$.  Moreover, since $H\cap L$ is the stabiliser in $L$ of a part, it follows that $L$ acts primitively on $\P$.
\Qed

We now describe some general methods for constructing transitive decompositions of Johnson graphs.

\begin{construction}
\label{con:classic}
{\rm Let $X$ be an $n$-set.
\begin{enumerate}
\item For each $(k-1)$-subset $Y$ of $X$, let $P_Y$ be the complete subgraph of $J(n,k)$ whose vertices are all the $k$-subsets containing $Y$. 
Then 
$${\mathcal P}_{\cap}=\{P_Y|Y \text{ is a $(k-1)$-subset of $X$}\}$$ 
is a decomposition of $J(n,k)$ with ${n\choose k-1}$ divisors, each 
isomorphic to $K_{n-k+1}$.

\item For each $(k+1)$-subset $W$ of $X$, let $Q_W$ be the complete subgraph
whose vertices are all the $k$-subsets contained in $W$. Then
$${\mathcal P}_{\cup}=\{Q_W|W \text{ is a $(k+1)$-subset of $X$}\}$$ 
is a decomposition of $J(n,k)$ with ${n\choose k+1}$ divisors, each isomorphic to
$K_{k+1}$.

\item For each $\{a,b\}\subseteq X$, let
$$M_{\{a,b\}}=\Big\{\big\{\{a\}\cup Y,\{b\}\cup Y\}\big\} \mid Y \text{ a
  $(k-1)$-subset of } X\backslash\{a,b\}\Big\}.$$
Then $$\mathcal{P}_{\ominus}=\{M_{\{a,b\}}\mid \{a,b\}\subseteq X\}$$ 
is a decomposition of $J(n,k)$ with ${n\choose 2}$ divisors, each of which
  is a matching with ${n-2 \choose k-1}$ edges.
\end{enumerate}
}%
\end{construction}

Given two sets $A$ and $B$ we denote the \emph{symmetric difference} of $A$
and $B$ by $A\ominus B$.
\begin{lemma}\label{lem:classic} Let $G\leqslant S_n$ such that 
$\Gamma =J(n,k)$ is $G$-arc-transitive. Let $A$ and $B$ be two
  adjacent vertices of $\Gamma$. Then $(\Gamma, {\mathcal P}_{\cap})$,
$(\Gamma, {\mathcal P}_{\cup})$,  $(\Gamma, {\mathcal P}_{\ominus})$ are $G$-transitive decompositions. Moreover, if $G_{A\cap B}$, $G_{A\cup B}$, or
$G_{A\ominus B}$ respectively is maximal in $G$, then the decomposition is $G$-primitive.
\end{lemma}
\proof Since $P_Y^g=P_{Y^g}$, $Q_W^g=Q_{W^g}$ and $M_{\{a,b\}}^g=M_{\{a,b\}^g}$, it follows that $G$ preserves ${\mathcal P}_{\cap}$, ${\mathcal P}_{\cup}$ and ${\mathcal P}_{\ominus}$. Since $G$ is arc-transitive, all three decompositions are $G$-transitive. The divisor of $\P_\cap$, $\P_\cup$ or $P_\ominus$ containing $\{A,B\}$ is $P_{A\cap B}$, $Q_{A\cup B}$ or $M_{A\ominus B}$ respectively, and the stabiliser of this divisor is $G_{A\cap B}$, 
$G_{A\cup B}$, or  $G_{A\ominus B}$ respectively. The last assertion follows.
\Qed

Another method for constructing transitive decompositions of $J(n,k)$
is to use  Steiner systems with multiply transitive automorphism groups.
A \emph{Steiner system} $S(t,k,v)=(X,\mathcal{B})$ is a collection $\mathcal{B}$ of $k$-subsets
(called \emph{blocks}) of a $v$-set $X$ such that each $t$-subset of $X$ is
contained in a unique block.

\begin{construction}
\label{con:design}
{\rm Let $\mathcal{D}=(X,\mathcal{B})$ be an $S(k+1,l,n)$ Steiner system with
automorphism group $G$ such that $G$ is transitive on
$\mathcal{B}$. For each $Y\in\mathcal{B}$, let $P_Y$ be the subgraph of 
$J(n,k)$ whose vertices are the $k$-subsets in $Y$ and let
$\P=\{P_Y\mid Y\in \mathcal{B}\}$.
}%
\end{construction}

\begin{lemma}
\label{lem:design}
The pair $(J(n,k),\P)$ yielded by Construction $\ref{con:design}$ is a
$G$-transitive decomposition with divisors isomorphic to $J(l,k)$.
Moreover, the decomposition is $G$-primitive if and only if the stabiliser
of a block of $\mathcal{D}$ is maximal in $G$.
\end{lemma}
\proof
Let $\{A,B\}$ be an edge of  $J(n,k)$. Then $A\cup B$ has size $k+1$ and
so is contained in a unique block $Y$ of $\mathcal{D}$, and hence
$\{A,B\}$ is contained in a unique part $P_Y$ of $\P$. Thus $(J(n,k),\mathcal{P})$ is a decomposition. Since $G$ is
transitive on $\mathcal{B}$ the pair $(J(n,k),\P)$ is $G$-transitive. Moreover, each $P_Y$ consists of all $k$-subsets of the
$l$-set $Y$ and so is isomorphic to $J(l,k)$. Since the stabiliser in $G$ of $P_Y$ is $G_Y$, the last statement follows.
\Qed

\begin{construction}
\label{con:design2}
{\rm Let $\mathcal{D}=(X,\mathcal{B})$ be an $S(k+1,l,n)$ Steiner system
with automorphism group $G$. Let $i=l-k-1$ and suppose that $G$ is
$i$-transitive on $X$. For each $i$-subset $Y$ of $X$ let
$$P_Y=\{\{A,B\}\mid |A|=|B|=k,|A\cap B|=k-1 \text{ and } A\cup B\cup Y \in
\mathcal{B}\}.$$
Define
$$\P=\{P_Y\mid Y \text{an $i$-subset of } X\}.$$
}%
\end{construction}

\begin{lemma}
\label{lem:design2}
The pair $(J(n,k),\P)$ yielded by Construction $\ref{con:design2}$ is a
$G$-transitive decomposition with divisors isomorphic to $mK_{k+1}$,
where $m$ is the number of blocks of $\mathcal{D}$ containing an
$i$-set.  Moreover, the decomposition is $G$-primitive if and only if the
stabiliser of an $i$-set is maximal in $G$.
\end{lemma}
\proof
Let $\{A,B\}$ be an edge of $J(n,k)$. Then $A\cup B$ is contained in a
unique block $W$ of $\mathcal{D}$ and the unique part of $\P$ containing $\{A,B\}$ is $P_Y$ where $Y=W\backslash (A\cup B)$.  Each block containing $Y$
contributes a copy of $J(k+1,k)\cong K_{k+1}$ to $P_Y$, and since each
$(k+1)$-subset is in a unique block, no two blocks containing $Y$ share a vertex of $P_Y$. Hence the $m$ copies of $K_{k+1}$ in $P_Y$,
are pairwise vertex-disjoint, that is $P_Y\cong mK_{k+1}$.
Since $G$ is $i$-transitive, it
follows that $(J(n,k),\P)$ is a $G$-transitive decomposition. Since the stabiliser in $G$ of $P_Y$ is $G_Y$, the last statement follows.
\Qed

\begin{construction}
\label{con:design3}
{\rm Let $\mathcal{D}=(X,\mathcal{B})$ be an $S(k+1,k+2,n)$ Steiner system
with automorphism group $G$ such that $G$ acts $3$-transitively on
$X$. For each $3$-subset $Y$ of $X$, let 
$$P_Y=\Big\{\big\{ Z\cup\{u\},Z\cup\{v\}\big\}\mid |Z|=k-1,Z\cup
  Y\in\mathcal{B}, u,v\in Y\Big\}$$
and let $\P=\{P_Y\mid Y \text{ a $3$-subset of }$X$\}$.
}%
\end{construction}

\begin{lemma}
The pair $(J(n,k),\P)$ yielded by Construction $\ref{con:design3}$ is a
$G$-transitive decomposition with divisors isomorphic to $mK_3$, where
$m$ is the number of blocks of $\mathcal{D}$ containing a given
$3$-set. Moreover, the decomposition is $G$-primitive if and only if the
stabiliser of a $3$-subset is maximal in $G$.
\end{lemma}
\proof
Let $\{A,B\}$ be an edge of $J(n,k)$. Then $A\cup B$ is contained in a
unique block $W$ of $\mathcal{D}$ and the unique part of $\P$ containing $\{A,B\}$ is
$P_Y$ where $Y=W\backslash (A\cap B)$.  Each block containing $Y$
contributes a copy of $K_3$ to $P_Y$, and since each $(k+1)$-subset is
in a unique block, no two blocks containing $Y$ share a vertex of $P_Y$. Hence the $m$ copies of $K_3$ in $P_Y$ are pairwise vertex-disjoint,
that is, $P_Y\cong mK_3$. Since $G$ is $3$-transitive, it follows that
$(J(n,k),\P)$ is a $G$-transitive decomposition.
 Since the stabiliser in $G$ of $P_Y$ is $G_Y$, the last statement follows. 
 \Qed

\begin{construction}
\label{con:design4}
{\rm Let $\mathcal{D}=(X,\mathcal{B})$ be an $S(k+1,k+2,n)$ Steiner system
with $k$-transitive automorphism group $G$. For each $k$-subset $Y$ of
$X$ let 
$$P_Y=\Big\{\big\{\{u\}\cup Z,\{v\}\cup Z\big\}\mid
Y\cup\{u,v\}\in\mathcal{B}, Z\subset Y, |Z|=k-1 \Big\}$$
and let $\P=\{P_Y\mid Y \text{ a k-subset of }$X$\}$.
}%
\end{construction}
\begin{lemma}
The pair $(J(n,k),\P)$ yielded by Construction $\ref{con:design4}$ is a
$G$-transitive decomposition with divisors isomorphic to $mkK_2$, where
$m$ is the number of blocks of $\mathcal{D}$ containing a given
$k$-set. Moreover, the decomposition is $G$-primitive if and only if the
stabiliser of a $k$-subset is maximal in $G$.
\end{lemma}
\proof
Let $\{A,B\}$ be an edge of $J(n,k)$. Then $A\cup B$ is contained in a
unique block $W$ of $\mathcal{D}$ and the unique part of $\P$ containing $\{A,B\}$ is
$P_Y$ where $Y=W\backslash (A\ominus B)$.  Each block containing $Y$
contributes a copy of $kK_2$ to $P_Y$, and since each
$(k+1)$-subset is in a unique block, no two blocks containing $Y$ share a vertex of $P_Y$. Hence the $m$ copies of $kK_2$ in $P_Y$,
are pairwise vertex-disjoint, that is $P_Y\cong mkK_2$.
Since $G$ is $k$-transitive, it
follows that $(J(n,k),\P)$ is a $G$-transitive decomposition. Since the stabiliser in $G$ of $P_Y$ is $G_Y$, the last statement follows.
\Qed

We end this section with a standard construction of arc-transitive
graphs.

Let $G$ be a group with corefree subgroup $H$ and let $g\in G$ such 
that $g^2\in H$ and $g\notin N_G(H)$. Define the graph 
$\Gamma=\Cos(G,H,HgH)$ with vertex set the set of right cosets of $H$ 
in $G$ and $Hx$ adjacent to $Hy$ if and only if $xy^{-1}\in HgH$. Then 
$G$ acts faithfully and arc-transitively on $\Gamma$ by right 
multiplication. We have the following lemma, see for example \cite{FP}.

\begin{lemma}
\label{lem:cosetgraph}
Let $\Gamma$ be a $G$-arc-transitive graph with adjacent vertices $v$ and $w$. Let $H=G_v$, and let $g\in G$ interchange $v$ and $w$. Then $\Gamma\cong \Cos(G,H,HgH)$. The connected component of $\Gamma$ containing $v$ consists of all cosets of $H$ contained in $\la H,g\ra$. In particular, $\Gamma$ is connected if and only if $\la H,g\ra=G$.
\end{lemma}


\section{Groups}
In this section, we determine the groups $G$ such that $J(n,k)$ is 
$G$-vertex-transitive and $G$-arc-transitive.

\begin{theorem}{\rm \cite[Theorem 9.1.2]{BCN}}
\label{thm:autgamma}
Let $n,k$ be positive integers and let $\Gamma=J(n,k)$. If $n>2k$ then $\Aut(\Gamma)=S_n$ with the action induced from the action of $S_n$ on $X$. For $n=2k\geq 4$, $\Aut(\Gamma)=S_n\times S_2=\la S_n,\tau \ra$ where $\tau$ acts on $V\Gamma$ by complementation in $X$.
\end{theorem}

Given a subset $A$ of $X$ we denote the complement of $A$ in $X$ by
$\overline{A}$. Also, if $|X|=n$ and $|A|=k$ then $\Gamma(A)$ denotes
the set of neighbours of $A$ in the graph $J(n,k)$, that is, vertices $B$ such that $\{A,B\}$ is an edge.

\begin{lemma}{\rm \cite[Proposition 3.2]{Cuaresma}}
\label{lemcond} 
Let $\Gamma=J(n,k)$ and $G\leqslant S_n$. The graph $\Gamma$ is $G$-arc-transitive if and only if $G$ is $k$-homogeneous on $X$ and, for a $k$-subset $A$, $G_A$ is transitive on $A\times \overline{A}$.
\end{lemma}
\proof
Note that $G$ is arc-transitive if and only if $G$ is
vertex-transitive and $G_A$ is transitive on $\Gamma(A)$. By definition,
$\Gamma$ is $G$-vertex-transitive if and only if $G$ is
$k$-homogeneous on $X$. Moreover, $G_A$ is transitive on $\Gamma(A)$  if and
only if $G_A$ is independently transitive on the set of $(k-1)$-subsets
of $A$ and on $\overline{A}$, that is, if and only if $G_A$ is
transitive on $A\times \overline{A}$.
\Qed

\begin{corollary}\label{Sncor} If $G\leqslant S_n$ is $(k+1)$-transitive, then $\Gamma$ is $G$-arc-transitive. If $\Gamma$ is $G$-arc-transitive and $G\leqslant S_n$, then $G$ is $k$- and $(k+1)$-homogeneous.
\end{corollary}

\begin{theorem}\label{Sn}
Let $n\geq 2k\geq 4$ and $G\leqslant S_n$. The graph $\Gamma=J(n,k)$ is $G$-arc-transitive if and only if $G$ is $(k+1)$-transitive on $X$ or $k=3$, $n=9$, and $G=\PGammaL(2,8)$.
\end{theorem}
\proof
If $G$ is $(k+1)$-transitive, then by Corollary \ref{Sncor}, $\Gamma$
is $G$-arc-transitive. If $k=3$, $n=9$, and $G=\PGammaL(2,8)$, then it is easy
to check that $G$ is arc-transitive.

Suppose now that $\Gamma$ is $G$-arc-transitive. By Corollary \ref{Sncor}, $G$ is $k$- and $(k+1)$-homogeneous on $X$. If $G$ is not $(k+1)$-transitive, then, by \cite{Kantor,LiWa} either $2k\leq n\leq 2k+1$, or $2\leq k\leq 3$ and $G$ is one of a small number of groups.

Suppose first that $k=2$. (This is an improvement on the proof of
\cite[Proposition 3.3]{Cuaresma}.) Since $G$ is $3$-homogeneous, it is
transitive on $X$. For $A=\{a,b\}$, Lemma \ref{lemcond} implies that $G_A$ is
transitive on $A\times \overline{A}$. Therefore using elements of $G_A$ we can map $(a,c)$ onto $(a,d)$ for any $c,d\in  \overline{A}$,
and so $G_{a,b}$ is transitive on  $\overline{A}$. Similarly,
$G_{a,c}$ is transitive on  $\overline{\{a,c\}}$ for any
$c\in \overline{\{a,b\}}$. Hence $G_a$ is transitive on
$\overline{\{a\}}$ and so $G$ is 3-transitive on $X$.

Next suppose that $k=3$. If $G$ is not 4-transitive then either
$n=6,7$, or by \cite{Kantor}, $G$ is one of $\PGL(2,8)$,
$\PGammaL(2,8)$ (with $n=9$),  or $\PGammaL(2,32)$ (with $n=33$).  Let $A=\{a,b,c\}$ and suppose that $G\neq \PGammaL(2,8)$.

Suppose first that $G=\PGL(2,8)$. Then $G_A\cong S_3$ and
$G_{A,a}=C_2$. Hence $G$ does not satisfy the arc-transitivity
condition given in Lemma \ref{lemcond}.  Next suppose that
$G=\PGammaL(2,32)$. Then $|G_{A,a}|=10$ and so again Lemma
\ref{lemcond} implies that $G$ is not arc-transitive.

If $n=6$, the only 3-homogeneous and 4-homogeneous group which is not
4-transitive is $\PGL(2,5)$. However, this does not satisfy the
condition in Lemma \ref{lemcond} for arc-transitivity. There are no
3-homogeneous and 4-homogeneous groups of degree 7 which are not
4-transitive.

Next suppose that $k=4$. If $G$ is not 5-transitive, then $n=8$ or $9$.
Since $G$ is $4$-homogeneous and 5-homogeneous, either $G$ is
$4$-transitive, or $G$ is one of $\PGL(2,8)$,
$\PGammaL(2,8)$. However, these two groups are not arc-transitive as the
stabiliser of a $4$-subset $A$ also stabilises a point in $\overline{A}$.
The only 4-transitive groups of degree $n$ are $A_n$ and $S_n$
and they are also $5$-transitive.

If $k=5$ and $G$ is not 6-transitive, then $n=10$ or $11$. Since $G$ is $5$-homogeneous it is
$5$-transitive and so $G$ contains $A_n$. Thus $G$ is also
$6$-transitive. Finally, let $k\geq 6$. Since $G$ is $k$-homogeneous it is
$k$-transitive. The only $k$-transitive groups for $k\geq 6$ are
$A_n$ and $S_n$, which are also $(k+1)$-transitive.
\Qed

We need a couple of results for the case $n=2k$.
\begin{theorem}
\label{thm:2karctrans}
Let $\Gamma=J(2k,k)$ and suppose that $G\leqslant \Aut(\Gamma)=S_{2k}\times \la \tau\ra$ and $\Gamma$ is $G$-arc-transitive. Then either $G\cap S_{2k}$ is arc-transitive on $\Gamma$, or $k=2$, $G=\la A_4,(1,2)\tau\ra$ and $G\cap S_4=A_4$ has two orbits on arcs.
\end{theorem}
\proof
Let $\hat{G}=G\cap S_{2k}$. If $\hat{G}=G$, we are done. Hence we can
assume $\hat{G}$ is an index $2$ subgroup of $G$. The graph $\Gamma$
is connected and is not bipartite, as it contains 3-cycles.  It
follows that $\hat{G}$ cannot have two orbits on vertices and so
$\hat{G}$ is vertex-transitive.

Suppose that $\hat{G}$ is not arc-transitive, and hence has two orbits
of equal size on $A\Gamma$. Let $(A,B)\in A\Gamma$. Then $\hat{G}_{(A,B)}\leqslant G_{(A,B)}$ and $|G_A:G_{(A,B)}|=|\Gamma(A)|=k^2=2|\hat{G}_A:\hat{G}_{(A,B)}|=|G_A:\hat{G}_{(A,B)}|$.
Hence $\hat{G}_{(A,B)}=G_{(A,B)}$ and $k$ is even.

Suppose first that $k\geq 6$. Since $\hat{G}$ is transitive on
$V\Gamma$, $\hat{G}$ is $k$-homogeneous and therefore also
$k$-transitive. Hence $A_{2k}\leqslant\hat{G}$, and so $\hat{G}$ is
$(k+1)$-transitive. It follows from Theorem \ref{Sn} that $\hat{G}$ is
transitive on $A\Gamma$, which is a contradiciton. Thus $k=2$ or 4.

If $k=4$, then $\hat{G}$ is $k$-homogeneous. The only 4-homogeneous
groups of degree 8 contain $A_8$, and so are also $5$-transitive. By
Theorem \ref{Sn}, $\hat{G}$ is transitive on
$A\Gamma$ in this case, and so $k=2$.

Since $\hat{G}$ is transitive on $V\Gamma$ and $(n,k)=(4,2)$ we have
that $6$ divides $|\hat{G}|$.  Since $\hat{G}$ is
2-homogeneous it follows that $A_4\leqslant \hat{G}$. Moreover, $S_4$
is arc-transitive and so $\hat{G}=A_4$.
There are two groups $G\leqslant S_n\times S_2$ such that
$\hat{G}=A_4$ and is of index 2 in $G$, namely $\la A_4,\tau \ra$
and $\la A_4,(1,2)\tau \ra$. It is easy to check that the second group
is transitive on $A\Gamma$ but not the first one.
\Qed

We also have the following theorem about primitivity.

\begin{theorem}
\label{thm:2kreducetoSn}
Let $\Gamma=J(2k,k)$ and $G\leqslant \Aut(\Gamma)=S_{2k}\times \la
\tau\ra$ such that both $G$ and $G\cap S_{2k}$ are
arc-transitive. Suppose that $(\Gamma,\P)$ is a $G$-primitive
decomposition. Then $(\Gamma,\P)$ is also $(G\cap S_{2k})$-primitive.
\end{theorem}
\proof
Let $\hat{G}=G\cap S_{2k}$, let $H$ be the stabiliser in $G$ of a
divisor and $\hat{H}=H\cap \hat{G}=H\cap S_{2k}$.  We may suppose that
$G\neq\hat{G}$. Moreover, as $\hat{G}$ is arc-transitive it acts
transitively on $\P$ and so $\hat{G}\not\leqslant H$.  Since $H$ is
maximal in $G$  it follows that $|H:\hat{H}|=2$.

Suppose first that $G=\hat{G}\times\la\tau\ra$. Now
$H=\la \hat{H},\sigma\tau\ra$ for some $\sigma\in \hat{G}$. Since $\hat{H}\norml H$, the element $\sigma\tau$ (and hence also $\sigma$) normalises $\hat{H}$ and $\hat{H}$ contains $(\sigma\tau)^2=\sigma^2$. This implies that 
$H\leqslant\la\hat{H},\sigma\ra\times\la\tau\ra\leqslant G$. Since
$H$ is maximal in $G$, either $H=\la\hat{H},\sigma\ra\times\la\tau\ra$ or $\la\hat{H},\sigma\ra\times\la\tau\ra=G$. The first implies that $\sigma\in\hat{H}$ and hence $H=\hat{H}\times \la\tau\ra$. Thus $\hat{H}$ is maximal in $\hat{G}$ and so by Lemma \ref{lem:restrict}, $\mathcal{P}$ is $\hat{G}$-primitive.  On the other hand, the second implies 
$\hat{G}=\la \hat{H},\sigma\ra$. Since $\sigma^2\in\hat{H}$,
we have $|\mathcal{P}|= |\hat{G}:\hat{H}|=2$ and so again $\hat{G}$ is primitive on $\P$.

Suppose now that $G=\la \hat{G},\sigma\tau\ra$ for some
$\sigma\in S_{2k}\backslash\{1\}$ and $\tau\notin G$. Then $\sigma$ normalises $\hat{G}$ and
$\sigma^2\in\hat{G}$.  Also, as $\tau\notin G$, we have $\sigma\notin\hat{G}$ and in particular $\hat{G}\neq S_{2k}$. By Theorem \ref{Sn} and the fact that $n=2k$, 
the classification of $(k+1)$-transitive groups (see for example
\cite[pp194--197]{Cameron}) implies that
$\hat{G}=A_{2k}$ and $k\geq 3$. Let $\phi:S_{2k}\times\la\tau\ra \rightarrow S_{2k}$ be the
projection of $\Aut(\Gamma)$ onto $S_{2k}$. Then $\phi_{\mid G}$ is an isomorphism.
Moreover, for an edge $\{A,B\}$ contained in the divisor stabilised by $H$, $\phi(G_{A,B})= S_{k-1}\times S_{k-1}$. 
Since $k\geq 3$, there is a transposition in $\phi(G_{A,B})$ and so
by \cite[Theorem 13.1]{wielandt} and since $\phi(G_{A,B})\subseteq \phi(H)$, $\phi(H)$ is not primitive. It follows that
$\phi(H)$ is a maximal intransitive subgroup of $S_{2k}$ or a maximal
imprimitive subgroup of $S_{2k}$ preserving a partition into at most 3
parts.  Thus by \cite{LPS} and since $\hat{H}=\phi(H) \cap A_{2k}$, it
follows that $\hat{H}$ is a maximal subgroup of
$\hat{G}=A_{2k}$. Hence again $\hat{G}$ is primitive on $\P$. 
\Qed

\section{Alternating and symmetric groups}\label{sec:AnSn}

We have already seen the $S_n$-transitive decompositions
$\P_\cap,\P_\cup$ and $\P_\ominus$. Since $n\geq 2k$ it follows that
$S_n$ always acts primitively on $\P_\cap$. Also, $S_n$ acts
primitively on $\P_\cup$ if and only if $n\neq 2k+2$. When $n=2k+2$, applying Construction \ref{con:primitive} to $\P_\cup$ we
obtain an $S_n$-primitive decomposition with divisors isomorphic to
$2K_{k+1}$. Finally $S_n$ acts primitively on $\P_\ominus$ if and only
if $(n,k)\neq (4,2)$. This justifies the first four lines of Table \ref{tab:Sn} below. We also have the following two examples.
\begin{example}
\label{eg:Sn}
{\rm \begin{enumerate}
\item Let $G=S_4$, $H=\la (1,2,3,4),(1,3)\ra\cong D_8$, $A=\{1,2\}$
  and $B=\{2,3\}$.   Then $P=\{A,B\}^H$ is the 4-cycle
 $$\Big\{\big\{\{1,2\},\{2,3\}\big\},\big\{\{2,3\},\{3,4\}\big\},
\big\{\{3,4\},\{1,4\}\big\},\big\{\{1,4\}, \{1,2\}\big\}\Big\}.$$
 Since $G_{\{A,B\}}=\la (1,3)\ra$ we have
$G_{\{A,B\}}<H<G$ and so by Lemma \ref{lem:general} $( (J(4,2),\P)$ is 
a $G$-primitive decomposition with $\P=\{P^g\mid g\in G\}.$ 

\item Let $G=S_6$ and  $H$ be the stabiliser in $G$ of the partition
  $$\big\{\{1,4\},\{2,3\},\{5,6\}\big\}$$ of $\{1,\ldots,6\}$. Let
  $A=\{1,2,3\}$ and $B=\{2,3,4\}$. Then $P=\{A,B\}^H$ is the matching
$$\Big\{\big\{\{1,2,3\},\{2,3,4\}\big\},\big\{\{2,5,6\},\{3,5,6\}\big\},\big\{\{1,4,5\},\{1,4,6\}\big\},$$
$$\big\{\{1,5,6\},\{4,5,6\}\big\},\big\{\{2,3,5\},\{2,3,6\}\big\},\big\{\{1,4,2\},\{1,4,3\}\big\}
\Big\}.$$ 
Since $G_{\{A,B\}}<H<G$ it follows from Lemma \ref{lem:general} that
$((J(6,3),\P)$ is a $G$-primitive decomposition with $\P=\{P^g\mid g\in G\}.$
\end{enumerate}
}%
\end{example}

We have now constructed all the $S_n$-primitive decompositions in
Table \ref{tab:Sn}. It remains to prove that these are the only ones.

\begin{table}
\begin{center}
\caption{$S_n$-primitive decompositions of $J(n,k)$}
\label{tab:Sn}
\begin{tabular}{|llll|}
\hline
$\P$ & $P$  & $G_P$  & $(n,k)$\\
\hline\hline 
$\P_{\cap}$          & $K_{n-k+1}$  & $(k-1)$-set stabiliser & \\
$\P_\cup$            & $K_{k+1}$    & $(k+1)$-set stabiliser & $n\neq 2k+2$\\
$\P_\ominus$         & ${n-2\choose k-1}K_2$ & $2$-set stabiliser& $(n,k)\neq(4,2)$\\
$\P_\cup$ and Construction \ref{con:primitive} & $2K_{k+1}$ & $S_{k+1}\Wr S_2$ & $n=2k+2$\\
Example \ref{eg:Sn}(1) & $C_4$  &  $D_8$ & $(n,k)=(4,2)$\\
Example \ref{eg:Sn}(2) & $6K_2$  &  $S_2\Wr S_3$ & $(n,k)=(6,3)$\\
\hline
\end{tabular}
\end{center}
\end{table}

\begin{theorem}\label{G=Sn}
If $(J(n,k),\P)$ is an $S_n$-primitive decomposition with
$n\geq 2k$ then $\P$ is given by one of the rows of Table
$\ref{tab:Sn}$.
\end{theorem}
\proof
Let $\Gamma=J(n,k)$, $X=\{1,\ldots, n\}$, and let $A=\{1,2,\ldots,k\}$ and
$B=\{2,\ldots,k+1\}$ be adjacent vertices of $\Gamma$. Then 
$G_{\{A,B\}}=$\mbox{$\Sym(\{1,k+1\})\times$}\\$\Sym(\{2,\ldots,k\})\times
\Sym(\{k+2,\ldots,n\})$. By Lemma \ref{lem:general}, to find all
$G$-primitive decompositions of $\Gamma$, we need to determine all
maximal subgroups $H$ of $G$ which contain $G_{\{A,B\}}$. Since $G_{\{A,B\}}$ contains a 2-cycle, \cite[Theorem 13.1]{wielandt}
implies that there are no proper primitive subgroups of $G$ containing
$G_{\{A,B\}}$.  Hence $H$ is either imprimitive or intransitive.

Suppose first that $H$ is intransitive.  Then $H$ is a maximal
intransitive subgroup and hence it has two orbits $U,W$ on $X$ and
$H=\Sym(U)\times \Sym(W)$. Since
$G_{\{A,B\}}\leqslant H$, the only possibilities for these two 
orbits are:
\[
\begin{array}{lll}
\{1,\ldots,k+1\} & \{k+2,\ldots,n\} & n\neq 2k+2\\
\{1,k+1\}        & X\backslash\{1,k+1\} & (n,k)\neq (4,2)\\
\{2,\ldots,k\}   & \{1,k+1,k+2,\ldots, n\}& \\
\end{array}
\]

When $H=\Sym(\{1,\ldots,k+1\})\times \Sym(\{k+2,\ldots,n\})=G_{A\cup B}$,
we obtain the decomposition $(\Gamma,\P_{\cup})$, while 
$H=\Sym(\{1,k+1\})\times\Sym(X\backslash\{1,k+1\})=G_{A\ominus B}$
yields the decomposition $(\Gamma,\P_\ominus)$.  Finally,
$H=\Sym(\{2,\ldots,k\})\times\Sym(\{1,k+1,k+2,\ldots,n\})=G_{A\cap B}$
gives us the decomposition $(\Gamma,\P_{\cap})$.

If $H$ is transitive but imprimitive, then the possible
systems of imprimitivity are:
\[
\begin{array}{ll}
\{1,\ldots,k+1\},\{k+2,\ldots,2k+2\}  &\text{when } n=2k+2 \\
\{1,4\},\{2,3\},\{5,6\} & \text{when } (n,k)=(6,3) \\
\{1,3\},\{2,4\}  & \text{when } (n,k)=(4,2)
\end{array}
\]

In the first case, $P=\{A,B\}^H$ is the union of two cliques each of size
$k+1$, and has as vertices all $k$-subsets of $\{1,\ldots,k+1\}$ and all
$k$-subsets of \mbox{$\{k+2,\ldots,2k+2\}$}, that is we get the decomposition
obtained from applying Construction \ref{con:primitive} to $\P_\cup$.  The last two cases give us the two decompositions from Example \ref{eg:Sn}.
\Qed

By Theorem \ref{Sn}, $A_n$ is arc-transitive on $J(n,k)$ if and only
if $n\geq 5$. Moreover, all the $S_n$-primitive decompositions in
Table \ref{tab:Sn} are $A_n$-primitive decompositions. We have the
following extra examples for alternating groups.
\begin{example}
\label{eg:An}
{\rm \begin{enumerate}
\item Let $(n,k)=(5,2)$, $G=A_5$, $A=\{1,2\}$ and $B=\{2,3\}$. Then
  $G_{\{A,B\}}=\la (1,3)(4,5)\ra$ and is contained in the maximal
  subgroup $H=\la (1,2,3,4,5),(1,3)(4,5)\ra\cong D_{10}$ of
  $G$. Letting $P=\{A,B\}^H$ and $\P=\{P^g\mid g\in G\}$,  Lemma
  \ref{lem:general} implies that  $(J(5,2),\P)$ is an $A_5$-primitive
  decomposition. Since $H_A\cong C_2$ it follows that the divisors are cycles 
of length $5$.

\item Let $(n,k)=(6,3)$, $G=A_6$, $A=\{1,2,3\}$ and
  $B=\{2,3,4\}$. Then $G_{\{A,B\}}=\la (2,3)(5,6),(1,4)(5,6)\ra$ and
  is contained in the maximal subgroup
  $H=\la(2,3)(5,6),(1,4,5)(2,3,6)\ra\cong \PSL(2,5)$ of $G$.  Letting
  $P=\{A,B\}^H$ and $\P=\{P^g\mid g\in G\}$,  Lemma \ref{lem:general}
  implies that  $(J(6,3),\P)$ is an $A_6$-primitive decomposition.
  Now $P$ is a graph on $10$ vertices with valency $3$ admitting an arc-transitive action of $H\cong A_5$. Hence $P$ is the Petersen graph.
\end{enumerate}
}%
\end{example}

\begin{lemma}
Let $\P$ be the decomposition of $J(6,3)$ given by Example
$\ref{eg:An}(2)$. Then $\P$ is $G$-primitive if and only if $G=A_6$ or
$\la A_6, (1,2)\tau\ra$ where $\tau$ is the complementation map as in Theorem \ref{thm:autgamma}.
\end{lemma}
\proof
As in the example, we take $A=\{1,2,3\}$, $B=\{2,3,4\}$ and  $P=\{A,B\}^H$ for $H=\la(2,3)(5,6),(1,4,5)(2,3,6)\ra\cong A_5$.

If $G\leq S_6$, by Theorem \ref{Sn}, $G$ must be 4-transitive, so $A_6\leq G$. We have seen above that $\P$ is $A_6$-primitive. However, $S_6$ does not preserve the partition $\mathcal{P}$ of Example \ref{eg:An}(2), since $(1,4)$ preserves $\{A,B\}$ but not $P$.
So assume  $G\not\leq S_6$. By Theorems \ref{thm:2karctrans} and
\ref{thm:2kreducetoSn}, $\P$ is a $(G\cap S_6)$-primitive
decomposition. Thus  $G\cap S_6=A_6$ and so $G=G_1=\la
A_6,\tau\ra$ or $G=G_2=\la A_6, (1,2)\tau\ra$.  Thus $|G|=2|A_6|$ and
so $|G_P:H|=2$. Then as $G_{\{A,B\}}\leqslant G_P$ it follows that
$G_{\{A,B\}}$ normalises $H$.  However, $(2,5)(3,6)\tau\in (G_1)_{\{A,B\}}$ and
does not normalise $H$, so $G\neq G_1$.  Now
$(G_2)_{\{A,B\}}=\la (1,4)(2,5)(3,6)\tau, H_{\{A,B\}}\ra$ does normalise $H$ and so fixes $P$. Thus
$\la H,(1,4)(2,5)(3,6)\tau\ra=(G_2)_P\cong S_5$  which is a maximal
subgroup of $G_2\cong S_6$. Hence $\P$ is a $G_2$-primitive
decomposition.
\Qed

We now show that Example \ref{eg:An} yields the only $A_n$-primitive decompositions which are not $S_n$-primitive.

\begin{theorem}
Let $(J(n,k),\P)$ be an $A_n$-primitive decomposition such that $A_n$ is arc-transitive and $n\geq 2k$. Then $\P$ is either an $S_n$-primitive decomposition, or $(n,k)=(5,2)$ or $(6,3)$ and $\P$ is isomorphic to a decomposition given in Example $\ref{eg:An}$.
\end{theorem}
\proof
Let $\Gamma=J(n,k)$. Since $G=A_n$ is arc-transitive it follows from Theorem \ref{Sn} that $n\geq 5$. Let  $X=\{1,\ldots,n\}$, $A=\{1,\ldots,k\}$ and $B=\{2,\ldots,k+1\}$. Then $$G_{\{A,B\}}=\big(\Sym(\{1,k+1\})\times \Sym(\{2,\ldots,k\})\times \Sym(\{k+2,\ldots,n\})\big)\cap A_n.$$ We need to consider all maximal subgroups $H$ such that $G_{\{A,B\}}<H<G$. For each such $H$,  $P=\{A,B\}^H$ is the edge-set of a divisor of the $G$-primitive decomposition.

Suppose first that $H$ is intransitive on $X$. Then
$G_{\{A,B\}}$ has the same orbits on $X$ as $(S_n)_{\{A,B\}}$ and so
$H$ is the intersection with $A_n$ of one of the maximal intransitive
subgroups which we considered in the $S_n$ case in the proof of Theorem \ref{G=Sn}. Moreover, we obtain
the decompositions in rows 1--3 in Table \ref{tab:Sn}, and so $(\Gamma,\P)$ is $S_n$-primitive.

Next suppose that $H$ is imprimitive on $X$. Since $G_{\{A,B\}}$ is
primitive on both $A\cap B$ and $\overline{A \cup B}$, the only
systems of imprimitivity preserved by $G_{\{A,B\}}$ are those
discussed in the $S_n$ case. Thus $H$ is the intersection with
$A_n$ of one of the maximal imprimitive subgroups considered
in the $S_n$ case and we obtain the decompositions in rows 4 and 6 in Table \ref{tab:Sn}. Thus $(\Gamma,\P)$ is $S_n$-primitive.

Finally, suppose that $H$ is primitive on $X$. If $k-1\geq 3$ or
$n-k-1\geq 3$, the edge stabiliser $G_{\{A,B\}}$, and hence $H$, contains a 3-cycle. Hence by \cite[Theorem 13.3]{wielandt}, $H=A_n$, contradicting $H$ being a proper subgroup. Note that if $k\geq 4$ then $k-1\geq 3$, and so $(n,k)$ is one of $(5,2)$ or $(6,3)$.

If $(n,k)=(5,2)$ then $G_{\{A,B\}}=\la(1,3)(4,5)\ra$ and $H\cong D_{10}$.
Since $A_5$ contains 15 involutions, $D_{10}$ contains 5 involutions
and there are six subgroups $D_{10}$ in $A_5$, it follows that there are 2
choices for $H$ and these are
$$H_1=\la(1, 2, 3, 4, 5),(1, 3)(4, 5)\ra$$
$$H_2=\la(1, 4, 5, 3, 2),(1, 3)(4, 5)\ra.$$
Note that $H_2=H_1^{(1,3)}$ and $(1,3)\in (S_n)_{\{A,B\}}$ and so by Lemma
\ref{lem:iso} the two decompositions obtained are
isomorphic. Moreover, $H_1$ is the stabiliser of the divisor containing
$\{A,B\}$ in the decomposition in Example \ref{eg:An}(1).

If $(n,k)=(6,3)$ then  $G_{\{A,B\}}=\la (2,3)(5,6), (1,4)(5,6)\ra$  and 
$H\cong \PSL(2,5)$. A computation using {\sc Magma} \cite{magma} showed that there are two choices for $H$ containing $G_{\{A,B\}}$ and these are:
$$H_1=\la (2, 3)(5, 6), (1, 4, 5)(2, 3, 6)\ra$$
$$H_2=\la (2, 3)(5, 6), (1, 4, 5)(3, 2, 6)\ra.$$
Note that $H_2=H_1^{(2,3)}$ and $(2,3)\in (S_n)_{\{A,B\}}$ and so the two decompositions obtained are isomorphic.  Moreover, $H_1$ is the stabiliser of the divisor containing $\{A,B\}$ in the decomposition in Example \ref{eg:An}(2).
\Qed

 We now look at the case where $n=2k$ and $G$ is not a subgroup of
$S_n$.
\begin{example}\label{eg:A4tau}
{\rm Let  $(n,k)=(4,2)$ and $G=\la A_4,(1,2)\tau\ra$. Let $A=\{1,2\}$ and
$B=\{2,3\}$. Then $G_{\{A,B\}}=\la (2,4)\tau\ra$.
\begin{enumerate}
\item Let $H_1=\la (1,2,4),(1,2)\tau\ra$ and
$$P=\{A,B\}^{H_1}=\bigg\{\big\{\{1,2\},\{2,3\}\big\},\big\{\{2,4\},\{3,4\}\big\},\big\{\{1,4\},\{1,3\}\big\}
\bigg\}.$$ Since $G_{\{A,B\}}\leqslant H_1$, it follows from Lemma
\ref{lem:general} that $(J(4,2),P^G)$ is a $G$-primitive
decomposition, with divisors isomorphic to $3K_2$.
\item  Let $H_2=\la(1,2)(3,4),(1,3)(2,4), (1,3)\tau\ra$ and
 $P=\{A,B\}^{H_2}=$
$$\bigg\{\big\{\{1,2\},\{2,3\}\big\},
\big\{\{2,3\},\{3,4\}\big\}\big\{\{3,4\},\{1,4\}\big\},
\big\{\{1,4\},\{1,2\}\big\} \bigg\}.$$ Since $G_{\{A,B\}}\leqslant H_2$, it follows from Lemma \ref{lem:general} that $(J(4,2),P^G)$ is a $G$-primitive decomposition, with divisors isomorphic to $C_4$.  Notice that this decomposition is the one in Example \ref{eg:Sn}(1) and so is also $S_4$-primitive.
 \end{enumerate}
}%
\end{example}

\begin{theorem} Let $\Gamma=J(n,k)$ with $n=2k$ and let
\mbox{$G\leqslant \Aut(\Gamma)=S_n\times S_2$} such that $G$ is not contained in $S_n$. Further, suppose that $(\Gamma,\P)$ is a $G$-primitive decomposition which is not $(G\cap S_n)$-primitive. Then $n=4$ and $\P$ is isomorphic to a decomposition given by Example $\ref{eg:A4tau}$.
\end{theorem}
\proof
By Theorems \ref{thm:2karctrans} and \ref{thm:2kreducetoSn}, it
follows that $k=2$ and $G=\la A_4,(1,2)\tau\ra$, where $\tau$ is
complementation in $X$. Let $A=\{1,2\}$ and $B=\{2,3\}$. Then
$G_{\{A,B\}}=\la (2,4)\tau\ra$.
It is not hard to see that the only maximal subgroups of $G$
containing $G_{\{A,B\}}$ are the groups $H_1$ and $H_2$ from Example
\ref{eg:A4tau}, and $H_3=\la(2,3,4),(2,3)\tau\ra$. The first two
give the two decompositions from Example \ref{eg:A4tau}.
Note that $(1,3)$ stabilizes $\{A,B\}$ and normalises $G$, and
$H_3=H_1^{(1,3)}$. So by Lemma \ref{lem:iso}, this yields a
decomposition isomorphic to the one in Example \ref{eg:A4tau}(1).
\Qed

\section{The case $k\geq 4$}

By Theorem \ref{Sn}, if $k\geq 4$ then  $G\leqslant S_n$ is arc-transitive on $J(n,k)$ if
and only if $G$ is $(k+1)$-transitive on the $n$-set $X$. Hence by the
Classificaton of finite 2-transitive permutation groups, other than $A_n$ or $S_n$,  the only possibilities for $(n,G)$ when $k\geq 4$ are $(12, M_{12})$ and $(24,M_{24})$ with $k=4$.

First we state the following well known lemmas.

\begin{lemma}\label{hexads-12}
Let $(X,\mathcal{B})$ be the Witt design $S(5,6,12)$. Then  $\mathcal{B}$ contains $132$ elements, called \emph{hexads}. Each point of $X$ is contained in $66$ hexads, each $2$-subset in $30$ hexads, each $3$-subset in $12$ hexads, each $4$-subset in $4$ hexads, and each $5$-subset in a unique hexad.
\end{lemma}
\proof
The number of hexads is given in \cite[p 31]{atlas} and then the number of hexads containing a given $i$-suset is calculated by counting $i$-subset--hexad pairs in two different ways.
\Qed
\begin{lemma}{\rm \cite[Lemma 2.11.7]{sashasbook}}
\label{hexad12stab}
Suppose that $(X,\mathcal{B})$ is a Witt design $S(5,6,12)$ preserved
by $G=M_{12}$ and let $h\in\mathcal{B}$ be a hexad. Then $G_h\cong
S_6$ and the actions of $G_h$ on $h$ and $X\backslash h$ are the two
inequivalent actions of $S_6$ on six points.
\end{lemma}

Since the stabiliser of a 3-set or a 2-set is maximal in $G=M_{12}$,
it follows from Lemma \ref{lem:classic} that $\P_\cap$ and
$\P_\ominus$ are $G$-primitive decompositions. Moreover, as $G$ acts
primitively on the point set $X$ of the Witt design, Construction
\ref{con:design2} yields a $G$-primitive decomposition of
$J(12,4)$. We also obtain a $G$-primitive decomposition from
Construction \ref{con:design3} as $G$ acts primitively on
3-subsets and one from Construction \ref{con:design4} as $G$ acts
primitively on 4-subsets. The $G$-transitive decomposition obtained
from Construction \ref{con:design} is not primitive as the stabiliser
of a hexad is
contained in the stabiliser of a pair of complementary
hexads. However, applying  Construction \ref{con:primitive} we obtain
a $G$-primitive decomposition with divisors isomorphic to $2J(6,4)$.

Before giving several more constructions arising from the Witt design, we need the following definition and lemma.
\begin{definition}
\label{def:linked3}
{\rm A \emph{linked three} in $S(5,6,12)$ is a set of four triads (or 3-sets)
such that the union of any two is a hexad.
}%
\end{definition}
\begin{lemma}\label{lem:LTM12} Let $A$, $B$ be two triads whose union is a hexad. Then there exists a unique linked three containing both $A$ and $B$.
\end{lemma}
\proof
By Lemma \ref{hexads-12}, there are exactly 12 hexads containing $A$. If such a hexad contains at least two points of $B$, then it is $A\cup B$. Let $b\in B$. Then there are 4 hexads containing $A$ and $b$, and so exactly 3 hexads meet $A\cup B$ in $A\cup \{b\}$. Therefore there are 9 hexads containing $A$ and meeting $A\cup B$ in a 4-set. Hence only two hexads contain $A$ and are disjoint from $B$. These yield two triads, $C$ and $D$, forming hexads with $A$. By Lemma \ref{hexad12stab}, the stabiliser of $A$ and $B$ is $S_3\times S_3$ which acts transitively on the remaining 6 points. Hence $C$ and $D$ must be disjoint. Since the complement of a hexad is a hexad, $C$ and $D$ must form hexads with $B$ too. It follows that $\{A,B,C,D\}$ is the unique linked three containing $A$ and $B$.
\Qed

\begin{construction}
\label{con:124wittdesign}
{\rm Let $(X,\mathcal{B})$ be the Witt design $S(5,6,12)$ and let $G=M_{12}$.
\begin{enumerate}
\item Let $T$ be a linked three as in Definition \ref{def:linked3}. Let 
$$P_T=\Big\{\big\{\{u\}\cup Y,\{v\}\cup Y\big\}\mid Y\in T, \{u,v\} \text{ contained in some triad of }T\setminus Y\Big\}$$ 
and $\P=\{P_T\mid T \text{ is a linked three}\}$.
Then $P_T\cong 12 K_3$, with each triad contributing $3K_3$.
If $\{A,B\}$ is an edge of $J(12,4)$ then $A\cup B$ is contained in a
unique hexad $A\cup B\cup\{x\}$ for some $x\in X$, and by Lemma \ref{lem:LTM12}, $\{A\cap B, \{x\}\cup (A\ominus B)\}$ is contained in a unique linked three $T$. For this $T$, $P_T$ is the unique part of $\P$ containing $\{A,B\}$. Since $G$ acts transitively on the set of linked threes and the stabiliser of a linked three is maximal, $(J(12,4),\P)$ is a $G$-primitive decomposition.

\item Let $T$ be a linked three.
A 4-set intersecting each triad of $T$ in a single point and such that its union with any triad is a hexad is called {\em special} for $T$.  For fixed triads $T_1,T_2$ of $T$ and points $x_1\in T_1$, $x_2\in T_2$, these conditions imply that there is at most one special 4-set containing $\{x_1,x_2\}$ and existence of such a 4-set was confirmed by {\sc Magma} \cite{magma}. Thus there are nine special 4-sets for $T$. Let 
$$P_T=\Big\{\big\{\{u,x,y,z\},\{v,x,y,z\}\big\}\mid \{x,y,z,t\}
 \text{sp. 4-set for }T, \{u,v,t\}\in T \Big\}$$ 
and $\P=\{P_T\mid T \text{ is a linked three}\}$.  
Then $P_T\cong 36 K_2$, with each special 4-set contributing $4K_2$. 
If $\{A,B\}$ is an edge of $J(12,4)$ then $A\cup B$ is contained in a
  unique hexad $A\cup B\cup\{x\}$ for some $x\in X$, and there is a
  unique linked three $T$ such that $(A\cap B) \cup \{x\}$ is special
  for $T$ and $\{x\}\cup (A\ominus B)$ is a triad 
of  $T$ (a {\sc Magma} \cite{magma} calculation).  Thus $P_T$ is the only part of $\P$ containing
  $\{A,B\}$. Since $G$ acts transitively on the set of linked threes
  and the stabiliser of a linked three is maximal, $(J(12,4),\P)$ is a
  $G$-primitive decomposition. 
\end{enumerate}
}%
\end{construction}

\begin{construction}
\label{con:M11graph}
{\rm Let $G=M_{12}<S_{12}$ and let $H=M_{11}$ be a $3$-transitive subgroup of
$G$. Then $H$ has an orbit of length $165$ on $4$-subsets and this orbit
forms a $3-(12,4,3)$ design. Let $\Sigma$ be the subgraph of $J(12,4)$
induced on the orbit of length $165$. The graph $\Sigma$ was studied in
\cite{DGLP}, where it is seen that  $\Sigma$ has valency $8$, is $H$-arc-transitive and given an edge
$\{A,B\}$ we have $H_{\{A,B\}}\cong S_2\times S_3=G_{\{A,B\}}$. Thus
Lemma \ref{lem:general} and the fact that $H$ is maximal in $G$,
imply that $\P=\Sigma^G$ is a $G$-primitive decomposition of $J(12,4)$.  
}%
\end{construction}

We have now seen all the $M_{12}$-primitive decompositions listed in
Table \ref{tab:M124}. It remains to prove that these are the only
ones.

\begin{table}
\begin{center}
\caption{$M_{12}$-primitive decompositions of $J(12,4)$}
\label{tab:M124}
\begin{tabular}{|lcc|}
\hline
$\P$ & $P$  & $G_P$ \\
\hline\hline 
$\P_{\cap}$          & $K_9$  & $M_9\rtimes S_3$ \\
$\P_\ominus$         & ${10\choose 3}K_2$ & $M_{10}.2$ \\
Constructions \ref{con:design} and \ref{con:primitive} & $2J(6,4)$ &  $M_{10}.2$\\
Construction \ref{con:design2} & $66K_5$ & $M_{11}$\\
Construction \ref{con:design3} & $12K_3$  & $M_9\rtimes S_3$\\ 
Construction \ref{con:design4} & $16K_2$  & $M_8\rtimes S_4$\\
Construction \ref{con:124wittdesign}(1) & $12K_3$  & $M_9\rtimes S_3$ \\
Construction \ref{con:124wittdesign}(2) &  $36K_2$ & $M_9\rtimes S_3$ \\
Construction \ref{con:M11graph} & $\Sigma$  &  $M_{11}$\\
\hline
\end{tabular}
\end{center}
\end{table}

\begin{proposition} 
If $(J(12,4),\P)$ is an $M_{12}$-primitive decomposition then $\P$ is given by one of the rows of Table $\ref{tab:M124}$.
\end{proposition}
\proof
Let $\Gamma=J(12,4)$ and $G=M_{12}$ acting on the point set $X$ of the Witt-design $S(5,6,12)$. Take adjacent vertices
$A=\{1,2,3,4\}$ and $B=\{2,3,4,5\}$ and suppose that
$h=\{1,2,3,4,5,6\}$ is the unique hexad containing $A\cup B$. Then 
$G_{\{A,B\}}=G_{\{1,5\},\{2,3,4\},\{6\}}\cong S_2\times S_3$, by Lemma
\ref{hexad12stab}. Since
transpositions in the action of $G_h$ on $h$ act as a product of
three transpositions on $X\backslash h$, and 3-cycles on $h$ act as a
product of two 3-cycles on $X\backslash h$ it follows that
$G_{1,5,6,\{2,3,4\}}\cong S_3$ acts regularly on $X\backslash h$, and
so $G_{\{A,B\}}$ acts transitively on $X\backslash h$.

Let $H$ be a maximal subgroup of $G$ such that 
$G_{\{A,B\}}\leqslant H<G$. The maximal subgroups of $G$ are given in
\cite[p 33]{atlas}.  The orbit lengths of $G_{\{A,B\}}$ imply that
$G_{\{A,B\}}$ does not preserve a system of imprimitivity on $X$ with blocks
of size 2 or 4 and so $H\not\cong C_4^2\rtimes D_{12}, A_4\times S_3$,
or $C_2\times S_5$. Moreover, $|H_6|$ is even and so 
$H\not\cong \PSL(2,11)$.

If $H$ is intransitive then $H$ is one of $G_{\{2,3,4,6\}}$,
$G_{\{2,3,4\}}$, $G_{\{1,5,6\}}$, $G_{\{1,5\}}$ or $G_6$. (Note that $G_h$ is not maximal.) The first is
the stabiliser of the divisor containing $\{A,B\}$ in the decomposition
yielded by Construction \ref{con:design4}. The second gives
$\P_\cap$ while the third is the stabiliser of the divisor of the
decomposition yielded by Construction \ref{con:design3}
containing $\{A,B\}$. If $H=G_{\{1,5\}}$ then we obtain the
decomposition $\mathcal{P}_{\ominus}$ while if $H=G_6$ we obtain the
decomposition yielded by Construction \ref{con:design2}. 

The only hexad pair fixed by $G_{\{A,B\}}$ is $\{h,X\backslash
h\}$. Now $G_h$ is the stabiliser of the divisor of the decomposition
yielded by Construction \ref{con:design} containing
$G_{\{A,B\}}$. Such a divisor is isomorphic to $J(6,4)$ and so
$G_{\{h,X\backslash h\}}$ yields the decomposition with divisors isomorphic
to $2J(6,4)$ obtained after applying Construction \ref{con:primitive}.

A calculation using {\sc Magma} \cite{magma}  shows that there is only one
transitive subgroup of $G$ isomorphic to $M_{11}$ which contains $G_{\{A,B\}}$ and this yields Construction \ref{con:M11graph}.

By the list of maximal subgroups of $G$ given in \cite[p 33]{atlas}, the only case left to consider is  $H$ being the stabiliser of a
linked three.
If  $T$ is a linked three preserved by $G_{\{A,B\}}$ then
$\{1,5,6\}$ is a triad of $T$ and either $\{2,3,4\}$ is also a triad
or $2$, $3$, and $4$ lie in distinct triads. Since a linked three is
uniquely determined by any two of its triads (Lemma \ref{lem:LTM12}), there is a
unique linked three $T$ containing  $\{1,5,6\}$ and 
$\{2,3,4\}$.  Then $G_T$ is the stabiliser of the divisor of the
decomposition yielded by Construction \ref{con:124wittdesign}(1)
containing $\{A,B\}$. If $2$, $3$ and $4$ are in distinct blocks, a calculation using {\sc Magma} \cite{magma} shows that there is a unique $H$ containing
$G_{\{A,B\}}$  and we obtain the decomposition in Construction
\ref{con:124wittdesign}(2).
\Qed

We need the following well known lemma to deal with the case where $G=M_{24}$.
\begin{lemma}\label{octads}
\cite[Lemma 2.10.1]{sashasbook}
Let $(X,\mathcal{B})$ be the Witt design $S(5,8,24)$.  Then
$\mathcal{B}$ contains $759$ elements, called \emph{octads}. Each point
of $X$ is contained in $253$ octads, each $2$-subset in $77$ octads, each
$3$-subset in $21$ octads, each $4$-subset in $5$ octads, and each $5$-subset in
a unique octad. Moreover, the stabiliser of an octad in $M_{24}$ is
$C_2^4\rtimes A_8$ where $C_2^4$ acts trivially on the octad and
transitively on its complement.
\end{lemma}
\proof
Then number of octads comes from \cite[Lemma 2.10.1]{sashasbook} and then the numbers of octads containing a given $i$-subset follows from basic counting. The statement about the stabiliser of an octad also comes from \cite[Lemma 2.10.1]{sashasbook}.
\Qed

Since the stabilisers of a 3-set, of a 2-set, and of an octad are maximal in $G$, applying Constructions \ref{con:classic}, \ref{con:design} and
\ref{con:design2}, we get the list of $M_{24}$-primitive decompositions
in Table \ref{tab:M244}.

\begin{table}
\begin{center}
\caption{$M_{24}$-primitive decompositions of $J(24,4)$}
\label{tab:M244}
\begin{tabular}{|lcc|}
\hline
$\P$ & $P$  & $G_P$ \\
\hline\hline 
$\P_{\cap}$          & $K_{21}$  & $\PGammaL(3,4)$ \\
$\P_\ominus$         & ${22\choose 3}K_2$ & $M_{22}.2$ \\
Construction \ref{con:design} & $J(8,4)$ &  $C_2^4\rtimes A_8$\\
Construction \ref{con:design2} & $21K_5$ & $\PGammaL(3,4)$\\
\hline
\end{tabular}
\end{center}
\end{table}

\begin{proposition}
If $(J(24,4),\P)$ is an $M_{24}$-primitive decomposition then $\P$ is given by one of the rows in Table $\ref{tab:M244}$.
\end{proposition}
\proof
Let $\Gamma=J(24,4)$ and $G=M_{24}$ acting on the point-set $X$ of the Witt-design $S(5,8,24)$. Take adjacent vertices $A=\{1,2,3,4\}$ and $B=\{2,3,4,5\}$ and suppose that $\Delta=\{1,2,3,4,5,6,7,8\}$ is the unique octad containing $A\cup B$. Then looking at the stabiliser of an octad given in Lemma \ref{octads}, we see that
$G_{\{A,B\}}=G_{\{1,5\},\{2,3,4\},\{6,7,8\}}= C_2^4\rtimes ((S_2\times S_3^2)\cap A_8)$ with orbits in $\Delta$ of lengths 2, 3, 3. Since $G_{\{A,B\}}$ contains the pointwise stabiliser of the octad $\Delta$, which by Lemma \ref{octads} acts regularly on $X\setminus\Delta$, it follows that $G_{\{A,B\}}$  is transitive on $X\setminus \Delta$.

Let $H$ be a maximal subgroup of $G$ such that
$G_{\{A,B\}}\leqslant H<G$. The maximal subgroups of $G$ are given in
\cite[p 96]{atlas}, and comparing orders we see that $H\not\cong \PSL(2,7)$
or $\PSL(2,23)$. Since $G_{\{A,B\}}$ has an orbit of length 16 and an
orbit of length 3 in $X$, it cannot fix a pair of dodecads. Similarly, if $H$
fixed a trio of disjoint octads, one of the three octads would be
$\Delta$ and $G_{\{A,B\}}$ would interchange the other
2. However, all index 2 subgroups of $G_{\{A,B\}}$ are transitive on $X\setminus \Delta$ (a {\sc Magma} calculation \cite{magma}) and so $H$
does not fix a trio of disjoint octads.  Suppose next that $H$ fixes a
sextet, that is, 6 sets of size 4 such that the union of any two is an
octad. Then the $G_{\{A,B\}}$-orbit $X\setminus \Delta$ is the union of four
of these sets. However, the remaining $G_{\{A,B\}}$-orbit lengths
are incompatible with $H$
fixing a partition of $\{1,\ldots,8\}$ into two sets of size 4. Thus
the list of maximal subgroups of $G$ in \cite[p 96]{atlas} implies that $H$ is
intransitive on $X$, and so $H=G_{\{1,5\}}, G_{\{2,3,4\}}, G_{\{6,7,8\}}$, or
$G_{\{1,2,3,4,5,6,7,8\}}$.  By Lemma \ref{lem:classic}, the first gives the
decomposition $\mathcal{P}_{\ominus}$ while the second
gives $\mathcal{P}_{\cap}$. The third is the stabiliser of the divisor of
the decomposition yielded by Construction \ref{con:design2} containing
$\{A,B\}$ while the fourth yields  the decomposition obtained from
Construction \ref{con:design}.
\Qed

\section{The case $k=3$}

By Theorem  \ref{Sn}, $G\leqslant S_n$ is arc-transitive on
$J(n,3)$ if and only if $G$ is 4-transitive or $G=\PGammaL(2,8)$ and
$n=9$. Thus other than $A_n$ or $S_n$ the only possibilites for
$(n,G)$ are  $(11,M_{11}),(12,M_{12})$, $(23,M_{23})$, $(24,M_{24})$
and $(9,\PGammaL(2,8))$.

Since the stabiliser of a 2-subset is maximal in $M_{24}$, it follows
that $\P_\cap$ and $\P_\ominus$ are $M_{24}$-primitive
decompositions with divisors  $K_{22}$ and ${22\choose 2}K_2$
respectively. We also have a construction involving sextets. 

\begin{construction}
\label{con:sextets}
{\rm Let $S$ be a sextet, that is, a set of six $4$-subsets such that the
union of any two is an octad, and define 
$P_S=\{\{A,B\}\mid A\cup B\in S\}$
and $\P=\{P_S\mid S\text{ a sextet}\}$.
Then $P_S\cong 6J(4,3)\cong 6K_4$ with one copy of $K_4$ for each
$4$-set in $S$. Let $\{A,B\}$ be an edge of $J(24,3)$. By 
\cite[Lemma 2.3.3]{sashasbook}, $A\cup B$ is a member of a unique
sextet $S$ and so $P_S$ is the only part of $\P$ containing
$\{A,B\}$. Since $G$ acts primitively on the set of sextets, it
follows that $(J(24,3),\P)$ is an $M_{24}$-primitive decomposition.
}%
\end{construction}

\begin{proposition}  If $(J(24,3),\P)$ is an $M_{24}$-primitive
  decomposition then either $\P=\P_\ominus$ or $\P_\cap$, or $\P$
  arises from Construction $\ref{con:sextets}$.
\end{proposition}
\proof
Let $\Gamma=J(24,3)$ and $G=M_{24}$ acting on the point set $X$ of the Witt-design $S(5,8,24)$. Let $A=\{1,2,3\}$ and $B=\{2,3,4\}$ be adjacent vertices in $\Gamma$. Then $G_{\{A,B\}}=G_{\{1,4\},\{2,3\}}$ which is the stabiliser in $\Aut(M_{22})$ of a 2-subset and so by \cite[p 39]{atlas},
$G_{\{A,B\}}\cong 2^5\rtimes S_5$. Since $G$ is 5-transitive on $X$, $G_{\{A,B\}}$ is transitive on $X\backslash \{1,2,3,4\}$.

Let $H$ be a  maximal subgroup of $G$ such
that $G_{\{A,B\}}\leqslant H<G$. The maximal subgroups of $G$ can be
found in \cite{atlas}. Comparing orders we see that 
$H\not\cong \PSL(2,7)$, $\PSL(2,23)$, or the
stabiliser of a trio of distinct octads. Now $G_{\{A,B\}}$ contains
$G_{1,2,3,4}$ which is transitive on the remaining 20 points. Thus
$G_{1,2,3,4}$ does not fix a pair of dodecads and so neither does
$H$. Hence by the list of maximal subgroups of $G$ in \cite[p 96]{atlas}, either $H$ is
intransitive, or fixes a sextet. If $H$ is intransitive, then
$H=G_{\{1,4\}}$ or $G_{\{2,3\}}$. By Lemma \ref{lem:classic}, the
first gives $\mathcal{P}_{\ominus}$ while the second gives
$\mathcal{P}_{\cap}$.

Suppose then that $H$ fixes a sextet.  The orbit lengths of
$G_{\{A,B\}}$ imply that $\{1,2,3,4\}$ is one of the 
blocks of the sextet. By \cite[Lemma 2.3.3]{sashasbook}, 
$\{1,2,3,4\}$ is contained in a unique sextet $S$.   Thus $H=G_S$ and
is the stabiliser in $G$ of the divisor of the decomposition obtained
from Construction \ref{con:sextets} containing $\{A,B\}$.
\Qed

Before dealing with $G=M_{23}$ we need the following well known result
which follows from Lemma \ref{octads}.

\begin{lemma}\label{heptads}
Let $(X,\mathcal{B})$ be the Witt design $S(4,7,23)$. Then  $\mathcal{B}$ contains $253$ elements, called \emph{heptads}. Each point of $X$ is contained in $77$ heptads,
each $2$-subset in $21$ heptads, each $3$-subset in $5$ heptads, and each
$4$-subset in a unique heptad. Moreover, the stabiliser of a heptad is
$C_2^4\rtimes A_7$ with the pointwise stabiliser of the heptad being
$C_2^4$ which acts regularly on the $16$ points not in the heptad.
\end{lemma}
\proof
Since $(X,\mathcal{B})$ is derived from the set of all blocks of the Witt design $S(5,8,24)$ containing a given point, this follows from Lemma \ref{octads}.
\Qed

Using the Witt design $S(4,7,23)$ and the fact that the stabiliser of
a 2-set is maximal in $M_{23}$ we get the $M_{23}$-primitive
decompositions in Table \ref{tab:M233}. These are in fact all such
decompositions.
\begin{table}
\begin{center}
\caption{$M_{23}$-primitive decompositions of $J(23,3)$}
\label{tab:M233}
\begin{tabular}{|lcc|}
\hline
$\P$ & $P$  & $G_P$ \\
\hline\hline 
$\P_{\cap}$          & $K_{21}$  & $\PSigmaL(3,4)$ \\
$\P_\ominus$         & ${21\choose 2}K_2$ & $\PSigmaL(3,4)$ \\
Construction \ref{con:design} & $J(7,3)$ &  $C_2^4\rtimes A_7$\\
Construction \ref{con:design2} & $5K_4$ & $C_2^4\rtimes(C_3\times
A_5)\rtimes C_2$\\
\hline
\end{tabular}
\end{center}
\end{table}

\begin{proposition}
If $(J(23,3),\P)$ is an $M_{23}$-primitive decomposition then $\P$ is as in one of the lines of Table $\ref{tab:M233}$.
\end{proposition}
\proof
Let $\Gamma=J(23,3)$ and $G=M_{23}$ acting on the point-set
$X$ of the Witt-design $S(4,7,23)$. Take adjacent vertices 
$A=\{1,2,3\}$ and $B=\{2,3,4\}$. By Lemma \ref{heptads}, $\{1,2,3,4\}$
is contained in a unique heptad, $h=\{1,2,3,4,5,6,7\}$ say, and so
$G_{\{A,B\}}$ fixes $h$. Since the stabiliser of a heptad is
isomorphic to $C_2^4\rtimes A_7$ (Lemma \ref{heptads}), it follows
that $G_{\{A,B\}}$ has order 192 and has orbits $\{1,4\}$, $\{2,3\}$,
$\{5,6,7\}$ and $X\backslash h$.

Let $H$ be a maximal subgroup of $G$ such
that $G_{\{A,B\}}\leqslant H<G$. The maximal subgroups of $G$ can be
found in \cite{atlas}. By comparing orders, $H\not\cong C_{23}\rtimes
C_{11}$ and so $H$ is intransitive. Thus $H=G_{\{1,4\}}, G_{\{2,3\}},
G_{\{5,6,7\}}$ or $G_h$. By Lemma \ref{lem:classic}, the first two 
give the decompositions $\mathcal{P}_{\ominus}$ and
$\mathcal{P}_{\cap}$ respectively. Also $G_{\{5,6,7\}}$ is the
stabiliser of the divisor of the decomposition obtained from Construction
\ref{con:design2} containing $\{A,B\}$ while $G_h$ is the stabiliser
of the divisor of the decomposition yielded by Construction
\ref{con:design}.
\Qed

Since 4-set stabilisers and 2-set stabilisers are maximal in $M_{12}$,
it follows from Lemma \ref{lem:classic} that $\P_{\cup}$, $\P_{\cap}$
and $\P_{\ominus}$ are $M_{12}$-primitive decompositions with divisors
isomorphic to $K_4$, $K_{10}$ and ${10 \choose 2}K_2$ respectively. We
also have the following construction.
\begin{construction}
\label{con:linkedfours}
{\rm Let $(X,\mathcal{B})$ be the Witt design $S(5,6,12)$.
Let $F$ be a \emph{linked four}, that is a set of three mutually disjoint tetrads (sets of size 4) admitting a refinement into six duads (called duads of $F$) such that the union of any three duads coming from any two tetrads is a hexad. 
Let $$P_F=\Big\{\big\{\{x,u,v\},\{y,u,v\}\big\}\mid \{x,y,u,v\}\in F,
  \{u,v\},\{x,y\} \text{ are duads of $F$} \Big\}$$
and let $\P=\{P_F\mid F \text{ a linked four}\}$.
Then $P_F\cong 6K_2$ with one copy of $2K_2$ for each tetrad in $F$. Let $\{A,B\}$ be an edge of $J(12,3)$. It turns out ({\sc Magma} calculation \cite{magma}) there is exactly one linked four $F$ having $A\cup B$ as a tetrad and $A\cap B$ as a duad of $F$, and so $P_F$ is the only part of $\P$ containing $\{A,B\}$.  Since $G$ acts primitively on the set of linked fours, it follows that $(J(12,3),\P)$ is an $M_{12}$-primitive decomposition.
}%
\end{construction}

\begin{proposition}
If $(J(12,3),\P)$ is an $M_{12}$-primitive decomposition then $\P=\P_\cup,\P_\cap$ or $\P_\ominus$ or  $\P$ is obtained from Construction $\ref{con:linkedfours}$.
\end{proposition}
\proof
Let $\Gamma=J(12,3)$ and $G=M_{12}$ acting on the point set $X$ of the Witt-design $S(5,6,12)$. Take adjacent vertices
$A=\{1,2,3\}$ and $B=\{2,3,4\}$.  The stabiliser in $G$ of a 4-set is
$M_8\rtimes S_4$ such that the pointwise stabiliser $M_8$ of the 4-set
acts regularly on the 8 remaining points. Hence
$G_{\{A,B\}}=G_{\{1,4\},\{2,3\}}=M_8\rtimes (S_2\times S_2)$ which has
order 32 and is transitive on the 8 points of
$X\backslash\{1,2,3,4\}$.

Let $H$ be a maximal subgroup of $G$ such that
$G_{\{A,B\}}\leqslant H<G$. The maximal subgroups of $G$ are given in
\cite{atlas}, and comparing orders we see that $H\not\cong M_{11}$,
$\PSL(2,11)$, $M_9\rtimes S_3$, $C_2\times S_5$ and $A_4\times S_3$. Moreover,
since $G_{\{A,B\}}$ has orbits of size 2,2 and 8 in $X$ it does not stabilise
a hexad pair. If $H$ is intransitive then $H=G_{\{1,2,3,4\}}$,
$G_{\{1,4\}}$ or $G_{\{2,3\}}$. These yield
$\P_\cup$, $\P_{\ominus}$ and $\P_{\cap}$
respectively. Thus by \cite[p 33]{atlas} we are left to consider the case where 
$H\cong 4^2\rtimes D_{12}$.  A {\sc Magma} \cite{magma} calculation shows
that there is a unique such $H$ containing $G_{\{A,B\}}$ and we obtain
the decomposition from Construction \ref{con:linkedfours}.
\Qed

Before dealing with $G=M_{11}$ we need the following couple of lemmas,
the first of which is well known.

\begin{lemma}\label{pentads}
Let $(X,\mathcal{B})$ be the Witt design $S(4,5,11)$.
Then  $\mathcal{B}$ contains $66$ elements, called \emph{pentads}. Each
point of $X$ is contained in $30$ pentads, each $2$-subset in $12$ pentads,
each $3$-subset in $4$ pentads, and each $4$-subset in a unique
pentad. Moreover, the stabiliser of a pentad is isomorphic to $S_5$,
which acts in its natural action on the pentad and as $\PGL(2,5)$ on
the complementary hexad.
\end{lemma}
\proof
Since $(X,\mathcal{B})$ can be derived from the set of blocks of the  Witt design $S(5,6,12)$ containing a given point, the first part follows from Lemma \ref{hexads-12}. By \cite[p 18]{atlas}, the stabiliser of a pentad is $S_5$ and has two orbits on $X$.
\Qed

\begin{lemma}\label{lem:calcpentads}
Let $(X,\mathcal{B})$ be the Witt design $S(4,5,11)$ and
$G=M_{11}$. Let $A=\{1,2,3\}$, $B=\{2,3,4\}$ and suppose that
$p=\{1,2,3,4,5\}$ is the unique pentad containing $A\cup B$. Then
$G_{\{A,B\}}\cong C_2^2$ and on $X\backslash p$ has an orbit $\{a,b\}$
of length $2$ and an orbit of length $4$. Moreover, $\{1,4,5,a,b\}$,
$\{2,3,5,a,b\}$  and $X\backslash\{1,2,3,4,a,b\}$ are pentads.
\end{lemma}
\proof
By Lemma \ref{pentads}, $G_p$ induces $S_5$ on $p$, and since
$G_{\{A,B\}}\leqslant G_p$ it follows that
$G_{\{A,B\}}=G_{\{2,3\},\{1,4\}}\cong C_2^2$ and fixes the point $5$.
By \cite{atlas}, each involution of $G$ fixes precisely three points
of $X$. Two of the involutions of $G_{\{A,B\}}$ fix three points of
$p$ and so are fixed point free on $X\backslash p$. The third
involution fixes the point 5 and fixes two points $a,b$ of
$X\backslash p$. It follows that $G_{\{A,B\}}$ has an orbit of length
two (namely, $\{a,b\}$) and an orbit of length 4 on $X\backslash p$.

Any four points lie in a unique pentad and by Lemma \ref{pentads}, any
3-subset is contained in 4 pentads. Hence $X\backslash p$ is 
divided into three sets of size two by the three pentads containing
$\{1,4,5\}$ other than $\{1,2,3,4,5\}$. Similarly,  $X\backslash p$ is
partitioned by the three pentads containing $\{2,3,5\}$. Since
$G_{\{A,B\}}$ fixes $\{1,4,5\}$ and $\{2,3,5\}$, it 
preserves both partitions and $\{a,b\}$ must be a block of both.
Hence $\{1,4,5,a,b\}$ and $\{2,3,5,a,b\}$ are pentads.  Moreover,
since $X\backslash (\{a,b\}\cup p)$ is an orbit of length 4 of
$G_{\{A,B\}}$ and is contained in a unique pentad, the fifth point of
this pentad must also be fixed by $G_{\{A,B\}}$ and hence is 5. Thus
$X\backslash\{1,2,3,4,a,b\}$ is a pentad.
\Qed

Since the stabiliser of a 2-set is maximal in $M_{11}$, it follows
from Lemma \ref{lem:classic} that $\P_\cap$ and $\P_\ominus$ are
$M_{11}$-primitive decompositions. We also obtain $M_{11}$-primitive
decompositions from Constructions \ref{con:design},
\ref{con:design2}, \ref{con:design3} and \ref{con:design4} by using
the Witt design $S(4,5,11)$, since the stabilisers of a block, of a point and of a 3-subset are maximal subgroups of $M_{11}$.

\begin{construction}
\label{con:M113}
{\rm Let $(X,\mathcal{B})$ be the Witt design $S(4,5,11)$ and
$G=M_{11}$. Let $A=\{1,2,3\}$ and $B=\{2,3,4\}$ be adjacent vertices
of $J(11,3)$ and let $\{a,b\}$ be the orbit of length 2 of
$G_{\{A,B\}}$ on $X\backslash \{1,2,3,4,5\}$ given by Lemma
\ref{lem:calcpentads}.
\begin{enumerate}
\item For each $3$-subset $Y$ of $X$ let 
$$P_Y=\Big\{\big\{\{x,u,v\},\{y,u,v\}\big\}\mid \{x,y\}\cup Y, \{u,v\}\cup Y\in \mathcal{B} \Big\}$$ 
and let $\P=\{P_Y\mid Y \text{ a $3$-subset}\}$. By Lemma \ref{pentads}, $Y$ is contained in 4 pentads, and so $12K_2$. Let $Y=\{5,a,b\}$. By Lemma \ref{lem:calcpentads}, $\{A,B\}\in P_Y$ and $G_{\{A,B\}}\leqslant G_Y=G_{P_Y}$, which is a maximal subgroup of $G$. Hence by Lemma \ref{lem:general}, $(J(11,3),\P)$ is an $M_{11}$-primitive decomposition.

\item Since $G$ is 4-transitive on $X$, Lemma \ref{lem:calcpentads} implies that the stabiliser in $G$ of two $2$-subsets of $X$ fixes a third. For each $2$-subset $Y$ let
$$P_Y=\Big\{\big\{\{x,u,v\},\{y,u,v\}\big\}\mid u,v,x,y\in X\backslash Y, G_{Y,\{x,y\}} =G_{Y,\{u,v\}}\Big\}$$
and let $\P=\{P_Y\mid Y \text{ a 2-subset}\}$. Then each $P_Y\cong {9\choose 2}K_2$. Moreover, by Lemma \ref{lem:calcpentads} any edge of $J(11,3)$ is contained in a unique part of $\P$ ($\{A,B\}\in P_{\{a,b\}}$) and so $(J(11,3),\P)$ is an $M_{11}$-primitive decomposition.

\item  For each $Y\in\mathcal{B}$ let 
$$P_Y=\Big\{\big\{\{x,u,v\},\{y,u,v\}\big\}\mid x,y\in Y, \{u,v\}\cup(Y\backslash\{x,y\}) \in\mathcal{B}\Big\}$$
and let $\P=\{P_Y\mid Y\in\mathcal{B}\}$. By Lemma \ref{pentads}, each $3$-subset of $Y$ is contained in three more pentads and so each part of $\P$ is isomorphic to $3{5\choose 2}K_2=30K_2$. By Lemma \ref{lem:calcpentads}, $\{A,B\}\in P_Y$ for $Y=\{1,4,5,a,b\}$. Moreover, $G_{\{A,B\}}$ fixes $Y$ and so $G_{\{A,B\}}<G_Y=G_{P_Y}$. Thus Lemma \ref{lem:general} and the fact that $G$ acts primitively on $\mathcal{B}$, imply that $(J(11,3),\P)$ is a $G$-primitive decomposition.

\item  For each $Y\in\mathcal{B}$ let 
$$P_Y=\Big\{\big\{\{x,u,v\},\{y,u,v\}\big\}\mid u,v\in Y, \{x,y\}\cup(Y\backslash\{u,v\}) \in\mathcal{B}\Big\}$$
and let $\P=\{P_Y\mid Y\in\mathcal{B}\}$. By Lemma \ref{pentads}, each $3$-subset of $Y$ is contained in three more pentads and so each part of $\P$ is isomorphic to
$3{5\choose 2}K_2=30K_2$. By Lemma \ref{lem:calcpentads}, $\{A,B\}\in P_Y$ for $Y=\{2,3,5,a,b\}$ and $G_{\{A,B\}}<G_Y=G_{P_Y}$. Thus Lemma \ref{lem:general} and the fact that $G$ acts primitively on $\mathcal{B}$, imply that $(J(11,3),\P)$ is a $G$-primitive decomposition.
\end{enumerate}
}%
\end{construction}

\begin{construction}
\label{con:PSL211}
{\rm Let $H=\PSL(2,11)<M_{11}=G$. Then $H$ has an orbit of length $55$ on
$3$-subsets and this orbit forms a $2-(11,3,3)$ design known as the
Petersen design.  The remaining $3$-subsets form an orbit of length
110 and a $2-(11,3,6)$ design \cite{PSL}.
\begin{enumerate}
\item Let $\Pi$ be the subgraph of $J(11,3)$ induced on
the orbit of length 55. The graph $\Pi$ was studied in
\cite{DGLP} and is $H$-arc-transitive of valency 6. Given an edge
$\{A,B\}$ of $\Pi$ we have $H_{\{A,B\}}=C_2^2=G_{\{A,B\}}$. Thus
letting $\P=\{\Pi^g\mid g\in G\}$, it follows by Lemma
\ref{lem:general} that $(J(11,3),\P)$ is a $G$-primitive
decomposition.
\item Let $\Delta$ be the subgraph of $J(11,3)$ induced on the orbit
  of length 110. Then $\Delta$ has valency $15$ and given a vertex
  $A$, $H_A\cong S_3$ has orbits of length $3$, $6$ and $6$ on the neighbours
  of $A$. Let $B$ be a neighbour of $A$ in the orbit of length $3$ and
  let $P=\{A,B\}^H$. Let $g\in H$ which interchanges $A$ and $B$. Then
  by Lemma \ref{lem:cosetgraph}, $P\cong
  \Cos(H,H_A,H_AgH_A)$. Moreover, $\la H_A,g\ra\cong A_5$ and so $P$
  has 11 connected components, each with 10 vertices and isomorphic to
  the Petersen graph. Since $|H_{\{A,B\}}|=4=|G_{\{A,B\}}|$, it
  follows from Lemma \ref{lem:general} that $(J(11,3),\P)$ is a
  $G$-primitive decomposition with $\P=P^G$.
\end{enumerate}
}%
\end{construction}

\begin{construction}
\label{con:petersen}
{\rm  Let $A=\{1,2,3\}$ and $B=\{2,3,4\}$. By Lemma
  \ref{lem:calcpentads}, $Y=X\backslash\{1,2,3,4,a,b\}$ is a pentad
  fixed by $G_{\{A,B\}}$. Let $H=G_Y$ and $P=\{A,B\}^H$. Then by Lemma
  \ref{pentads}, $H$ induces $S_5$ on $Y$ and $\PGL(2,5)$ on
  $\{1,2,3,4,a,b\}$. Thus $H_A\cong S_3$ and is a maximal subgroup of
  $A_5\cong \PSL(2,5)$. Moreover, \mbox{$g\in H_{\{A,B\}}$} which interchanges $A$ and
  $B$ induces even permutations on $Y$ and so for such a $g$ we have
  $\la H_A,g\ra =A_5$. By Lemma \ref{lem:cosetgraph}, 
  $P\cong \Cos(H,H_A,H_AhH_A)$. Since $|H:H_A|=20$ and $\la H_A,g\ra
  \cong A_5$, it follows that $P$ has two disconnected components with
  10 vertices each. Since $|H_A:G_{A,B}|=3$ it follows that $P$ is a
  copy of two Petersen graphs. Let
  $\P=P^G$. Then as $G_{\{A,B\}}<H$, it follows from Lemma
  \ref{lem:general} that $(J(11,3),\P)$ is a $G$-primitive
  decomposition.
}%
\end{construction}

\begin{table}
\begin{center}
\caption{$M_{11}$-primitive decompositions of $J(11,3)$}
\label{tab:M113}
\begin{tabular}{|lcc|}
\hline
$\P$ & $P$  & $G_P$ \\
\hline\hline 
$\P_{\cap}$          & $K_9$  & $M_9\rtimes C_2$ \\
$\P_\ominus$         & ${9\choose 2}K_2$ & $M_9\rtimes C_2$ \\
Construction \ref{con:design} & $J(5,3)\cong J(5,2)$ &  $S_5$\\
Construction \ref{con:design2} & $30K_4$ & $M_{10}$\\
Construction \ref{con:design3} & $4K_3$ & $M_8\rtimes S_3$\\
Construction \ref{con:design4} & $12K_2$ & $M_8\rtimes S_3$\\
Construction \ref{con:M113}(1)    & $12K_2$ & $M_8\rtimes S_3$ \\
Construction \ref{con:M113}(2)    & ${9\choose 2}K_2$ & $M_9\rtimes C_2$ \\
Construction \ref{con:M113}(3)    & $30K_2$ & $S_5$ \\
Construction \ref{con:M113}(4)    & $30K_2$ & $S_5$ \\
Construction \ref{con:PSL211}(1) & $\Pi$ &  $\PSL(2,11)$ \\
Construction \ref{con:PSL211}(2) & 11 Petersen graphs &  $\PSL(2,11)$ \\
Construction \ref{con:petersen}  &2 Petersen graphs & $S_5$\\
\hline
\end{tabular}
\end{center}
\end{table}

\begin{proposition} If $(J(11,3),\P)$ is an $M_{11}$-primitive
  symmetric decomposition  then $\P$ is given by Table $\ref{tab:M113}$.
\end{proposition}
\proof
Let $\Gamma=J(11,3)$ and $G=M_{11}<\Sym(X)$, and consider $X$ as the point set of the Witt-design $S(4,5,11)$ with automorphism group $G$.
Let $A=\{1,2,3\}$ and $B=\{2,3,4\}$ be adjacent vertices. Suppose
that $p=\{1,2,3,4,5\}$ is the unique pentad of the Witt design
containing $\{1,2,3,4\}$ and let $H$ be a maximal subgroup of $G$
containing $G_{\{A,B\}}=G_{\{2,3\},\{1,4\}}$. The maximal subgroups of $G$ are given in \cite[p 18]{atlas}.

If $H$ is the stabiliser of a point then $H=G_5$ and so we obtain the
decomposition yielded by Construction \ref{con:design2}.
Next suppose that $H$ is the stabiliser of a duad. Then $H$ is one of
$G_{\{2,3\}}, G_{\{1,4\}}$ or $G_{\{a,b\}}$ where $\{a,b\}$ is the
orbit of length two of $G_{\{A,B\}}$ on $\{6,7,\ldots,11\}$. The first
gives $\P_{\cap}$ the second gives $\P_{\ominus}$.
Finally, if $H=G_{\{a,b\}}$ then $H$ is the stabiliser of the divisor of
the decomposition obtained from Construction \ref{con:M113}(2)
containing $\{A,B\}$.

Next suppose that $H$ is the stabiliser of a triad. Then $H$
stabilises $\{1,4,5\}$, $\{2,3,5\}$ or $\{5,a,b\}$. If
$H=G_{\{1,4,5\}}$ then $H$ is the stabiliser of the divisor of the
decomposition from Construction \ref{con:design3} containing 
$\{A, B\}$. Also $H=G_{\{2,3,5\}}$ is the stabiliser of the divisor of
the decomposition yielded by Construction \ref{con:design4} containing
$\{A,B\}$. Finally, $H=G_{\{5,a,b\}}$ is the stabiliser of the divisor of
the decomposition obtained from Construction \ref{con:M113}(1) containing
$\{A,B\}$.

Next suppose that $H$ is the stabiliser of a pentad. Since
$G_{\{A,B\}}$ has only one orbit  of odd length, it follows that $5$ is in the pentad. Combining $5$ with two orbits of $G_{\{A,B\}}$ of length two we
get that $G_{\{A,B\}}$ fixes the  pentads $\{1,2,3,4,5\}$,
$\{1,4,5,a,b\}$, $\{2,3,5,a,b\}$ and $X\backslash\{1,2,3,4,a,b\}$
(by Lemma \ref{lem:calcpentads}, these 5-sets are actually
pentads). Thus there are four choices for $H$. 
If $H=G_{\{1,2,3,4,5\}}$ then we obtain the decomposition from
Construction \ref{con:design}. If $H=G_{\{1,4,5,a,b\}}$, then $H$ is
the stabiliser of the divisor of the decomposition from Construction
\ref{con:M113}(3) containing $\{A,B\}$ while
$H=G_{\{2,3,5,a,b\}}$ is the stabiliser of the divisor of the
decomposition yielded by Construction \ref{con:M113}(4).
Finally, if $H=G_{X\backslash\{1,2,3,4,a,b\}}$ then $H$ is the
stabiliser of the divisor of the decomposition produced by Construction
\ref{con:petersen} containing $\{A,B\}$.

We are left to consider $H\cong \PSL(2,11)$. By a calculation using {\sc Magma}
\cite{magma}, there are two such $H$ containing $G_{\{A,B\}}$. These
give us the two decompositions in Construction \ref{con:PSL211}.
\Qed

We now give constructions for $\PGammaL(2,8)$-primitive decompositions of
$J(9,3)$. 

\begin{construction}
\label{con:pgammal}
{\rm Let $G=\PGammaL(2,8)$ and $X=\GF(8)\cup\{\infty\}$, where $\GF(8)$ is
defined by the relation $i^3=i+1$.
\begin{enumerate}
\item By Theorem \ref{Sn}, $T=\PSL(2,8)$ is not arc-transitive on
  $J(9,3)$ and so as $T\norml G$ and has index three, $T$ has three
  equal sized orbits on edges. Thus the partition $\P=\{P_1,P_2,P_3\}$
  given by these three orbits is a $G$-primitive decomposition. Since
  $T$ is vertex-transitive, this is in fact a homogeneous
  factorisation and appears in \cite{Cuaresma}.

\item Let $x\in X$. Then $G_x=\AGammaL(1,8)$ and preserves the
  structure of an affine space $\AG(3,2)$ (with plane-set $\mathcal{B}$) on $X\backslash\{x\}$.
Let $$P_x=\Big\{\{A,B\}\mid A\cup B\in\mathcal{B}\Big\}$$
and $\P=\{P_x\mid x\in X\}$. Then since each 3-subset lies in a unique
  plane, $P_x\cong 14K_4$. Moreover, $G_x$ acts transitively on the
  set $\mathcal{B}$ of affine planes and for $Y\in\mathcal{B}$ we have
  $G_{x,Y}$ induces $A_4$ on $Y$. Thus $G_x$ acts transitively on the
  set of edges in $P_x$ and so given $\{A,B\}\in P_x$ we have
  $|G_{x,\{A,B\}}|=2=|G_{\{A,B\}}|$.  Thus $G_{\{A,B\}}\leqslant H$
  and so by Lemma \ref{lem:general}, $\P=P_x^G$ is a $G$-primitive
  decomposition of $J(9,3)$.

\item Let $A=\{\infty,0,1\}$ and $B=\{\infty,0,i\}$. Then
  $G_{\{A,B\}}=\la g\ra\cong C_2$ where $x^g=ix^{-1}$ and has
  orbits $\{0,\infty\}$, $\{1,i\}$, $\{i^2,i^6\}$, $\{i^3,i^5\}$ and
  $\{i^4\}$. Thus $G_{\{A,B\}}\leqslant G_{\{i^2,i^6\}}=H$ ($H$ has order 42) 
  and so by Lemma \ref{lem:general}, letting $P=\{A,B\}^H$ and $\P=P^G$ we
  obtain a $G$-primitive decomposition of $J(9,3)$. Now $H_A=\la h\ra$
  where $x^h=x+1$, which has order two and so $P$ has $21$ vertices
  and valency $2$. Moreover, $\la H_A,g\ra=D_{14}$ and so by Lemma
  \ref{lem:cosetgraph}, $P$ has three connected components. Thus
  $P\cong 3C_7$.

\item  Let $A=\{\infty,0,1\}$ and $B=\{\infty,0,i\}$. Then
  $G_{\{A,B\}}\leqslant G_{\{i^3,i^5\}}=H$ and so by
  Lemma \ref{lem:general}, letting $P=\{A,B\}^H$ and $\P=P^G$ we
  obtain a $G$-primitive decomposition of $J(9,3)$. Then 
$H_A=\la h\ra$ where $x^h=(x^4+1)^{-1}$, which has order three. Thus
  $P$ has $14$ vertices and valency 3.  Since $g$ and $h$ do not
  commute, $\la H_A,g\ra=H$ and so $P$ is a connected graph. Moreover,
  $P$ is $H$-arc-transitive and so by \cite[p167]{graphatlas}, $P$ is
  the Heawood graph.
\end{enumerate}
}%
\end{construction}

\begin{construction}
\label{con:D18}
{\rm Let $K=\GF(64)$, with $\xi$ a primitive element of $K$, and let
$F=\{0\}\cup\{(\xi^{9})^l|l=0,1,\ldots ,6\}\cong \GF(8)$. One can consider
the projective line $X$ on which $G$ acts as the elements of $K$
modulo $F$.  Then $H=\la\hat{\xi},\sigma,\tau\ra\cong D_{18}\rtimes C_3$
where $\hat{\xi}:x\rightarrow \xi x \pmod F$, 
$\sigma:x\rightarrow x^8=x^{-1} \pmod F$, and 
$\tau:x\rightarrow x^4 \pmod F$. 

\begin{enumerate}
\item Let $A=\{1,\xi,\xi^2\}$ and $B=\{\xi,\xi^2,\xi^3\}$. Then
  $\{A,B\}$ is an edge of $J(9,3)$ whose ends are interchanged by
  $\hat{\xi}^6\sigma\in H$. Thus letting $P=\{A,B\}^H$ and $\P=P^G$, Lemma
  \ref{lem:general} implies that $(J(9,3),\P)$ is a $G$-primitive
  decomposition. Now $H_A=\la\hat{\xi}^7\sigma\ra$  and so $P$ has 27
  vertices. Moreover, $H_{A,B}=1$ and so $P$ has valency 2. Since 
$\la \hat{\xi}^6\sigma,\hat{xi}^7\sigma \ra=D_{18}$ it follows from Lemma
  \ref{lem:cosetgraph} that $P$ has 3 connected components and so
  $P\cong 3C_9$.
\item  Let $A=\{1,\xi,\xi^3\}$ and $B=\{1,\xi,\xi^7\}$. Then
  $\{A,B\}$ is an edge of $J(9,3)$ whose ends are interchanged by
  $\hat{xi}^8\sigma\in H$. Thus letting $P=\{A,B\}^H$ and $\P=P^G$, Lemma
  \ref{lem:general} implies that $(J(9,3),\P)$ is a $G$-primitive
  decomposition. Now $\mid H_A\mid=1$  and so $P$ is a matching of 27
  edges. 
\item  Let $A=\{1,\xi,\xi^3\}$ and $B=\{\xi,\xi^3,\xi^4\}$. Then
  $\{A,B\}$ is an edge of $J(9,3)$ whose ends are interchanged by
  $\hat{xi}^5\sigma\in H$. Thus letting $P=\{A,B\}^H$ and $\P=P^G$, Lemma
  \ref{lem:general} implies that $(J(9,3),\P)$ is a $G$-primitive
  decomposition. Now $\mid H_A\mid=1$  and so $P$ is a matching of 27
  edges. 
\item  Let $A=\{1,\xi,\xi^3\}$ and $B=\{1,\xi^2,\xi^3\}$. Then
  $\{A,B\}$ is an edge of $J(9,3)$ whose ends are interchanged by
  $\hat{xi}^6\sigma\in H$. Thus letting $P=\{A,B\}^H$ and $\P=P^G$, Lemma
  \ref{lem:general} implies that $(J(9,3),\P)$ is a $G$-primitive
  decomposition. Now $\mid H_A\mid=1$  and so $P$ is a matching of 27
  edges. 
\end{enumerate}
}%
\end{construction}

\begin{table}
\begin{center}
\caption{$\PGammaL(2,8)$-primitive decompositions of $J(9,3)$}
\label{tab:pgammal}
\begin{tabular}{|lcc|}
\hline
$\P$ & $P$  & $G_P$ \\
\hline\hline 
$\P_{\cap}$          & $K_7$  & $D_{14}\rtimes C_3$ \\
$\P_\ominus$         & ${7\choose 2}K_2$ & $D_{14}\rtimes C_3$ \\
Construction \ref{con:pgammal}(1) & $\PSL(2,8)$-orbits & $\PSL(2,8)$\\
Construction \ref{con:pgammal}(2) &  $14K_4$ &$\AGammaL(1,8)$\\
Construction \ref{con:pgammal}(3) & $3C_7$  & $D_{14}\rtimes C_3$\\
Construction \ref{con:pgammal}(4) & Heawood graph &  $D_{14}\rtimes C_3$\\
Construction \ref{con:D18}(1)& $3C_9$ & $D_{18}\rtimes C_3$\\
Construction \ref{con:D18}(2)& $27K_2$ & $D_{18}\rtimes C_3$\\
Construction \ref{con:D18}(3)& $27K_2$ & $D_{18}\rtimes C_3$\\
Construction \ref{con:D18}(4)& $27K_2$ & $D_{18}\rtimes C_3$\\

\hline
\end{tabular}
\end{center}
\end{table}

\begin{proposition} If $(J(9,3),\P)$ is a $\PGammaL (2,8)$-primitive
  decomposition then $\P$ is as in Table $\ref{tab:pgammal}$.
 \end{proposition}
\proof
Let $G=\PGammaL (2,8)$ act on $\{\infty\}\cup \GF(8)$ and suppose that $\GF(8)$ has primitive element $i$ such that $i^3=i+1$.
Let $A=\{\infty,0,1\}$ and $B=\{\infty,0,i\}$ be adjacent vertices in
$\Gamma=J(9,3)$. Then
$G_{\{A,B\}}=G_{\{\infty,0\},\{1,i\}}=\la g\ra\cong C_2$, where
$x^g=ix^{-1}$, which fixes the point $i^4$ and has 4 orbits of size
2. Let $H$ be a maximal subgroup of $G$ containing $G_{\{A,B\}}$. The
maximal subgroups of $G$ are given in \cite[p 6]{atlas}.

If $H=\PGL(2,8)$ then we obtain the decomposition in Construction
 \ref{con:pgammal}(1) while if $H$ is a point stabiliser then $H=G_{i^4}$
 and we obtain the decomposition in Construction \ref{con:pgammal}(2).

Suppose now that $H\cong D_{14}\rtimes C_3$ is the stabiliser of a
2-subset. Then $H=G_{\{\infty,0\}}$, $H=G_{\{1,i\}}$, $H=G_{\{i^2,i^6\}}$,
or $H=G_{\{i^3,i^5\}}$. In the first case  we get the decomposition
$\mathcal{P}_{\cap}$, while the second yields $\mathcal{P}_{\ominus}$.
The third case gives Construction \ref{con:pgammal}(3) and the fourth
gives the decomposition in Construction \ref{con:pgammal}(4).

Let $H=\la\hat{\xi},\sigma,\tau\ra\cong D_{18}\rtimes C_3$ as given in
Construction \ref{con:D18}. Instead of finding all conjugates of $H$
containing $G_{\{A,B\}}$, we (equivalently) find all edge orbits $\{C,D\}^H$ such that $H$
contains $G_{\{C,D\}}$. Note that, for such an edge, $C$ and $D$ lie in
the same $H$-orbit on vertices. One sees easily that $H$ has three orbits
on vertices of $J(9,3)$, of sizes 3
($\{1,\xi^3,\xi^6\}^{\la\hat{xi}\ra}$), 27
($\{1,\xi,\xi^2\}^{\la\hat{xi}\ra}\cup\{1,\xi^2,\xi^4\}^{\la\hat{xi}\ra}\cup\{1,\xi^4,\xi^8\}^{\la\hat{xi}\ra}$),
and 54 (all the other vertices). The orbit of size 3 contains no
edges. In the orbit of size 27, if we fix the vertex
$C=\{1,\xi,\xi^2\}$, we find two vertices $D$, namely
$\{1,\xi,\xi^8\}$ and  $\{\xi,\xi^2,\xi^3\}$, such that the unique
involution switching $C$ and $D$ is in $H$. Moreover, these two vertices are
interchanged by $H_C$. Hence this vertex orbit yields one orbit of edges
whose stabilisers are contained in $H$ and we get the decomposition
in Construction \ref{con:D18}(1).

In the orbit of size 54, if we fix the vertex $C=\{1,\xi,\xi^3\}$, we
find three vertices $D$, namely $\{1,\xi,\xi^7\}$,
$\{\xi,\xi^3,\xi^4\}$ and  $\{1,\xi^2,\xi^3\}$, such that the unique
involution switching $C$ and $D$ is in $H$.  Since $H$ acts regularly
on this orbit, each choice of $D$ gives a different $H$-orbit on edges
and we get the three decompositions of  Constructions
\ref{con:D18}(2,3,4).
\Qed

\section{The case $k=2$}
By Theorem  \ref{Sn}, a subgroup $G$ of $S_n$ is arc-transitive on
$J(n,2)$ if and only if $G$ is 3-transitive. Hence other than $A_n$ or $S_n$, the possibilities for $(n,G)$ are $(11,M_{11})$,  $(12,M_{11})$, $(12,M_{12})$,
$(22,M_{22})$,  $(22,\Aut(M_{22}))$, $(23,M_{23})$, $(24,M_{24})$, $(2^d,\AGL(d,2))$ for $d>2$ , $(16,C_2^4\rtimes A_7)$, and $(q+1,G)$
where $G$ is a 3-transitive subgroup of $\PGammaL(2,q)$ with $q\geq 4$. We treat all but the last case in this section and deal with the subgroups of $\GammaL(2,q)$  in Section \ref{sec:PSL}.

\begin{proposition} If $(J(11,2),\P)$ is an $M_{11}$-primitive
  decomposition then $\P$ is $\P_{\cap}$, $\P_{\cup}$, or $\P_{\ominus}$.
\end{proposition}
\proof
Let $G=M_{11}$ act on the point set $X$ of the Witt design
$S(4,5,11)$, and let $A=\{1,2\}$, $B=\{2,3\}$ be adjacent
vertices. Then $G_{\{A,B\}}=G_{2,\{1,3\}}$ and 
since $G$ is strictly 4-transitive it follows that $|G_{\{A,B\}}|=16$ and has 
one orbit on the 8 remaining points.  Let $H$ be a maximal subgroup of
$G$ containing $G_{\{A,B\}}$. Comparing orders and the
maximal subgroups of $G$ given in \cite[p 18]{atlas} we see that 
$H\not\cong \PSL(2,11)$ or $S_5$.
It follows that $H$ stabilises either a point, a pair or a 3-subset.
In the first case $H=G_2$ and so $\P=\mathcal{P}_{\cap}$. In the second case,
$H=G_{\{1,3\}}$ and we obtain the decomposition
$\mathcal{P}_{\ominus}$, while in the last case $H=G_{\{1,2,3\}}$ and
so we get the decomposition $\mathcal{P}_{\cup}$.
\Qed

Since the stabilisers of a point and a 2-subset are maximal in $M_{11}$ it
follows from Lemma \ref{lem:classic} that $\P_\cap$ and $\P_\ominus$
are $M_{11}$-primitive decompositions of $J(12,2)$.
In order to give more constructions for $M_{11}$-primitive decompositions of
$J(12,2)$, we will need the following lemma.

\begin{lemma}\label{lem:linked3} Let $G=M_{11}$ act $3$-transitively on the point set $X$ of the Witt design $S(5,6,12)$. As seen in Construction \ref{con:M11graph}, $G$ has an orbit of length $165$ on $4$-subsets, 
forming a $3-(12,4,3)$ design with block set $\mathcal{D}$. In this
design, each $3$-set $S$ determines uniquely another $3$-set $S_
\mathcal{D}$, namely the set of fourth points of the $3$ blocks of
$\mathcal{D}$ containing $S$. We have $(S_ \mathcal{D})_
\mathcal{D}=S$ and $S\cup S_\mathcal{D}$ is a hexad of $S(5,6,12)$. Moreover if $\{S, S_\mathcal{D}, U,V\}$ is the unique linked three containing $S$ and $S_\mathcal{D}$ as triads (see Lemma $\ref{lem:LTM12}$), then $U_\mathcal{D}=V$.
\end{lemma}
\proof
For any 3-set $S$, the set $S_\mathcal{D}$ is obviously well defined
by the pro\-per\-ties of the  $3-(12,4,3)$ design. Now, an element of $G$
stabilising $S$ must also stabilise $S_\mathcal{D}$. Therefore
$G_S\leqslant G_{S_\mathcal{D}}$. Since $S_\mathcal{D}$ is also a
3-set and $G$ is 3-transitive, we must have 
$|G_S|=|G_{S_\mathcal{D}}|$. Therefore $G_S=G_{S_\mathcal{D}}$.  By a computation using {\sc Magma} \cite{magma} we find that $G_S\cong S_3\times S_3$ has orbits of lengths 3, 3 and 6 on $X$. Hence $(S_\mathcal{D})_\mathcal{D}=S$.

Let $u,v$ be two points of $S_\mathcal{D}$. Then $S\cup\{u,v\}$ is contained in a unique hexad $h$. Since $G_S$ preserves the set of hexads containing $S$, and acts
transitively on the 3 points of $S_\mathcal{D}$ and on the 6 points of $X\backslash (S\cup S_\mathcal{D})$, it follows that the sixth point of $h$ must also lie in $S_\mathcal{D}$. Hence $S\cup S_\mathcal{D}$ is a hexad.
Let $T=\{S, S_\mathcal{D}, U,V\}$ be the unique linked three
containing $S$ and $S_\mathcal{D}$ as triads (Lemma \ref{lem:LTM12}).
Since $G_S$ preserves $T$ and is transitive on $U\cup V$, it follows
that $G_S$ has an index 2 subgroup $G_{S,U}$ with orbits
$S, S_\mathcal{D}, U$ and $V$. Since the orbits of $G_{S,U}$ are a
refinement of the orbits of $G_U$, $U_\mathcal{D}$ must be one of
these orbits of size 3. Since $U_\mathcal{D}$ cannot be $S$ nor
$S_\mathcal{D}$, it follows that  $U_\mathcal{D}=V$. 
\Qed

\begin{construction}
\label{con:M11122}
{\rm Let $G=M_{11}$ act 3-transitively on the point set $X$ of the Witt design
$S(5,6,12)$. We use the notation of Lemma \ref{lem:linked3}.
\begin{enumerate}
\item 
Let $Y\in\mathcal{D}$.
Let  $$P_Y=\Big\{ \big\{\{u,x\},\{x,v\}\big\}\mid \{x,u,v\}_\mathcal{D}=Y\backslash\{x\} \Big\}$$
and $\P=\{P_Y\mid Y\in \mathcal{D}\}$.
Then $P_Y\cong 4K_2$.  Let $\{\{u,x\},\{x,v\}\}$ be an edge of  $J(12,2)$. Then it is in a unique $P_Y$, with $Y=\{x\} \cup \{x,u,v\}_\mathcal{D}$.
 Since $G_Y$ is maximal in $G$, it
 follows that $(J(12,2),\P)$ is a $G$-primitive decomposition.

\item 
 Let $T$ be a $\mathcal{D}$-linked three, that is, a linked three for the
 $S(5,6,12)$ such that, for any $X\in T$, $X_\mathcal{D}$ is a triad of $T$. 
Let $$P_T=\Big\{\big\{\{u,x\},\{x,v\}\big\}\mid \{x,u,v\} \in T \Big\}$$ 
and $\P=\{P_T\mid T \text{ is a $\mathcal{D}$-linked three}\}$.  
Then $P_T\cong 4K_3$, with each triad contributing $K_3$. 
Let $\{\{u,x\},\{x,v\}\}$ be an edge of  $J(12,2)$. Then $\{u,v,x\}$
and $\{u,v,x\}_\mathcal{D}$ must be  triads of $T$. By Lemma \ref{lem:linked3}, the unique linked three containing these two triads is a $\mathcal{D}$-linked three. It follows that
there is exactly one $\mathcal{D}$-linked three $T$ such that $P_T$
contains a given edge. Since the stabiliser in $G$ of a
$\mathcal{D}$-linked three is maximal in $G$, it follows that $(J(12,2),\P)$ is a $G$-primitive decomposition.
\end{enumerate}
}%
\end{construction}

Thus we have the $M_{11}$-primitive decompositions listed in Table \ref{tab:M112}.

\begin{table}
\begin{center}
\caption{$M_{11}$-primitive decompositions of $J(12,2)$}
\label{tab:M112}
\begin{tabular}{|lcc|}
\hline
$\P$ & $P$  & $G_P$ \\
\hline\hline 
$\P_{\cap}$          & $K_{11}$  & $\PSL(2,11)$ \\
$\P_\ominus$         & $10K_2$ & $S_5$ \\
Construction \ref{con:M11122}(1) & $4K_2$ & $M_8\rtimes S_3$ \\
Construction \ref{con:M11122}(2) & $4K_3$ &  $M_9\rtimes C_2$ \\
\hline
\end{tabular}
\end{center}
\end{table}

\begin{proposition} If $(J(12,2),\P)$ is an $M_{11}$-primitive
  decomposition then $\P$ is given by Table $\ref{tab:M112}$.
\end{proposition}
\proof
Let $G=M_{11}$ act transitively on the point set $X$ of the Witt design $S(5,6,12)$ and let $\mathcal{D}$ be the block set of the  $3-(12,4,3)$ design described in Construction \ref{con:M11graph} (see above). Take adjacent vertices
$A=\{1,2\}$ and $B=\{2,3\}$.   Then $G_{\{A,B\}}=G_{2,\{1,3\}}\cong
D_{12}$ which has an orbit of length 3 (namely, $\{1,2,3\}_\mathcal{D}$) and an orbit of length 6 on the
remaining 9 points of $X$. Let $H$ be a maximal subgroup of $G$
containing $G_{\{A,B\}}$. Since $M_{10}$ contains no elements of order
6, it follows that $H\not\cong M_{10}$. If $H$ is a point stabiliser,
then $H=G_2$ and we get the decomposition $\mathcal{P}_{\cap}$. If $H$
is a pair stabiliser then $H=G_{\{1,3\}}$, and we get the
decomposition $\mathcal{P}_{\ominus}$.  If $H\cong M_8\rtimes S_3$
then $H$ is the stabiliser of a block in $\mathcal{D}$. There is a unique such block, namely the union of $\{2\}$ with $\{1,2,3\}_\mathcal{D}$. Hence $H$ is  the stabiliser of the divisor of the decomposition obtained from Construction
\ref{con:M11122}(1) containing $\{A,B\}$.

Now let $H\cong M_9\rtimes S_3$. Then $H$ is a $\mathcal{D}$-linked three stabiliser, namely the only one containing $\{1,2,3\}$ as a triad (see the construction). Hence $H$ is  the stabiliser of the divisor of the decomposition obtained from Construction \ref{con:M11122}(2) containing $\{A,B\}$.
\Qed

\begin{proposition} If $(J(12,2),\P)$ is an $M_{12}$-primitive
  decomposition, then $\P$ is $\P_{\cup}$, $\P_{\cap}$ or $\P_{\ominus}$.
\end{proposition}
\proof
Let  $G=M_{12}$ act on the point set $X$ of the Witt-design
$S(5,6,12)$ and take adjacent vertices
$A=\{1,2\}$ and $B=\{2,3\}$. Then
$G_{\{A,B\}}=G_{2,\{1,3\}}$ which has order
144 and is 2-transitive on the 9 remaining points since $G$ is
5-transitive on $X$. Let $H$ be a maximal subgroup of $G$ containing
$G_{\{A,B\}}$. The maximal subgroups of $G$ are given in
\cite{atlas}, and comparing orders we see that $H\not\cong \PSL(2,11)$, 
$2\times S_5$, $4^2: D_{12}$, $M_8.S_4$ or $A_4\times S_3$. Since $G_{\{A,B\}}$ fixes a point but not a hexad it follows that $H$ is not the stabiliser of a hexad pair, and since $G_{\{A,B\}}$ is 2-transitive on $X\backslash\{1,2,3\}$ we
also  have that $H$ is not the stabiliser of a linked three. In the
action of $M_{11}$ on 12 points, $\PSL(2,11)$ is the stabiliser of a
point. Since $144$ does not divide $|\PSL(2,11)|$ and $G_{\{A,B\}}$
fixes the point 2, it follows that $H$ is not a transitive copy of
$M_{11}$. Thus $H=G_2, G_{\{1,3\}}$ or $G_{\{1,2,3\}}$. In the first
case we get the decomposition $\P_\cap$, the second case yields
$\P_\ominus$ while the third gives $\P_\cup$.
\Qed

Before dealing with $G=M_{22}$ we need the following well known result
which follows from Lemma \ref{heptads}.

\begin{lemma}\label{hexads-22}
Let $(X,\mathcal{B})$ be the Witt design $S(3,6,22)$. Then  $\mathcal{B}$ contains $77$ elements, called \emph{hexads}. Each point of $X$ is contained in $21$ hexads, each $2$-subset in $5$ hexads, and each $3$-subset in a unique hexad. Moreover, the stabiliser of a hexad is
$C_2^4\rtimes A_6$ with the pointwise stabiliser of the hexad being
$C_2^4$ which acts regularly on the $16$ points not in the hexad.
\end{lemma}
\proof
Since $(X,\mathcal{B})$ can be derived from the set of blocks of the Witt design $S(4,5,23)$ containing a given point, this follows from Lemma \ref{heptads}.
\Qed

\begin{proposition} If $(J(22,2),\P)$ is an $M_{22}$-primitive
  decompositions then $\P=\P_{\cap}$ or $\P_{\ominus}$,  or $\P$ is
  obtained from Construction \ref{con:design} and has divisors isomorphic
  to $J(6,2)$.
\end{proposition}
\proof
Let $G=M_{22}$ act on the point-set $X$ of the Witt design
$S(3,6,22)$ and take adjacent vertices $A=\{1,2\}$ and $B=\{2,3\}$.
Moreover, suppose that $h=\{1,2,3,4,5,6\}$ is the unique hexad of the
Witt design containing $\{1,2,3\}$. By Lemma \ref{hexads-22}, 
$G_h=C_2^4\rtimes A_6$, where $C_2^4$ acts trivially on $h$ and transitively on 
$X\backslash h$. It follows that $G_{\{A,B\}}=G_{2,\{1,3\},\{4,5,6\}}$
had order 96 and acts transitively on $X\backslash h$.

Let $H$ be a maximal subgroup of $G$ containing $G_{\{A,B\}}$.
 Comparing orders and the
maximal subgroups of $G$ given in \cite{atlas} we see that 
$H\not\cong \PSL(2,11)$, $A_7$ or $M_{10}$.
Since $G_{\{A,B\}}$ does not stabilise an octad, it follows that $H$ is
either $G_2$, $G_{\{1,3\}}$ or $G_h$. The first gives
the decomposition $\mathcal{P}_{\cap}$, while the second yields
$\mathcal{P}_{\ominus}$.  Finally $G_h$ is the stabiliser of the part
 of the decomposition obtained from Construction \ref{con:design}
 containing $\{A,B\}$ and has divisors isomorphic to $J(6,2)$.
\Qed

\begin{proposition} All $\Aut (M_{22})$-primitive decompositions of
  $J(22,2)$ are $M_{22}$-primitive decompositions.
\end{proposition}
\proof
By \cite{atlas}, a maximal subgroup of $\Aut(M_{22})$ is either
$M_{22}$ or arises from a maximal subgroup of $M_{22}$. Since $M_{22}$
is arc-transitive it does not give a decomposition. In all other
cases, Lemma \ref{lem:restrict} implies that we get $M_{22}$-primitive
decompositions.
\Qed

\begin{proposition} If $(J(23,2),\P)$ is an  $M_{23}$-primitive
 decomposition  then $\P$ is $\P_{\cap}$, $\P_{\ominus}$ or $\P_{\cup}$.
\end{proposition}
\proof
Let $G=M_{23}$ act on the point-set $X$ of the Witt design $S(4,7,23)$
and take adjacent vertices $A=\{1,2\}$ and $B=\{2,3\}$. Then 
$G_{\{A,B\}}=G_{2,\{1,3\}}\cong 2^4\rtimes S_5$ (see \cite[p 71]{atlas}) and since $G$ is 4-transitive, $G_{\{A,B\}}$ is transitive on
$X\backslash\{1,2,3\}$. Let $H$ be a maximal subgroup of $G$
containing $G_{\{A,B\}}$.
Since $|G_{\{A,B\}}|$ does not divide $23.11$, it follows from
\cite[p 71]{atlas} that $H$ is intransitive. Hence $H$ is $G_2$,
$G_{\{1,3\}}$ or $G_{\{1,2,3\}}$. These give us the decompositions
$\mathcal{P}_{\cap}$, $\mathcal{P}_{\ominus}$ and $\mathcal{P}_{\cup}$
respectively.
\Qed

\begin{proposition} If $(J(24,2),\P)$ is an $M_{24}$-primitive
  symmetric decompositions then $\P$ is $\P_\cap$, $\P_\ominus$
  or $\P_\cup$.
 \end{proposition}
\proof
Let $G=M_{24}$ acting on the point-set $X$ of the Witt design
$S(5,8,24)$ and take adjacent vertices $A=\{1,2\}$ and $B=\{2,3\}$. Then
$G_{\{A,B\}}=G_{2,\{1,3\}}\cong \PSigmaL(3,4)$  (see \cite[p 96]{atlas}). Note that $G_{\{A,B\}}$ is transitive on $X\backslash\{1,2,3\}$ since $G$ is
5-transitive on $X$. Let $H$ be a maximal subgroup of $G$ containing
$G_{\{A,B\}}$. Looking at the maximal subgroups of $G$ in
\cite{atlas}, it follows that $H$ is either $G_2$, $G_{\{1,3\}}$ or
$G_{\{1,2,3\}}$.  Thus we obtain the decompositions $\mathcal{P}_{\cap}$, $\mathcal{P}_{\ominus}$ and $\mathcal{P}_{\cup}$ respectively.
\Qed

Let $G=\AGL(d,2)$ acting on the set $X$ of vectors of a $d$-dimensional vector space over $\GF(2)$.
Since the stabiliser of a vector is maximal in $G$, Lemma
\ref{lem:classic} implies that $\P_\cap$ is a $G$-primitive
decomposition. The set of affine planes in the affine space
$\AG(d,2)$ yields an $S(3,4,2^d)$ Steiner system with each point
contained in $\frac{(2^d-1)(2^{d-1}-1)}{3}$ planes. In both cases, $G$ acts transitively on planes hence we can use Construction \ref{con:design}. However, $G$ is not primitive on planes as it preserves parallelness. Applying now Construction \ref{con:primitive} yields line 2 of Table \ref{tab:AGL}.
As $G$ is transitive on points and the stabiliser of a point is maximal in $G$, applying Construction \ref{con:design2} yields line 3 of Table \ref{tab:AGL}.
As $G$ is 2-transitive, we can use Construction \ref{con:design4}. However, $G$ acts
imprimitively on $2$-subsets as $2$-subsets correspond to lines and
again $G$ preserves parallelness. Thus we also apply Construction \ref{con:primitive} and obtain line 4 of Table \ref{tab:AGL}.
Indeed the divisors are indexed by lines of the
affine plane and are isomorphic to $2^{d-2}K_2$. Each pair $Y_1,Y_2$ of parallel lines yields a $C_4$ in the $J(4,2)$ induced on $Y_1\cup Y_2$. As a parallel class of lines contains $2^{d-1}$ lines, we have $\frac{2^{d-1}(2^{d-1}-1)}{2}$ pairs of parallel lines in the imprimitivity class. 

When $d=4$ the group $\overline{G}=C_2^4\rtimes A_7<\AGL(4,2)$ is 3-transitive on $X$ and hence, by Corollary \ref{Sncor}, is arc-transitive on $J(2^4,2)$. Thus the four $G$-primitive decompositions in Table \ref{tab:AGL} are also $\overline{G}$-transitive. The stabiliser in $\overline{G}$ of a point is $A_7$ which is maximal in $\overline{G}$. Hence the partitions in Rows 1 and 3 are $\overline{G}$-primitive. The stabilisers of 2-spaces and 1-spaces in $A_7$ are maximal in $A_7$ and so the remaining two partitions are also $\overline{G}$-primitive.

Before showing that these are the only $G$-primitive decompositions with $G\leqslant \AGL(d,2)$ we need a lemma. 

\begin{lemma}
\label{lem:affinemax}
Let $G=N\rtimes G_0$ where $N\cong C_p^d$ for some prime $p$ and $G_0$
acts irreducibly on $N$. Suppose that $H$ is a maximal subgroup of
$G$. Then either $H$ is a complement of $N$, or $M=N\rtimes H_0$ for
some maximal subgroup $H_0$ of $H$.
\end{lemma}
\proof
Since $H$ normalises $N$ we have
$H\leqslant NH\leqslant G$. Thus as $H$ is maximal, either $NH=H$ or
$NH=G$. The first case implies that $N\leqslant H$ and so 
$H=N\rtimes H_0$ for some maximal subgroup $H_0$ of $G_0$. Suppose now
that $NH=G$. Then $H/(H\cap N)\cong G_0$, and so for each $g\in G_0$,
there exists $n\in N$ such that $ng\in H$. Since $N$ is abelian, it
follows that $H$ induces $G_0$ in its action on $N$ by conjugation.
Since $G_0$ acts irreducibly on $N$ and $H$ normalises $H\cap N$, it
follows that $H\cap N=1$ or $N$. However, $H\cap N=N$ implies that $H=G$ which is not the case. Hence $H\cap N=1$ and $H\cong G_0$, that is $H$ is a complement
of $N$.
\Qed

\begin{table}
\begin{center}
\caption{$G$-primitive decompositions of $J(2^d,2)$ for $G=\AGL(d,2)$ with $d\geq 3$, or $G=C_2^4\rtimes A_7$ with $d=4$}
\label{tab:AGL}
\begin{tabular}{|lcc|}
\hline
$\P$ & $P$  & $G_P$ \\
\hline\hline 
$\P_{\cap}$          & $K_{2^d-1}$  & $G_0$\\
Constructions \ref{con:design} and \ref{con:primitive} &
$2^{d-2}J(4,2)\cong 2^{d-2}K_{2,2,2}$ &  $C_2^d\rtimes (G_0)_{\la v,w\ra}$\\
Construction \ref{con:design2} & $\frac{(2^d-1)(2^{d-1}-1)}{3}K_3$ &
$G_{v+w}$\\
Construction \ref{con:design4} and \ref{con:primitive} &
$2^{d-2}(2^{d-1}-1)C_4$ & $C_2^d\rtimes (G_0)_{\la v+w\ra}$\\
\hline
\end{tabular}
\end{center}
\end{table}

\begin{proposition}
Let $d\geq 3$ and $G=\AGL(d,2)$, or $d=4$ and $G=C_2^4\rtimes A_7$. If $(J(2^d,2),\P)$ is a $G$-primitive decomposition then $\P$ is given by Table $\ref{tab:AGL}$.
\end{proposition}
\proof
We can identify $X$ with the vectors of a $d$-dimensional vector space over $\GF(2)$. Let $A=\{0,v\}$ and $B=\{0,w\}$ where $v,w$ are distinct non-zero elements of $X$.
Thus $G_{\{A,B\}}=(G_0)_{\{v,w\}}$ which is an index 3 subgroup of 
$(G_0)_{\la v,w\ra}$. Moreover, $G_{\{A,B\}}$ fixes the vector $v+w$ and is 
transitive on  $X\setminus\la v,w\ra$.

Let $H$ be a maximal subgroup of $G$ containing $G_{\{A,B\}}$. By
Lemma \ref{lem:affinemax}, either $H$ is a complement of $N=\soc(G)$ or
$H=N\rtimes H_0$ for some maximal subgroup $H_0$ of $G_0$.

Suppose first that $H$ is a complement.  By a 
{\sc Magma} \cite{magma} calculation, $C_2^4\rtimes A_7$ has a unique conjugacy class of complements. If $d\geq 4$ then there is a unique class of complements of $N$ in $\AGL(d,2)$, while
in $\AGL(3,2)$ there are two classes (see for example \cite{dempwolff}). Hence either $H$ is the
stabiliser of a vector or $d=3$ and $H$ is transitive. In the second case $H=\PSL(2,7)$ acting transitively on $V$. However, a Sylow 2-subgroup of $H$ is then regular on $V$, and hence $H$ cannot contain $G_{\{A,B\}}\cong D_8$ (fixing the point 0). Thus $H$ is the stabiliser of a vector and so $H=G_0$ or $G_{v+w}$. The first case yields the decomposition $\P_{\cap}$, while the second is the stabiliser of the divisor of the decomposition obtained from Construction \ref{con:design2} containing $\{A,B\}$.

Suppose now that $H=N\rtimes H_0$ for some maximal subgroup $H_0$ of $G_0$. 
First let $G=\AGL(d,2)$. Since $G_{\{A,B\}}$ is an index 3 subgroup of the stabiliser in $\GL(d,2)$ of the 2-space $\la v,w\ra$, it contains a Sylow 2-subgroup of $\GL(d,2)$. Thus $H_0$ contains a Sylow 2-subgroup of $\GL(d,2)$ and it follows from a Lemma of Tits (see for example \cite[(1.6)]{seitz}) that $H_0$ is a parabolic subgroup and hence is a subspace stabiliser. Now let $G=C_2^4\rtimes A_7$. Since $G_{\{A,B\}}\cong S_4$ fixes a nonzero vector it is contained in a subgroup $\PSL(2,7)$ of $A_7$ and hence by \cite[p 10]{atlas}, the elements of order 3 in $G_{\{A,B\}}$ are from the conjugacy class $3B$, that is, in the representation of $A_7$ on 7 points they are products of two 3-cycles. By \cite[p 10]{atlas}, $A_7$ has 5 conjugacy classes of maximal subgroups. The elements of order 3 in a maximal $S_5$ subgroup are from the conjugacy class $3A$ (\cite[p 10]{atlas}), instead of $3B$ and so $H_0\not\cong S_5$. If $H_0\cong A_6$ then $A_6\cong \PSp(4,2)'$ and contains two conjugacy classes of $S_4$ subgroups. One is the stabiliser of a vector and has orbit lengths 1, 6 and 8 on nonzero vectors and the other is the stabiliser of a totally isotropic 2-space with orbit sizes 3 and 12. Hence none of them stabilises the pair $\{v,w\}$ and so $H_0\not\cong A_6$. The remaining three conjugacy classes of maximal subgroups of $A_7$ are stabilisers of subspaces. Thus for both groups $G$, $H_0$ is a subspace stabiliser. The only proper, nontrivial subspaces fixed by $G_{\{A,B\}}$ are $\la v+w\ra$ and $\la v,w\ra$.   If $H_0=(G_0)_{\la v,w\ra}$ then $H$ is the stabiliser of the class of planes parallel to $\la v,w\ra$ and so $H$ is the stabiliser of the divisor containing $\{A,B\}$  of the decomposition in Row 2 of Table \ref{tab:AGL}. Similarly, if $H_0= (G_0)_{\la v+w\ra}$ then $H$ is the stabiliser of the class of lines parallel to $\la v+w\ra$ and so is the stabiliser of the divisor containing $\{A,B\}$ of the decomposition in Row 4 of Table \ref{tab:AGL}.
\Qed

\section{Completing the case $k=2$: $G\leqslant\PGammaL(2,q)$}
\label{sec:PSL}

In this section we determine all $G$-primitive decompositions of
$J(q+1,2)$ where $G$ is a 3-transitive subgroup of $\PGammaL(2,q)$ for
$q=p^f\geq 4$ with $p$ a prime. The group $\PGL(2,q)$ is the group of
all fractional linear transformations  
$$t_{a,b,c,d}: z\mapsto \frac{az+b}{cz+d}, \hspace{1cm} ad-bc\neq 0$$
of the projective line $X=\{\infty\}\cup\GF(q)$ with the conventions 
$1/0=\infty$ and $(a\infty+b)/(c\infty+d)=a/c$.
Note that $t_{a,b,c,d}=t_{a',b',c',d'}$ if and only if
$(a,b,c,d)=\lambda(a',b',c',d')$ for some $\lambda\neq 0$. The group
$\PSL(2,q)$ is then the set of all $t_{a,b,c,d}$ such that $ad-bc$ is
a square in $\GF(q)$. The Frobenius map $\phi:z\mapsto z^p$ also acts
on $X$ and $\phi^{-1}t_{a,b,c,d}\phi=t_{a^p,b^p,c^p,d^p}$.  Then
$\PGammaL(2,q)=\la \PGL(2,q),\phi\ra$.  Another interesting family of subgroups
of $\PGammaL(2,q)$ occurs when $p$ is odd and $f$ is even. In this
case we can define for each divisor $s$ of $f/2$, the group
$M(s,q)=\la \PSL(2,q), \phi^st_{\xi,0,0,1}\ra$, where  $\xi$ is a primitive element of $\GF(q)$. Each  $g \in\PGL(2,q)\setminus\PSL(2,q)$ can be
written as $t_{\xi,0,0,1}h$ for some $h\in\PSL(2,q)$, and so
$\phi^sg\in M(s,q)$. It was shown in 
\cite[Theorem 2.1]{GPZ} that $G$ is a 3-transitive subgroup of
$\PGammaL(2,q)$ if and only if either $G$ contains $\PGL(2,q)$, or $G=M(s,q)$ for
some $s$.

We begin with the following construction.
\begin{construction}\label{con:PGL1} 
{\rm \cite{Cuaresma} Let $X=\{\infty\}\cup \GF(q)$ be the projective line, $H=\PSL(2,q)$
and $q\equiv 1\pmod 4$. Then $H$ is has two equal sized orbits on
edges, namely $P_\square=\{\{\infty,0\},\{\infty,1\}\}^H$, and
$P_{\not\square}=\{\{\infty,0\},\{\infty,t\}\}^H$, with  $t$ not a
square in $ \GF(q)$.  Thus the partition $\P=\{P_\square,
P_{\not\square}\}$ is a $G$-primitive decomposition of  $J(q+1,2)$ for
any 3-transitive subgroup $G$ of $\PGammaL(2,q)$.  The divisors are
complementary spanning graphs $\Theta$ of valency $q-1$.
}%
\end{construction}
 
\begin{proposition}\label{prop:PSL}
Let $G$ be a $3$-transitive subgroup of $\PGammaL(2,q)$ and let $\P$ be
a $G$-primitive decomposition of $J(q+1,2)$ such that $\PSL(2,q)$
fixes a part. Then  $q\equiv 1\pmod 4$ and $\P$ is obtained from Construction $\ref{con:PGL1}$.
\end{proposition}
\proof
The graph  $J(q+1,2)$ contains $\frac{q(q^2-1)}{2}$ edges. For $q$ even, $|\PSL(2,q)| =q(q^2-1)$ and an edge stabiliser has order 2, so $\PSL(2,q)$ is transitive on edges. Thus $q$ is odd and so
$|\PSL(2,q)|=\frac{q(q^2-1)}{2}$. Whenever $(q-1)/2$ is odd, the stabiliser in $\PSL(2,q)$ of a point of $X$ has odd order. Since the stabiliser of
the edge $\{\{x,y\},\{x,z\}\}$ fixes $x$ and interchanges $y$ and $z$,
it follows that no nontrivial element of $\PSL(2,q)$ fixes an edge and
so $\PSL(2,q)$ is edge-transitive. Hence
$(q-1)/2$ is even and $\PSL(2,q)$ has two equal length orbits on edges, giving the $G$-primitive decomposition of Construction \ref{con:PGL1} for any 3-transitive subgroup $G$ of $\PGammaL(2,q)$.
\Qed

To classify all $G$-primitive decompositions with $G$ a 3-transitive
subgroup of $\PGammaL(2,q)$ we require knowledge of the maximal
subgroups of all such $G$. First we note the following theorem.
\begin{theorem} {\rm \cite[Corollary 1.2]{PSLmax}}
\label{thm:dropdown}
Let $\PGL(2,q)\leqslant G\leqslant \PGammaL(2,q)$ and suppose that $H$ is a maximal subgroup of $G$ not containing $\PSL(2,q)$. Then $H\cap \PGL(2,q)$ is maximal in $\PGL(2,q)$.
\end{theorem}

Theorem \ref{thm:dropdown} and Lemma \ref{lem:restrict} imply that we only need to find all
$\PGL(2,q)$-primitive and all $M(s,q)$-primitive decompositions. We
now state all maximal subgroups of these two groups. The first is well known and follows from Dickson's classification \cite{dickson} of subgroups of $\PSL(2,q)$, see also \cite{PSLmax}.
\begin{theorem}
\label{thm:PGLmax}
Let $G=\PGL(2,q)$ with $q\geq 4$ a power of the prime $p$. Then the maximal subgroups of $G$ are:
\begin{enumerate}
\item $[q]\rtimes C_{q-1}$.
\item $D_{2(q-1)}$, $q\neq 5$.
\item $D_{2(q+1)}$.
\item $S_4$ if $q=p\equiv \pm 3\pmod 8$.
\item $\PGL(2,q_0)$ where $q=q_0^r$ with $q_0> 2$, $r$ is a prime and $r$ is odd if $q$ odd.
\item $\PSL(2,q)$, $q$ odd.
\end{enumerate}
\end{theorem}

\begin{theorem}{\rm \cite[Theorem 1.5]{PSLmax}}
\label{thm:Msqmax}
Let $G=M(s,q)$ with $q=p^f\geq 3$ for $p$ odd and $f$ even, and $s$ a
divisor of $f/2$. Then the maximal subgroups of $G$ which do not
contain $\PSL(2,q)$ are:
\begin{enumerate}
\item stabiliser of a point of the projective line,
\item $N_G(D_{q-1})$,
\item $N_G(D_{q+1})$,
\item $N_G(\PSL(2,q_0))$ where $q=q_0^r$ with $r$ an odd prime.
\end{enumerate} 
\end{theorem}

We require the following knowledge about the stabiliser of an
edge.
\begin{lemma}
\label{lem:Gedge}
Let $e=\{\{\infty,0\},\{\infty,1\}\}$. Then 
\begin{enumerate}
\item $\PGL(2,q)_e=\la t_{-1,1,0,1}\ra$,
\item $\PGammaL(2,q)_e=\la t_{-1,1,0,1},\phi\ra$ of order $2f$, and
\item $M(s,q)_e= \la t_{-1,1,0,1},\phi^{2s}\ra$ of order $f/s$.
\end{enumerate}
\end{lemma}
\proof
Since $\PGL(2,q)$ is sharply 3-transitive,  $\PGL(2,q)_e=\la g\ra$ where
$g$ fixes $\infty$ and interchanges $0$ and $1$. Thus $\PGL(2,q)_e$ is
as in the lemma.  Since $\phi$ fixes $\infty$, 0 and 1, the second claim
follows. By \cite[Corollary 2.2]{GPZ},
$M(s,q)_{\infty,0,1}=\la \phi^{2s}\ra$ and since $q$ is an even power
of a prime we have $q\equiv 1\pmod 4$. Thus 
$t_{-1,1,0,1}\in\PSL(2,q)$ and so $M(s,q)_e$ is as given by the
lemma. 
\Qed

Instead of finding all maximal subgroups $H$ containing the stabiliser
of a fixed edge $\{A,B\}$ we solve the equivalent problem of choosing a representative $H$ from each conjugacy class of maximal subgroups and finding all edges whose edge stabiliser is contained in $H$. See Remark \ref{rem:2ways}.

\begin{construction}\label{con:PGL2}
{\rm Let $X=\{\infty\}\cup \GF(q)$ be the projective line with $q$ odd and
let $H=\PGammaL(2,q)_\infty=\AGammaL(1,q)$. Let $e=\{\{0,1\},\{0,-1\}\}$. The
 stabiliser in $\PGammaL(2,q)$ of $e$ is $\la \phi,t_{-1,0,0,1}\ra$, which is
 contained in $H$. Moreover $H$ is a maximal subgroup of $\PGammaL(2,q)$. Thus 
by Lemma \ref{lem:general}, letting 
$$P=e^H=\Big\{\big\{\{i,i+j\},\{i,i-j\}\big\}\mid
i,j\in\GF(q), i\neq j\Big\}$$ and $\P=P^{\PGammaL(2,q)}$, we obtain a
$\PGammaL(2,q)$-primitive decomposition of \mbox{$J(q+1,2)$}. The divisors have valency
2 and hence are a union of cycles. Since $\GF(q)$ has characteristic
$p$ it follows that each cycle has length $p$ and so the divisors are
isomorphic to $\frac{q(q-1)}{2p}C_p$. For any 3-transitive group $G$
with socle $\PSL(2,q)$,  $H\cap G$ is maximal in $G$ and so $\P$ is
$G$-primitive by  Lemma \ref{lem:restrict}.
}%
\end{construction}

\begin{proposition}\label{lem:AGL}
Let $(J(q+1,2),\mathcal{P})$ be a $G$-primitive decomposition with $G$ a $3$-transitive subgroup of $\PGammaL(2,q)$ such that, for $P\in\P$, $G_P$ is the stabiliser of a point of the projective line. Then either $\P=\P_{\cap}$ with divisors $K_q$ or $q$ is a power of an odd prime $p$ and $\P$ is obtained by Construction $\ref{con:PGL2}$.
\end{proposition} 
\proof
Let $P\in\P$ and $\Gamma=J(q+1,2)$. Then without loss of generality we may suppose that $H=G_P$ is the stabiliser of the point $\infty$ of 
$X=\{\infty\}\cup \GF(q)$. We recall that $G$ either contains $\PGL(2,q)$ or is $M(s,q)$ for some $s$. Thus $H$ acts 2-transitively on $\GF(q)$ and so the orbits of $H$ on $V\Gamma$ are  $O_1=\{\{\infty,x\}\mid x\in\GF(q)\}$ and $O_2=\{\{x,y\}\mid x,y\in \GF(q)\}$. If $\{A,B\}\in P$ then $H$ contains the stabiliser in $G$ of $\{A,B\}$ and so either $\{A,B\}\subseteq O_1$ or $\{A,B\}\subseteq O_2$.  Note that $P=\{A,B\}^H$.

Since $H$ is 2-transitive on $\GF(q)$ it follows that $H$ acts
transitively on the set of arcs between vertices of $O_1$ and so $H$
contains the stabiliser in $G$ of every edge between vertices of
$O_1$. Thus if $\{A,B\}\subseteq O_1$ then
$$\{A,B\}^H=\Big\{\big\{\{\infty,x\},\{\infty,y\}\big\}\mid
x,y\in\GF(q)\Big\} \cong K_q.$$
Hence $\P= \mathcal{P}_{\cap}$.

Suppose now that $\{A,B\}\subseteq O_2$. We may suppose that $A=\{0,1\}$ and $B=\{0,b\}$ for some $b\in\GF(q)\backslash\{0,1\}$. Let $g=t_{0,b,1-b,b}\in \PGL(2,q)$. Then $g$ maps $\infty\rightarrow 0\rightarrow1\rightarrow b$ and so $G_{\{A,B\}}=G_{\{\{\infty,0\},\{\infty,1\}\}}^g$ (this is obvious if $G$ contains  $\PGL(2,q)$ and follows from the fact that $M(s,q)\norml\PGammaL(2,q)$ for $G=M(s,q)$). By Lemma \ref{lem:Gedge}, $t_{-1,1,0,1}^g\in G_{\{A,B\}}\leqslant H=G_\infty$, and since $g$ does not fix $\infty$ and the only fixed points of $t_{-1,1,0,1}$ are $\infty$ and $2^{-1}$ (only if $q$ is odd), it follows that $q$ is odd and $g:2^{-1}\rightarrow\infty$. This implies that $b=-1$.  Hence $\phi^g$ fixes $\infty$ and so  by Lemma \ref{lem:Gedge},  $G_{\{\{0,1\},\{0,-1\}\}}\leqslant H$ in all cases. Hence $\P$ is the decomposition of  Construction \ref{con:PGL2}.
\Qed

\subsection{$D_{q-1}$ subgroups}
\begin{construction}\label{con:PGL3}
{\rm Let $X=\{\infty\}\cup \GF(q)$ be the projective line where $q=p^f$ for
some odd prime $p$ and let $\xi$ be a primitive element of
$\GF(q)$. Then  $\PGammaL(2,q)_{\{0,\infty\}}=\la
t_{\xi,0,0,1},t_{0,1,1,0},\phi\ra\cong D_{2(q-1)}\rtimes C_f$.
\begin{enumerate}
\item  Let $H=\PGammaL(2,q)_{\{0,\infty\}}$ and $e=\{\{0,1\},\{0,-1\}\}$. Then $t_{-1,0,0,1}\in H$  interchanges the two vertices of $e$ while $\phi$ fixes each of the vertices of $e$. Hence $H$ contains the stabiliser in $\PGammaL(2,q)$ of $e$ and $H$ is a maximal subgroup of $\PGammaL(2,q)$ for $q\neq 5$. Thus by Lemma \ref{lem:general}, letting
$$P=e^H= \Big\{\big\{\{x,y\},\{x,-y\}\big\}\mid 
x\in\{0,\infty\},y\in\GF(q)\backslash\{0\}\Big\}$$ and
$\P=P^{\PGammaL(2,q)}$, we obtain a $\PGammaL(2,q)$-primitive decomposition of
$J(q+1,2)$. The divisors are isomorphic to $(q-1)K_2$ since the stabiliser of the vertex $\{0,1\}$ in $H$ is $\la\phi\ra$, which fixes $\{0,-1\}$. For any $3$-transitive subgroup $G$ of $\PGammaL(2,q)$, we have $H\cap G$ is maximal in $G$ and so $\P$ is a $G$-primitive decomposition  by  Lemma \ref{lem:restrict}.

\item Let $i<\frac{q-1}{2}$ and $l$ be an integer such that $\phi^l$ fixes the set $\{\xi^i,\xi^{-i}\}$. Let $G=\la \PGL(2,q),\phi^l\ra$ and $H=G_{\{\infty,0\}}=\la t_{\xi,0,0,1},t_{0,1,1,0},\phi^l\ra$.  The automorphism of $\PGL(2,q)$ switching the vertices of the edge $e=\{\{1,\xi^i\},\{1,\xi^{-i}\}\}$ is $t_{0,1,1,0}$, while either $\phi^l$ or $t_{0,1,1,0}\phi^l$ fixes both vertices of $e$. Hence $G_e<H$ and $H$ is a maximal subgroup of $G$ for $q\neq 5$. Hence by Lemma \ref{lem:general}, letting  $$P=e^H=\Big\{\big\{ \{x,\xi^ix\}, \{x,\xi^{-i}x\} \big\}\mid x\in\GF(q)\backslash\{0\}\Big\}$$ and $\P=P^G$, we obtain a $G$-primitive decomposition of $J(q+1,2)$. The divisors have valency 2 and  hence are a union of cycles. These cycles have length the order of $\xi^i$, which is $\frac{q-1}{(q-1,i)}$. Thus each divisor is isomorphic to $(q-1,i)C_{\frac{q-1}{(q-1,i)}}$. In fact for any 3-transitive subgroup $\overline{G}$ of $G$, $H\cap \overline{G}$ is maximal in $\overline{G}$ and so $\P$ is a $\overline{G}$-primitive decomposition.
\end{enumerate}
}%
\end{construction}

\begin{proposition}\label{lem:D2(q-1)}
Let $(J(q+1,2),\mathcal{P})$ be a $G$-primitive decomposition such that $\PGL(2,q)\leqslant G \leqslant \PGammaL(2,q)$ and for $P\in\P$ we have $G_P=N_G(D_{2(q-1)})$. Then either $\P=\P_{\ominus}$,  or $q$ is odd and  $\P$ is obtained by Construction $\ref{con:PGL3}(1)$, or $\P$ is obtained by Construction $\ref{con:PGL3}(2)$.
\end{proposition}
\proof
Let $P\in\P$. Since $G_P\cap \PGL(2,q)$ is a maximal subgroup of
$\PGL(2,q)$, by Lemma \ref{lem:restrict}, $\P$ is a
$\PGL(2,q)$-primitive decomposition. Thus we may suppose that
$G=\PGL(2,q)$ and $H=G_P=\la t_{\xi,0,0,1},t_{0,1,1,0}\ra\cong
D_{2(q-1)}$. The orbits of $H$ on vertices are
$\{\{0,\infty\}\}$, 
$$O_0=\{\{x,y\}\mid x\in\{0,\infty\},y\in\GF(q)\backslash\{0\}\}$$ 
and
$$O_i=\{\{x,\xi^ix\}\mid x\in\GF(q)\backslash\{0\}\}$$ for each
$i\leq \frac{q-1}{2}$. Note that $|O_0|=2(q-1)$. When $q$ is even
there are $q/2-1$ orbits $O_i$, each having length $q-1$. When $q$ is odd there are $\frac{q-3}{2}$ of length $q-1$ and one, $O_{\frac{q-1}{2}}$, of length $\frac{q-1}{2}$.

If $\{A,B\}\in P$ then $H$ contains the stabiliser in $G$ of $\{A,B\}$
and so $\{A,B\}$ is contained in one of the orbits of $H$ on vertices.
Note that $P=\{A,B\}^H$.

Suppose first that $\{A,B\}\subseteq O_0$. Without loss, let
$A=\{0,1\}$. Then the neighbours of $A$ in $O_0$ are $\{\infty,1\}$ and
$\{0,y\}$ such that $y\in\GF(q)\backslash\{0\}$. The only ones which can be interchanged with $A$ by an element of $H$ are
$\{\infty,1\}$, by $t_{0,1,1,0}$ and $\{0,-1\}$, by $t_{-1,0,0,1}$,
when $q$ is odd. Thus the only edges between vertices of $O_0$ whose stabiliser
in $G$ is contained in $H$ are those in the orbits
$\{A,\{\infty,1\}\}^H$ and $\{A,\{0,-1\}\}^H$. The first gives the
matching 
$\{\{\{0,y\},\{\infty,y\}\}\mid y\in\GF(q)\backslash\{0\}\}$  and
hence the decomposition $\mathcal{P}_{\ominus}$ while the second
gives the matching 
$\{\{\{x,y\},\{x,-y\}\}\mid x\in\{0,\infty\},y\in\GF(q)\backslash\{0\}\}$ and hence Construction \ref{con:PGL3}(1). Both matchings have
$q-1$ edges and the second only occurs for $q$ odd. Note also that
both orbits are preserved  by $\PGammaL(2,q)_{\{0,\infty\}}$ and so
both decompositions are also $\PGammaL(2,q)$-decompositions.

Note that when $q$ is odd the orbit $O_{\frac{q-1}{2}}$ contains no
edges. Thus suppose next that $\{A,B\}\subseteq O_i$ for
$i<\frac{q-1}{2}$. Without loss of generality, let $A=\{1,\xi^i\}$. Then the 
neighbours of $A$ in $O_i$ are $\{1,\xi^{-i}\}$ and
$\{\xi^i,\xi^{2i}\}$ and these are interchanged by $H_A=\la
t_{0,\xi^i,1,0}\ra\cong C_2$. Hence $H$ acts transitively on the set
of edges between vertices of $O_i$. Moreover, $\la t_{0,1,1,0}\ra$ is
the stabiliser $H$ of the edge $\{\{1,\xi^i\},\{1,\xi^{-i}\}\}$ and so
$H$ contains the stabiliser in $G$ of an edge between two vertices of
$O_i$. Thus $\P$ is obtained by Construction
\ref{con:PGL3}(2). Moreover, an overgroup $\overline{G}=\la
\PGL(2,q),\phi^l\ra$ of $\PGL(2,q)$ in
$\PGammaL(2,q)$ preserves $\P$ if and only if
$\overline{G}_{\{0,\infty\}}=\la H,\phi^l\ra$ fixes $O_i$. Since
$\phi^l$ fixes $1$, it follows that $\phi^l$ fixes $O_i$ if and only
if $\phi^l$ fixes $\{\xi^i,\xi^{-1}\}$ and so $\overline{G}$ is as
stated in Construction \ref{con:PGL3}(2).
\Qed

\begin{construction}\label{con:MsqDq-1}
{\rm Let $G=M(s,q)$ and $\xi$ be a primitive element of $\GF(q)$ with
$q=p^f$ for some odd prime $p$ and even integer $f$. Let $i$ be an
integer and assume that either 
\begin{itemize}
\item $s=f/2$ and $(\xi^i)^{\la\phi^s\ra}$ has length 2 and does not
  contain $\xi^{-i}$, or
\item $s=f/4$ and $(\xi^i)^{\la\phi^s\ra}$ has length 4 and does
  contain $\xi^{-i}$.
\end{itemize}
Let $H=G_{\{0,\infty\}}=\la
\PSL(2,q)_{\{0,\infty\}},\phi^st_{\xi,0,0,1}\ra$ and note that
$\PSL(2,q)_{\{0,\infty\}}=\la t_{\xi^2,0,0,1},t_{0,1,1,0}\ra$.
\begin{enumerate}
\item Suppose that $i$ is even and let 
$e=\{\{1,\xi^i\},\{1,\xi^{-i}\}\}$ and $P=e^H$. Then
\begin{align*}
P =&\Big\{\big\{\{x^2,x^2\xi^i\},\{x^2,x^2\xi^{-i}\}\big\}\mid x\in\GF(q)\setminus\{0\}\Big\}\\
 &\cup \Big\{\big\{\{y,y\xi^{ip^s}\},\{y,y\xi^{-ip^s}\}\big\}\mid 
                      y= \not\square\Big\}
\end{align*}
Then $P$ has valency 2 (as the two neighbours of $\{1,\xi^i\}$ are
$\{1,\xi^{-i}\}$ and $\{\xi^i,\xi^{2i}\})$ and so is a union of
cycles. Each cycle has length the order of $\xi^i$ and so $P\cong 
(q-1,i)C_{\frac{q-1}{(q-1,i)}}$.

Now $|\{1,\xi^i\}^H|=q-1$ and by Lemma \ref{lem:Gedge},
$|G_e|=f/s$. Since $|H|=(q-1)f/s$ it follows that $|H_e|=f/s$ and so
$H_e=G_e$. Hence by Lemma \ref{lem:general} and the fact that $H$ is
maximal in $G$, letting $\P=P^G$ we get that $\P$ is a $G$-primitive
decomposition.

\item Suppose now that $i$ is odd and let
  $e=\{\{1,\xi^i\},\{1,\xi^{-i}\}\}$ and $P=e^H$. Then
\begin{align*}
P =&\Big\{\big\{\{x^2,x^2\xi^i\},\{x^2,x^2\xi^{-i}\}\big\}\mid
x\in\GF(q)\backslash\{0\}\Big\} \\
 &\cup \Big\{\big\{y,y\xi^{ip^s}\},\{y,y\xi^{-ip^s}\}\big\}
                     \mid y=\not\square\Big\}
\end{align*}
Then $|P|=q-1$ and so $|H_e|=f/s=|G_e|$, by Lemma \ref{lem:Gedge}. 
The only neighbour of $\{1,\xi^i\}$ in $P$ is $\{1,\xi^{-i}\}$  and so
$P=(q-1)K_2$.  By Lemma \ref{lem:general} and the fact that $H$ is
maximal in $G$, letting $\P=P^G$ we get that $\P$ is a $G$-primitive
decomposition.
\end{enumerate}
}%
\end{construction}

\begin{proposition}\label{lem:MsqD2(q-1)}
Let $(J(q+1,2),\mathcal{P})$ be a $G$-primitive decomposition with
$G=M(s,q)$ for some $s$ such that for $P\in\P$, $G_P=N_G(D_{q-1})$.
Then either $\P=\P_{\ominus}$,  or $\P$ arises from Construction $\ref{con:PGL3}(1)$, $\ref{con:PGL3}(2)$ or $\ref{con:MsqDq-1}$.
\end{proposition}
\proof
A subgroup $N_G(D_{q-1})$ of $G$ is a pair-stabiliser in $G$.
Without loss of generality we may suppose that $H=G_{\{0,\infty\}}=\la
\PSL(2,q)_{\{0,\infty\}},\phi^s t_{\xi,0,0,1}\ra$. Note that 
$q\equiv 1\pmod 4$ and so 
$\PSL(2,q)_{\{0,\infty\}}=\la t_{\xi^{2},0,0,1},t_{0,1,1,0}\ra$. Since
$G$ is 3-transitive it follows that
$$O_0=\{\{x,y\}\mid x\in\{0,\infty\},y\in\GF(q)\backslash\{0\}\}$$ is
an $H$-orbit on vertices and as in the proof of Lemma
\ref{lem:D2(q-1)}, if $\{A,B\}\subset O_0$ is an edge whose stabiliser
in $G$ is contained in $H$  we obtain either  $\P=\P_{\ominus}$ or
$\P$ is obtained by Construction \ref{con:PGL3}(1).

Now suppose $\{A,B\}\not\subset O_0$. Since $H$ is transitive on
$\GF(q))\backslash\{0\}$, we can assume that $A=\{1,\xi^i\}$ where
$1\leq i\leq q-2$ and that $A\cap B=\{1\}$, say $B=\{1,t\}$. We need to find the neighbours $B$ of $A$ such that
$G_{\{A,B\}}\leqslant H$. Let $g\in \PGL(2,q)$ map
$\{\{\infty,0\},\{\infty,1\}\}$ onto $\{A,B\}$. Then $G_{\{A,B\}}=\la
t_{-1,1,0,1},\phi^{2s}\ra^g$ by Lemma \ref{lem:Gedge}. Hence
$t_{-1,1,0,1}$ and $\phi^{2s}$ must stabilise $\{0,\infty\}^{g^{-1}}$. Note that $\infty^g\neq \infty$ (since $\infty\notin A$) and $\infty^g\neq 0$ (since $O\notin A$).

Since $B=\{1,t\}$, we can take $g=t_{a,\xi^i,a,1}$ where
$a=\frac{\xi^i-t}{t-1}$, and then
$\{0,\infty\}^{g^{-1}}=\{-\frac{\xi^i}{a},-\frac{1}{a}\}$. Recall that $t_{-1,1,0,1}$ stabilises this set. Now $t_{-1,1,0,1}$ fixes only the points $\infty, 2^{-1}$, and if $\{0,\infty\}^{g^{-1}}=\{\infty,2^{-1}\}$ we would have \mbox{$\infty^g\in\{0,\infty\}$} which is not the case. Hence $t_{-1,1,0,1}$ interchanges $-\frac{\xi^i}{a}$ and $-\frac{1}{a}$. Thus $-\frac{\xi^i}{a}=1+\frac{1}{a}$, that is
$a=-1-\xi^i=\frac{\xi^i-t}{t-1}$, and so $t=\xi^{-i}$.  For this value of $t$, $\{0,\infty\}^{g^{-1}}=\{\frac{\xi^i}{\xi^i+1},\frac{1}{\xi^i+1}\}$ and this set is stabilised by $t_{-1,1,0,1}$ and $\phi^{2s}$. The equality
$\{\frac{\xi^i}{1+\xi^i},\frac{1}{1+\xi^i}\}^{\phi^{2s}}=\{\frac{\xi^i}{1+\xi^i},\frac{1}{1+\xi^i}\}$, is equivalent to either $\frac{\xi^{ip^{2s}}}{1+\xi^{ip^{2s}}}=\frac{\xi^i}{1+\xi^i}$ and $\frac{1}{1+\xi^{ip^{2s}}}=\frac{1}{1+\xi^i}$, or $\frac{\xi^{ip^{2s}}}{1+\xi^{ip^{2s}}}=\frac{1}{1+\xi^i}$ and $\frac{1}{1+\xi^{ip^{2s}}}=\frac{\xi^i}{1+\xi^i}$. In the first case $\xi^{ip^{2s}}=\xi^i$; in the second case $\xi^{ip^{2s}}=\xi^{-i}$.
That means $O=(\xi^i)^{\la \phi^s\ra}$ has length 1,2 or 4.  

Set $e=\{A,\{1,\xi^{-1}\}\}$. If $O$ has length 1, or $O$ has length 2 and $(\xi^i)^{\phi^s}=\xi^{-i}$, then $e^H$ yields a decomposition in Construction \ref{con:PGL3}(2). If $O$ has length 2 and $(\xi^i)^{\phi^s}\neq \xi^{-i}$, or $O$ has length 4 and $\xi^{ip^{2s}}=\xi^{-i}$, then $e^H$ yields a decomposition Construction \ref{con:MsqDq-1}(1) if $i$ is even and in Construction \ref{con:MsqDq-1}(2) if $i$ is odd.
\Qed

\subsection{$D_{q+1}$ subgroups}
Before dealing with the case where $H\cap \PSL(2,q)=D_{q+1}$ we need a
new model for the group action. Let $K=\GF(q^2)$ for $q=p^f$ with primitive element $\xi$, and let
$F=\{0\}\cup\{(\xi^{q+1})^l\mid l=0,1,\ldots,q-2\}\cong\GF(q)$. 
The element $\xi$ acts on $K$ by multiplication and induces an $F$-linear map. Moreover, under the induced action of$F$, $K$ is a 2-dimensional vector space over $F$. The field automorphism $\varphi$ of $K$ of order $2f$ mapping each element of $K$ to its $p^{\mathrm{th}}$ power is $F$-semilinear, that is, $\varphi$ preserves addition and for each $x\in K$, $\lambda\in F$, we have
$(\lambda x)^{\varphi}=\lambda^p x^{\varphi}$.  Then $\GammaL(2,q)=\la\GL(2,q),\varphi\ra$. Note that $\varphi^f$ is an $F$-linear map so $\varphi^f\in\GL(2,q)$.

We can identify the projective line $X$ on which $\PGL(2,q)$ acts with
the elements of $K$ modulo $F$, that is, $X=\{\xi^iF\mid i=0,1,\ldots,q\}$.  Then 
$\PGammaL(2,q)=\la \PGL(2,q),\varphi\ra$. Multiplication by $\xi$ 
induces the map $\hat\xi$ of order $q+1$ and $\la \hat{\xi}\ra$ is
normalised by $\varphi$. Moreover, for each $i$, 
$(\xi^i F)^{\varphi^f}=\xi^{iq}F=\xi^{-i}F$ and so $\varphi^f$ inverts
$\hat\xi$. Hence $\la \hat{\xi},\varphi^f\ra\cong D_{2(q+1)}$.

\begin{construction}\label{con:PGL4}
{\rm Let $X$ be the projective line modelled as above.
Let  $1\leq i<\frac{q+1}{2}$ and $e=\{\{1F,\xi^iF\},\{1F,\xi^{-i}F\}\}$ and let
$s$ be a positive integer dividing $f$ such that $\la \varphi^s\ra$ has
$\{\xi^iF,\xi^{-i}F\}$ as an orbit on $X$. Let 
$G=\la \PGL(2,q),\varphi^s\ra$ and 
$H=\la \hat{\xi},\varphi^s\ra\cong C_{q+1}\rtimes C_{2f/s}$. Now
$\la\varphi^s\ra$ fixes $e$ and has order $2f/s$, which by Lemma
\ref{lem:Gedge} is the order of $G_e$. Hence $G_e<H$ and $H$ is a
maximal subgroup of $G$. Thus by Lemma \ref{lem:general}, letting 
$$P=e^H=\Big\{\big\{\{xF,x\xi^iF\},\{xF,x\xi^{-i}F\}\big\}\mid
x\in\GF(q)\backslash\{0\}\Big\}$$ and $\P=P^G$, we obtain a
$G$-primitive decomposition of $J(q+1,2)$.  The divisors have valency 2
and hence are unions of cycles. These cycles have length the order
of $\xi^iF$, which is $\frac{q+1}{(q+1,i)}$. Thus each divisor is
isomorphic to $(q+1,i)C_{\frac{q+1}{(q+1,i)}}$.
}%
\end{construction}

\begin{proposition}\label{lem:D2(q+1)}
Let $(J(q+1,2),\mathcal{P})$ be a $G$-primitive decomposition such that
$\PGL(2,q)\leqslant G \leqslant\PGammaL(2,q)$ and, for $P\in\P$, $G_P=N_G(D_{2(q+1)})$. Then $\P$ is obtained by Construction $\ref{con:PGL4}$.
\end{proposition}
\proof
Since $\PGammaL(2,q)=\la\PGL(2,q),\varphi\ra$ and $\varphi^f\in\PGL(2,q)$
we have $G=\la\PGL(2,q),\varphi^s\ra$ for some $s$ dividing $f$. Let
$L=\la \hat{\xi},\varphi^f\ra\cong D_{2(q+1)}$. Then 
$N_G(L)=\la \hat{\xi},\varphi^s\ra\cong C_{q+1}\rtimes C_{2f/s}$ and
we may assume that $H=G_P=N_G(L)$. Let $e\in P$. Since $H$ is
transitive on $X$ we may also assume that
$e=\{\{1F,\xi^iF\},\{1F,\xi^jF\}\}$ for some integers $i$ and
$j$. Since $H_{1F}=\la \varphi^s\ra$ and by Lemma \ref{lem:Gedge},
$|G_e|=2f/s$, it follows that $G_e\leqslant H$ if and only if
$\la\varphi^s\ra$ has $\{\xi^iF,\xi^jF\}$ as an orbit on $X$. Since
$\varphi^f\in\la\varphi^s\ra$ and maps $\xi^iF$ to $\xi^{-i}F$ it
follows that $j=-i$. Since $\xi^{-i}F=\xi^{q+1-i}F$ we may assume that $1\leq i\leq (q+1)/2$. Moreover, if $i=(q+1)/2$ then $q$ is odd and
$\xi^{-(q+1)/2}F=\xi^{(q+1)/2}F$. Thus we may further assume that
$1\leq i<(q+1)/2$. Hence $\P$ arises from Construction \ref{con:PGL4}.
\Qed

Next we need the following lemma about the normaliser in $M(s,q)$ of a subgroup  $D_{q+1}$ in $\PSL(2,q)$.
\begin{lemma}
\label{lem:msqq+1dihedral}
Suppose $q=p^f$ where $f$ is even and $p$ is an odd prime. Let $L=\la
\hat{\xi},\varphi^f\ra\cap \PSL(2,q)$ and $G=M(s,q)$ for some divisor
$s$ of $f/2$. Then
\begin{enumerate}
\item $L=\la \hat{\xi}^2,\varphi^f\ra\cong D_{q+1}$.
\item If $p\equiv 1\pmod 4$ or $s$ is even then 
$N_G(L)= \la \hat{\xi}^2,\varphi^s\hat{\xi}\ra$, and is transitive
  on the projective line.
\item If $p\equiv 3\pmod 4$ and $s$ is odd then
$N_G(L)=\la \hat{\xi}^2,\varphi^s \ra$, and has two equal
  sized orbits on the projective line.
\end{enumerate}
\end{lemma}
\proof
Now $\{1,\xi^{(q+1)/2}\}$ is a basis for $K$ over $F$. Define $\phi:K\rightarrow K$ such that, for all $\lambda_1,\lambda_2\in F$,
$(\lambda_1+\lambda_2\xi^{(q+1)/2})^{\phi}=\lambda_1^p+\lambda_2^p\xi^{(q+1)/2}$.
Then $\GammaL(2,q)=\la \GL(2,q),\phi\ra$. Since also $\GammaL(2,q)=\la\GL(2,q),\varphi\ra$, we must have $\varphi=\phi g$ for
some $g\in \GL(2,q)$. Since $\varphi$ and $\phi$ fix 1, so does
$g$. Moreover, $\phi$ fixes $\xi^{(q+1)/2}$ while
$(\xi^{(q+1)/2})^{\varphi}=\xi^{p(q+1)/2}=\xi^{\frac{(p-1)(q+1)}{2}}\xi^{\frac{q+1}{2}}$. Note that
$\xi^{\frac{(p-1)(q+1)}{2}}\in F$ and so $\xi^{(q+1)/2}$ is an
eigenvector for $g$. Thus with respect to the basis
$\{1,\xi^{(q+1)/2}\}$, the element $g$ is represented by the matrix
$$\left(\begin{array}{cc} 1 & 0 \\ 
0 & \xi^{\frac{(p-1)(q+1)}{2}}\end{array}\right),$$
and $\det(g)=\xi^{\frac{(p-1)(q+1)}{2}}$ is a square in $\GF(q)$ if and only if $p\equiv 1\pmod 4$. Furthermore, $\varphi^f$ is represented by the matrix 
$$\left(\begin{array}{cc} 1 & 0 \\ 
0 & -1\end{array}\right).$$

Recall that an element of $\GL(2,q)$ induces an element of $\PSL(2,q)$ if and only if its determinant is a $\GF(q)$-square. Since $q\equiv 1\pmod 4$ it follows that $\varphi^f\in\PSL(2,q)$. Now $\la \hat{\xi^2}\ra\cong C_{(q+1)/2}$ and $\hat{\xi}^2\in\PSL(2,q)$, and since $\varphi^f$ inverts $\hat{\xi}$ it also inverts $\hat{\xi}^2$. Hence $L$ is as in part (1) of the lemma. Moreover, $L$ has two orbits on the projective line $X$, these being $\{1F,\xi^2F,\ldots,\xi^{q-1}F\}$ and
$\{\xi F,\xi^3F,\ldots,\xi^qF\}$.

Now $\varphi=\phi g$ and $g\in\PSL(2,q)$ if and only if $p\equiv 1\pmod 4$. By definition it follows that $G=M(s,q)=\la \PSL(2,q),\phi^st\ra$ for any $t\in\PGL(2,q)\setminus\PSL(2,q)$. Suppose first that $p\equiv 1\pmod 4$. Then $\varphi=\phi g$ with $g\in\PSL(2,q)$ and so
$G=\la \PSL(2,q),\varphi^s\hat{\xi}\ra$. When $p\equiv 3\pmod 4$ we have $\varphi=\phi g$ with $g\in\PGL(2,q)\setminus \PSL(2,q)$. Thus for odd $s$ we have $G=\la \PSL(2,q),\varphi^s\ra$ while for even $s$ we have
$G=\la \PSL(2,q),\varphi^s\hat{\xi}\ra$. Now $(\varphi^f)^{\varphi^s\hat{\xi}}=(\varphi^f)^{\hat{\xi}}=\varphi^f\hat{\xi}^{-p^s+1}\in L$. Hence for $p\equiv 1\pmod 4$ or $s$ even we have $N_G(L)=\la \hat{\xi}^2,\varphi^s\hat{\xi}\ra$. Since $\varphi^s\hat{\xi}$ interchanges the two $L$-orbits on $X$, $N_G(L)$ is transitive on $X$ and so we have proved part (2). For $p\equiv 3\pmod 4$ and $s$ odd we have
$N_G(L)=\la \hat{\xi}^2,\varphi^s \ra$.   Since $\varphi^s$ fixes each $L$-orbit it follows that $N_G(L)$ has two orbits and the proof is complete.
\Qed

\begin{construction}
\label{con:MsqDq+11}
{\rm Let $q=p^f$ where $p$ is odd and $f$ even and let $G=M(s,q)$ for
some divisor $s$ of $f/2$. Suppose that either $p\equiv 1\pmod 4$ or
$s$ is even. Let $1\leq i<(q+1)/2$ such that $\la \varphi^{2s}\ra$ has
$\{\xi^iF,\xi^{-i}F\}$ as an orbit on $X$. Let
$H=\la \hat{\xi}^2,\varphi^s\hat{\xi}\ra$ and
$e=\{\{1F,\xi^iF\},\{1F,\xi^{-i}F\}\}$. Now $\la \varphi^{2s}\ra$  fixes
$e$, lies in $G$, and has order $f/s$. Since this is the same order as $G_e$ (Lemma \ref{lem:Gedge}) it follows that $G_e<H$. Hence by Lemma \ref{lem:general}, letting $P=e^H$ and $\P=P^G$ we obtain a $G$-primitive decomposition.
\begin{enumerate}
\item Suppose first that $i$ is even. Then $H_{\{1F,\xi^iF\}}=\la
  \varphi^f\hat{\xi}^i,\varphi^{4s}\ra$ whose orbit containing
  $\{1F,\xi^{-i}F\}$ is $\{\{1F,\xi^{-i}F\},\{\xi^iF,\xi^{2i}F\}\}$.  Thus
  $P$ has valency 2 and so is a union of cycles of length the order of
  $\hat{\xi}^i$, that is, $P\cong (q+1,i)C_{\frac{q+1}{(q+1,i)}}$.

\item Suppose now that $i$ is odd. An element of $H$ mapping $1F$ to $\xi^iF$ is of the form $h=\varphi^{st}\hat{\xi}^i$ with $t$ odd. Since $\la \varphi^{2s}\ra$ has $\{\xi^iF,\xi^{-i}F\}$ as an orbit on $X$, we have that $h$ maps  $\xi^iF$ onto $\xi^{i(1+p^s)}F$ or onto  $\xi^{i(1-p^s)}F$, according as $t\equiv 1$ or $3\pmod 4$ respectively. Hence, for $h$ to map  $\xi^iF$ onto $1F$, we need $q+1$ to divide $i(1+p^s)$ or $i(1-p^s)$ respectively. Since $p^{2s}-1$ divides $p^f-1=q-1$, it follows that $\gcd(q+1,p^s+1)=2$ and $\gcd(q+1,p^s-1)=2$, and so $\frac{q+1}{2}$ must divide $i$, which is not possible since $1\leq i < \frac{q+1}{2}$. Thus $(\xi^iF)^h\neq 1F$.  Hence $H_{\{1F,\xi^iF\}}=H_{1F,\xi^iF}=\la\varphi^{4s}\ra$, which also fixes $\xi^{-i}F$ and hence fixes $e$. Thus $P$ is a matching with $q+1$ edges.
\end{enumerate}
}%
\end{construction}

\begin{construction}
\label{con:MsqDq+12}
{\rm Let $p\equiv 3\pmod 4$ and let $G=M(s,q)$ for $q=p^f$ and $s$ an odd
divisor of $f/2$.  Let  $1\leq i< (q+1)/2$ such that
$\la\varphi^{2s}\ra$ has $\{\xi^iF,\xi^{-i}F\}$ as an orbit on $X$.
Let $H= \la \hat{\xi}^2,\varphi^s \ra$ and
$e=\{\{1F,\xi^iF\},\{1F,\xi^{-i}F\}\}$. Now $\la \varphi^{2s}\ra$  fixes
$e$, lies in $G$ and has order $f/s$. Since this is the same order as $G_e$ (Lemma
\ref{lem:Gedge}) it follows that $G_e<H$ and so by Lemma
\ref{lem:general}, letting $P=e^H$ and $\P=P^G$, we obtain a
$G$-primitive decomposition.
\begin{enumerate}
\item Suppose first that $i$ is even. Then 
$H_{\{1F,\xi^iF\}}=\la \varphi^f\hat{\xi}^i, \varphi^{4s}\ra$ and the $H$-orbit containing $\{1F,\xi^{-i}F\}$ has length 2. Thus $P$ is a union of cycles of length the order of $\hat{\xi}^i$, so $P\cong (q+1,i)C_{\frac{q+1}{(q+1,i)}}$.

\item If $i$ is odd then $1F$ and $\xi^iF$ lie in different $H$-orbits
and so $H_{\{1F,\xi^iF\}}=H_{1F,\xi^iF}=\la\varphi^{4s}\ra$ which also fixes $\xi^{-i}F$ and hence fixes $e$. Thus $P$ is a matching with $q+1$ edges.
\end{enumerate}
}%
\end{construction}

\begin{construction}
\label{con:MsqDq+13}
{\rm Let $p\equiv 3\pmod 4$ and let $G=M(s,q)$ for $q=p^f$ and $s$ an odd
divisor of $f/2$.  Let $1\leq i<\frac{q+1}{2}$ such that
$\la \hat{\xi}^{-1}\varphi^{2s}\hat{\xi}\ra$ has $\{\xi^{i+1}F,\xi^{-i+1}F\}$ as an orbit on $X$. Let $H= \la \hat{\xi}^2,\varphi^s \ra$ and
$e=\{\{\xi F,\xi^{i+1}F\},\{\xi F,\xi^{-i+1}F\}\}$. Now
$\la \hat{\xi}^{-1}\varphi^{2s}\hat{\xi}\ra\leqslant H$, fixes $e$, and has the
same order as $G_e$.  Thus $G_e<H$ and so by  Lemma
\ref{lem:general}, letting $P=e^H$ and $\P=P^G$, we obtain a
$G$-primitive decomposition.
\begin{enumerate}
\item Suppose first that $i$ is odd. Then $\xi F$ and $\xi^{i+1}F$ lie
  in different $H$-orbits. Hence 
$H_{\{\xi F,\xi^{i+1}F\}}=H_{\xi F,\xi^{i+1}F}=\la \hat{\xi}^{-1}\varphi^{4s}\hat{\xi}\ra$  which also fixes
  $\xi^{-i+1}F$ and so $P$ is a matching with $q+1$ edges. 

\item If $i$ is even then $\varphi^f\hat{\xi}^{i+2}\in H$ interchanges
  $\xi F$ and $\xi^{i+1} F$, and so $H_{\{\xi F,\xi^{i+1}F\}}$ $=\la
  \hat{\xi}^{-1}\varphi^{4s}\hat{\xi},\varphi^f\hat{\xi}^{i+2}\ra$,
  whose orbit containing $\{\xi F,\xi^{-i+1} F\}$ has size 2. 
Hence $P$ is a union of cycles of length the
  order of $\hat{\xi}^i$. Thus
  $P=(q+1,i)C_{\frac{q+1}{(q+1,i)}}$.
\end{enumerate}
}%
\end{construction}

\begin{proposition}
Let $\P$ be an $M(s,q)$-primitive decomposition of $J(q+1,2)$ with
divisor stabiliser $N_{M(s,q)}(D_{q+1})$. Then
$\P$ can be obtained from Construction $\ref{con:MsqDq+11}$,
$\ref{con:MsqDq+12}$ or $\ref{con:MsqDq+13}$.
\end{proposition}
\proof
Let $G=M(s,q)$ and suppose first that $q=p^f$ where $p\equiv 1\pmod 4$
or $s$ is even. We may assume that 
$H=\la \hat{\xi}^2,\varphi^s\hat{\xi}\ra$ by Lemma \ref{lem:msqq+1dihedral}. Let $e\in P\in\P$. By Lemma \ref{lem:msqq+1dihedral} again, $H$ is transitive on $X$ and so we can
assume that $e=\{\{1F,\xi^iF\},\{1F,\xi^jF\}\}$ for some $i$ and
$j$. Now $H_{1F}=\la \varphi^{2s}\ra$, which has order $f/s$. By Lemma
\ref{lem:Gedge}, this is the same order as $G_e$. Hence  $G_e<H$ if
and only if $H_{1F}=G_e$, which holds if and only if
$\{\xi^iF,\xi^jF\}$ is an orbit of $\la \varphi^{2s}\ra$. Since
$\varphi^f\in\la  \varphi^{2s}\ra$ and maps $\xi^iF$ to $\xi^{-i}F$ it
follows that $j=-i$ and we may assume as before that $1\leq
i<(q+1)/2$. Thus $\P$ comes from Construction \ref{con:MsqDq+11}.  

Suppose now that $p\equiv 3\pmod 4$ and $s$ is odd.  Then by Lemma
\ref{lem:msqq+1dihedral}, we may
assume that $H=\la \hat{\xi}^2,\varphi^s\ra$. Let $e\in P\in\P$. By
Lemma \ref{lem:msqq+1dihedral}, $H$ has 2 orbits on $X$ and so we may
assume that $e=\{\{1F,\xi^iF\},\{1F,\xi^jF\}\}$ or
$\{\{\xi F,\xi^{i+1}F\},\{\xi F,\xi^{j+1}F\}\}$. Suppose that $e$ is the
first edge. Now $H_{1F}=\la \varphi^s\ra$ which has order $2f/s$ while $G_e$
has order $f/s$ by Lemma \ref{lem:Gedge}. Since $H_{1F}$ has a unique
subgroup of order $f/s$ it follows that $G_e<H$ if and only if
$G_e=\la \varphi^{2s}\ra$, that is, if and only if $\la
\varphi^{2s}\ra$ has $\{\xi^i F,\xi^j F\}$ as an orbit on $X$. Since
$\varphi^f\in\la \varphi^{2s}\ra$ we have $j=-i$ and may assume 
$1\leq i<(q+1)/2$. It
follows that $\P$ is as constructed in Construction
\ref{con:MsqDq+12}. If on the other hand 
$e=\{\{\xi F,\xi^{i+1}F\},\{\xi F,\xi^{j+1}F\}\}$, then 
$H_{\xi F}=\la \hat{xi}^{-1}\varphi^s\hat{\xi}\ra$ which has order
$2f/s$.   Its only index two subgroup is 
$\la \hat{\xi}^{-1}\varphi^{2s}\hat{\xi}\ra$ and so by order arguments
again this must have $\{\xi^{i+1} F,\xi^{j+1} F\}$ as an orbit. Since
$\hat{\xi}^{-1}\varphi^f\hat{\xi}\in \la
\hat{\xi}^{-1}\varphi^{2s}\hat{\xi}\ra$ and maps $\xi^{i+1}F$ to 
$\xi^{-i+1}F$ it follows that $j=-i$. Once again we have 
$1\leq i<\frac{q+1}{2}$. Hence $\P$ is as given by Construction
\ref{con:MsqDq+13}.
\Qed

\subsection{$S_4$-subgroups}
First we have the following lemma on the orbit lengths of a
subgroup $S_4$ of $\PGL(2,q)$ which we have adapted from \cite{PGLsubs}.
\begin{lemma}\label{lem:S4orbits}{\rm \cite[Lemma 10]{PGLsubs}}
Let  $q=p\equiv\pm 3\pmod 8$, $q>3$, $G=\PGL(2,q)$ acting on the
projective line $X$, and $H$ a subgroup of $G$ isomorphic to $S_4$.
Then $H$ has the following orbits of length less than $24$ on $X$.
\begin{enumerate}
\item If $q\equiv 5\pmod {24}$, then $H$ has one orbit of length $6$.
\item If $q\equiv 11\pmod {24}$, then $H$ has one orbit of length $12$.
\item If $q\equiv 13\pmod {24}$, then $H$ has one orbit of length $6$
  and one of length $8$.
\item If $q\equiv 19\pmod {24}$, then $H$ has one orbit of length $8$
  and one of length $12$.
\end{enumerate}
\end{lemma}

\begin{construction}\label{con:PGL5}
{\rm Let $X=\{\infty\}\cup \GF(q)$ be the projective line.
\begin{enumerate}
\item Let $q\equiv\pm 3\pmod 8$ be a prime ($q>3$) and $H=S_4$.
  Choose $x,y_1,y_2\in X$ such that $(|x^H|,|y_1|^H)=(6,8), (6,24),
  (12,8)$ or $(12,24)$, and there exists in $H_x$ an element switching
  $y_1$ and $y_2$.  Let $P=\{\{\{x,y_1\},\{x,y_2\}\}^H$ and $\P=P^{\PGL(2,q)}$. Then  by Lemma \ref{lem:general}, $(J(q+1,2),\P)$ is a $\PGL(2,q)$-primitive
  decomposition. Since \mbox{$|\{x,y_1\}|^H=24$}, the stabiliser in $H$ of
  $\{x,y_1\}$ is trivial. Hence the divisors are isomorphic to $12K_2$. 
\item Let $q\equiv 5\pmod 8$ be a prime and $H=S_4$. Let
  $P=\{\{x,y_1\},\{x,y_2\}\}^H$ where $x,y_1,y_2$ all lie in an
  $H$-orbit of length 6 and there exists in $H_x$ an element switching
  $y_1$ and $y_2$. By Lemma \ref{lem:S4orbits}, there is a unique
  orbit of $O_6$ of length 6. The group $H$ acts imprimitively on
  $O_6$ with blocks of size $2$, and $H_x\cong C_4$ contains an
  element interchanging $y_1$, $y_2$ if and only if $\{y_1,y_2\}$ is a
  block not containing $x$. Moreover, $P\cong 3C_4$. Let
  $\P=P^{\PGL(2,q)}$. Then by Lemma \ref{lem:general} $(J(q+1,2),\P)$,
  is a $\PGL(2,q)$-primitive decomposition. 

\item Let $q\equiv 3\pmod 8$ be a prime and $H=S_4$. Let 
$P=\{\{x,y_1\},\{x,y_2\}\}^H$ where $x,y_1,y_2$ all lie in an
  $H$-orbit of length 12 and
and there exists in $H_x$ an element switching $y_1$ and $y_2$.  By
Lemma \ref{lem:S4orbits},  there is a unique orbit $O_{12}$ of length
12. We can see this action as $S_4$ acting on ordered pairs, denoted
by $[a,b]$.  Then for $x=[1,2]\in O_{12}$, $H_x$ is the transposition
$(3,4)$ in $S_4$. It fixes one 
remaining point of $O_{12}$, namely $[2,1]$ and interchanges the 5 pairs 
$\{[2,3],[2,4]\}$, $\{[3,1],[4,1]\}$, $\{[1,3],[1,4]\}$,
$\{[3,2],[4,2]\}$, and $\{[3,4],[4,3]\}$. If we take $\{y_1,y_2\}$ as in the first two cases, then the
stabiliser in $H$ of $\{x,y_1\}$ is trivial and so we get a matching
$12K_2$ in each case. In the last three cases, the stabiliser in $H$
of $\{x,y_1\}$ has order 2, and we get unions of cycles. It is easy to
see that in the third and fourth case, we get $4C_3$, while in the
last case we get $3C_4$. Let $\P=P^{\PGL(2,q)}$. Then by Lemma
\ref{lem:general}, $(J(q+1,2),\P)$ is a $\PGL(2,q)$-primitive
decomposition.
\end{enumerate}
}%
\end{construction}

\begin{proposition}\label{lem:S4}
Let $(J(q+1,2),\mathcal{P})$ be a $G$-primitive decomposition with
$G=\PGL(2,q)$ for $q=p\equiv\pm 3\pmod 8$ with $q\geq 5$ and given
$P\in\P$ we have $G_P\cong S_4$.  Then $P$ is obtained by Construction
$\ref{con:PGL5}(1)$, $(2)$ or $(3)$.
\end{proposition}
\proof
Let $P\in \P$ and $H=G_P\cong S_4$. If $\{x,y\}\subseteq X$ with $x$ and $y$
in different $H$-orbits of length 24 then $|\{x,y\}^H|=24$ and that orbit
contains no edges of $J(q+1,2)$. Thus if $x$ and $y$ come from
different $H$-orbits $O_1$ and $O_2$ respectively, we may assume by
Lemma \ref{lem:S4orbits},
that $|O_1|<|O_2|$ and so $\{x,y\}^H$ has length $\lcm(|O_1|,|O_2|)$
and contains edges. Moreover, $H$ contains the stabiliser in $G$ of
such an edge $\{\{x,y_1\},\{x,y_2\}\}$ if and only if $H_x$ contains
an element interchanging $y_1$ and $y_2$. If $x$ is in an orbit of
size 8 then $|H_x|=3$ and so no such element exists, and if $x$ is in
an orbit of size 24 then $|H_x|=1$ 
and so no such element exists.  Thus the possibilities for
$(|O_1|,|O_2|)$ are $(6,8)$, $(6,24)$, $(8,12)$ or
$(12,24)$. In the first two cases $x$ must be in the orbit of length 6 
and in the last two cases $x$ must be in the orbit of length 12. Thus
we get the decomposition of Construction \ref{con:PGL5}(1). 

Suppose now  $e=\{\{x,y_1\},\{x,y_2\}\}$ is an
edge such that $x,y_1,y_2\}$ lie in the same $H$-orbit $O_i$. Then $H$
contains $G_e$ if and only if $H_x$ interchanges $y_1$ and $Y_2$. Thus
$|H_x|$ is even and so $|O_i|\neq 8,24$. If $q\equiv 5\pmod 8$ and
$O_i$ is the unique orbit of size 6 then we obtain the decomposition
in Construction \ref{con:PGL5}(2). If $q\equiv 3\pmod 8$ and $O_i$ is
the unique orbit of size 12 then we obtain the decompositions in
Construction \ref{con:PGL5}(3).
\Qed

\subsection{Subfield subgroups}
Suppose now that $q=q_0^r$. Then $S=\{\infty\}\cup \GF(q_0)$ is a
subset of the projective line $X=\{\infty\}\cup \GF(q)$ which is an
orbit of the subgroup $\PGammaL(2,q_0)$ of $\PGammaL(2,q)$.
Notice that $\phi$ fixes the set $S$.
 Moreover, by
\cite[I, Example 3.23]{handbook}, if $\mathcal{B}=S^{\PGL(2,q)}$ then
$(X,\mathcal{B})$ is a $S(3,q_0+1,q+1)$ Steiner system. Since $\phi$
fixes $S$ and $\PGammaL(2,q)=\la\PGL(2,q),\phi\ra$ it follows that
$\mathcal{B}=S^{\PGammaL(2,q)}$. Thus by Lemma
\ref{lem:design}, we can construct a $\PGammaL(2,q)$-transitive
decomposition of $J(q+1,2)$ with divisors isomorphic to
$J(q_0+1,2)$. The stabiliser of a divisor is $\PGammaL(2,q_0)$. Moreover,
this decomposition is $G$-transitive for any 3-transitive subgroup $G$  of
$\PGammaL(2,q)$.  For further constructions we need the orbits of
$\PGL(2,q_0)$ on $\GF(q)\setminus \GF(q_0)$.

\begin{lemma}
\label{lem:subfieldorbits}
{\rm \cite[Lemma 14]{PGLsubs}}
Let $q=q_0^r$ for some prime $r$ and let
$H=\{t_{a,b,c,d}\mid a,b,c,d\in\GF(q_0),ad-bc\neq 0\}$. If $r$ is odd
then $H$ acts semiregularly on $\GF(q)\setminus\GF(q_0)$, while if
$r=2$ then $H$ is transitive on $\GF(q)\setminus \GF(q_0)$.
\end{lemma}

\begin{construction}
\label{con:PGL6}
{\rm Let $X=\{\infty\}\cup \GF(q)$ be the projective line. Let $q=q_0^r$, where $q_0>2$, $r$ is a prime and $r$ is odd if $q$ is odd.
Let $e=\{\{\infty,w_1\},\{\infty,w_2\}\}$ such
that $w_1,w_2\in \GF(q)\setminus\GF(q_0)$ but $w_1+w_2\in\GF(q_0)$.
Let $l$ be a positive integer such that $\phi^l$ fixes
$\{w_1,w_2\}$. Then let $G=\la \PGL(2,q),\phi^l\ra$ and
$H=\la\PGL(2,q_0),\phi^l\ra$. Let $P=e^H$ and
$\P=P^G$. Then by Lemma \ref{lem:Gedge},
$G_e=\la t_{-1,w_1+w_2,0,1},\phi^l\ra$ which is in $H$. Therefore by
Lemma \ref{lem:general}, $(J(q+1,2),\P)$ is a $G$-primitive
decomposition.  The stabiliser $H_{\{\infty,w_1\}}$ fixes $\infty$ and
$w_1$ as they are in different $H$-orbits. We claim that
$\PGL(2,q_0)_{\infty,w_1}=1$. Indeed, an element in that subgroup must
be of the form $t_{a,b,0,1}$ with $a,b\in \GF(q_0)$, whose only fixed point is $\frac{b}{1-a}\in\GF(q_0)$ if it is not the
identity. Hence  there is a unique element of $\PGL(2,q_0)_{\infty}$
interchanging $w_1$ and $w_2$, this being $ t_{-1,w_1+w_2,0,1}$. Then
as $\phi^l$ fixes $\{w_1,w_2\}$ and $\infty$, it 
follows that $H_{\infty,w_1}$ fixes $w_2$. Hence $P$ is isomorphic to
$\frac{q_0(q_0^{2}-1)}{2}K_2$.
}%
\end{construction}

\begin{proposition}\label{lem:PGLq0}
Let $(J(q+1,2),\mathcal{P})$ be a $G$-primitive decomposition such that
$\PGL(2,q)\leqslant G\leqslant \PGammaL(2,q)$ and for $P\in\P$, $G_P\cong N_G(\PGL(2,q_0))$ where $q=q_0^r$, $q_0>2$, $r$ is a prime and $r$ is odd if $q$ is odd. Then  $\P$ is obtained by Construction $\ref{con:design}$ or $\ref{con:PGL6}$.
\end{proposition}
\proof
By Theorem \ref{thm:dropdown}, $\mathcal{P}$ is also a
$\PGL(2,q)$-primitive decomposition so we may suppose that $G=\PGL(2,q)$ and $H=G_P=\{t_{a,b,c,d}\mid a,b,c,d\in\GF(q_0), ad-bc\neq 0\}$.
We have already seen that $H$ has an orbit $\{\infty\}\cup \GF(q_0)$
of length $q_0+1$ on $X$.  Moreover, by Lemma \ref{lem:subfieldorbits}, when $r$ is odd, $H$ has $q_0^{r-3}+q_0^{r-5}+\cdots+q_0^{2}+1$ other orbits, all of length $q_0(q_0^{2}-1)$, while when $r=2$, $H$ is transitive on $\GF(q)\backslash\GF(q_0)$.

Suppose that $H$ contains the stabiliser in $G$ of the edge
$e=\{\{v,w_1\},\{v,w_2\}\}$. Then $H_v$ contains the unique nontrivial element
interchanging $w_1$ and $w_2$ (see Lemma \ref{lem:Gedge}). Now $v$ must lie in the unique $H$-orbit of length $q_0+1$. For, if $r$ is odd and $v$ lies in an $H$-orbit of length $q_0(q_0^{2}-1)$ then $H_v=1$, while if $r=2$ and $v$ lies in $\GF(q)\setminus \GF(q_0)$, then $|H_v|=q_0+1$ which is odd. Without loss of generality we may suppose that $v=\infty$.

Then $G_e=\la t_{-1,w_1+w_2,0,1}\ra$, so $G_e\leqslant H$ if and only if $w_1+w_2\in \GF(q_0)$. If $w_1$ and $w_2$ lie in the orbit of length $q_0+1$, that is, are in $\GF(q_0)$ then we obtain the decomposition from Construction \ref{con:design}, which is in fact preserved by $\PGammaL(2,q)$. If $w_1\notin\GF(q_0)$ and $w_2=a-w_1$ with $a\in\GF(q_0)$, then we get a decomposition obtained from Construction \ref{con:PGL6}.
\Qed

For a primitive element $\mu$ of $\GF(q_0)$,  
$t_{\mu,0,0,1}\in \PGL(2,q)\setminus \PSL(2,q)$. Thus
$\phi^st_{\mu,0,0,1}\in M(s,q)$ and normalises $\PSL(2,q_0)$. Hence $N_{M(s,q)}(\PSL(2,q_0))=\la \PSL(2,q_0),\phi^st_{\mu,0,0,1}\ra$.

We will need the following lemma.

\begin{lemma}
\label{lem:subfieldedge}
Let $G=M(s,q)$ with $q=q_0^r=p^f$ for some odd primes $r$ and $p$, and even integer $f$, and let
$H=\la \PSL(2,q_0),\phi^st_{\mu,0,0,1}\ra$ where $\mu$ is a primitive
element of $\GF(q_0)$. Let $e=\{\{\infty,w_1\},\{\infty,w_2\}\}$. \begin{enumerate}
\item Then $G_e\leqslant H$ if and only if both $w_1+w_2$ and $(w_2-w_1)^{p^{2s}-1}$ lie in $\GF(q_0)$. 
\item There exist $w_1,w_2\notin\GF(q_0)$ such that $w_1+w_2$ and $(w_2-w_1)^{p^{2s}-1}$ lie in $\GF(q_0)$ if and only if $\gcd(\frac{q-1}{q_0-1},p^{2s}-1)\neq 1$.
\end{enumerate}
\end{lemma}
\proof
\begin{enumerate}
\item By Lemma \ref{lem:Gedge}, $G_e=\la t_{-1,w_1+w_2,0,1},(\phi^{2s})^g\ra$ where $g=t_{w_2-w_1,w_1,0,1}$.  Since $f$ is even and $q=q_0^r$ with $r$ odd, $q_0$ is an even power of $p$ and hence $-1$ is a square in $\GF(q_0)$. Thus $t_{-1,w_1+w_2,0,1}\in H$ if and only if $w_1+w_2\in\GF(q_0)$. Moreover, 
\begin{align*}
g^{-1}\phi^{2s}g &=t_{1,-w_1,0,w_2-w_1}\phi^{2s}t_{w_2-w_1,w_1,0,1}\\
 &=\phi^{2s}t_{1,-w_1^{p^{2s}},0,(w_2-w_1)^{p^{2s}}}t_{w_2-w_1,w_1,0,1}\\
&=\phi^{2s}t_{w_2-w_1,-(w_2-w_1)w_1^{p^{2s}}+w_1(w_2-w_1)^{p^{2s}},0,(w_2-w_1)^{p^{2s}}}\\
&=\phi^{2s}t_{1,w_1(w_2-w_1)^{p^{2s}-1}-w_1^{p^{2s}},0,(w_2-w_1)^{p^{2s}-1}}.
\end{align*}
Let $h=t_{1,w_1(w_2-w_1)^{p^{2s}-1}-w_1^{p^{2s}},0,(w_2-w_1)^{p^{2s}-1}}$.
As $\phi^{2s}\in H$, it follows that  $(\phi^{2s})^g \in H$ if and only if $h\in\PSL(2,q_0)$. Now if $h\in\PSL(2,q_0)$ then $(w_2-w_1)^{p^{2s}-1}\in\GF(q_0)$. Thus if $G_e\leqslant H$ then both $w_1+w_2$ and $(w_2-w_1)^{p^{2s}-1}$ lie in $\GF(q_0)$. Conversely, suppose that $w_1+w_2=a\in\GF(q_0)$ and $w_2-w_1=u$ with
$u^{p^{2s}-1}=b\in\GF(q_0)$. Then writing $\frac{1}{2}$ for $2^{-1}\in\GF(p)$ and noting that $2^{p^{2s}}=2$,
$w_1(w_2-w_1)^{p^{2s}-1}-w_1^{p^{2s}}
=\frac{a-u}{2}b-\frac{a^{p^{2s}}-u^{p^{2s}}}{2^{p^{2s}}}=\frac{ab-a^{p^{2s}}}{2}\in\GF(q_0)$. Thus $h\in\PGL(2,q_0)$, and since $p^{2s}-1$ is even $h\in\PSL(2,q_0)$. 
\item  Let $\xi$ be a primitive element of $\GF(q)$. Then $\GF(q_0)=\{0\}\cup\{\xi^{i\frac{q-1}{q_0-1}}|i=1,\ldots,q_0-1\}$ and we can choose $\mu=\xi^{\frac{q-1}{q_0-1}}$.
For $w_2-w_1=\xi^j\in \GF(q)\setminus\{0\}$, if $(w_2-w_1)^{p^{2s}-1}$ lies in $\GF(q_0)$, that means that  $\xi^{j(p^{2s}-1)}=\xi^{i\frac{q-1}{q_0-1}}$ for some integer $i$.
If  $\gcd(\frac{q-1}{q_0-1},p^{2s}-1)=1$, we must have $j$  a multiple of $(q-1)/(q_0-1)$, and so $w_2-w_1\in \GF(q_0)$. If we also have  $w_1+w_2\in\GF(q_0)$, then this implies that $w_1,w_2\in\GF(q_0)$. Hence if $w_1,w_2\notin\GF(q_0)$ such that $w_1+w_2$ and $(w_2-w_1)^{p^{2s}-1}$ lie in $\GF(q_0)$ then $\gcd(\frac{q-1}{q_0-1},p^{2s}-1)\neq 1$.
Conversely, suppose $\gcd(\frac{q-1}{q_0-1},p^{2s}-1)=d\neq 1$ and choose $j=(q-1)/d(q_0-1)$.
Then take $w_2=\xi^j/2$ and  $w_1=-\xi^j/2$. We obviously have $w_1+w_2\in\GF(q_0)$ and $w_1,w_2\notin\GF(q_0)$. Moreover $(w_2-w_1)^{p^{2s}-1}=\xi^{\frac{p^{2s}-1}{d}\frac{q-1}{q_0-1}}\in\GF(q_0)$.\Qed
\end{enumerate}

\begin{construction}
\label{con:Msqsubfield}
{\rm Let $G=M(s,q)$ and let $X=\{\infty\}\cup \GF(q)$ be the projective line. Let $q=q_0^r=p^f$ for some odd primes $r$ and $p$, and $f$ an even integer, and let $H=\la \PSL(2,q_0), \phi^st_{\mu,0,0,1}\ra$ where $\mu$ is a primitive element of $\GF(q_0)$. Assume $\gcd(\frac{q-1}{q_0-1},p^{2s}-1)\neq 1$, so that by Lemma \ref{lem:subfieldedge}, there exist $w_1,w_2\notin \GF(q_0)$ such that $w_1+w_2,(w_2-w_1)^{p^{2s}-1} \in\GF(q_0)$.  Let $e=\{\{\infty,w_1\},\{\infty,w_2\}\}$.  By Lemma \ref{lem:Gedge}, $G_e=\la t_{-1,w_1+w_2,0,1},(\phi^{2s})^g\ra$, where $g=t_{w_2-w_1,w_1,0,1}$, and by Lemma \ref{lem:subfieldedge}, $G_e\leqslant H$. Thus letting $P=e^H$ and $\P=P^G$, $(J(q+1,2),\P)$ is a $G$-primitive decomposition by Lemma \ref{lem:general}.

We claim that the divisors of $\P$ are either matchings or unions of cycles. Since $\PSL(2,q_0)_\infty\norml H_\infty=\la \PSL(2,q_0)_\infty,\phi^st_{\mu,0,0,1}\ra$, $H_\infty$ (of order $\frac{q_0(q_0-1)}{2}\frac{f}{s}$) acts on the set of $\PSL(2,q_0)_\infty$-orbits.  Now $t_{-1,w_1+w_2,0,1}\in \PSL(2,q_0)_{\infty}$ interchanges $w_1$ and $w_2$, and hence $w_1,w_2$ lie in the same $\PSL(2,q_0)_{\infty}$-orbit, $\theta$ say.
By Lemma \ref{lem:subfieldorbits}, $\PGL(2,q_0)$ acts semiregularly on $\GF(q)\backslash\GF(q_0)$ and hence $|\theta|=|\PSL(2,q)_{\infty}|=\frac{q_0(q_0-1)}{2}$. Note that $H_{\{\infty,w_1\}}=H_{\infty,w_1}$ and $H_{\{\infty,w_2\}}=H_{\infty,w_2}$. Also $H_{\{\infty,w_1\},\{\infty,w_2\}}=\la(\phi^{2s})^g\ra$ has order $\frac{f}{2s}$.
Notice that $(\phi^st_{\mu,0,0,1})^2=\phi^{2s}t_{\mu^{p^s+1},0,0,1}$. Hence 
$\la \PSL(2,q_0)_\infty,(\phi^{2s})^g\ra$ has index 2 in $H_\infty$ and fixes $\theta$. Therefore $H_\infty$ either fixes $\theta$ or switches it with another $\PSL(2,q_0)_\infty$-orbit $\theta'$.
In the first case, $H_{\{\infty,w_1\}}$ has order $\frac{f}{s}$, while $H_{\{\infty,w_1\},\{\infty,w_2\}}$ has order $\frac{f}{2s}$, hence the divisor has valency 2 and is a union of cycles.
In the second case, $H_{\{\infty,w_1\}}$ and $H_{\{\infty,w_1\},\{\infty,w_2\}}$ both have order $\frac{f}{2s}$, and so the divisor is a matching $\frac{q_0(q_o^2-1)}{2}K_2$.
}%
\end{construction}

\begin{remark}
{\rm We have not determined the length of the cycles occuring in the first case of Construction \ref{con:Msqsubfield}. This case happens if and only if there exists $w\in \GF(q)\setminus\GF(q_0)$ such that $w^{\phi^st_{\mu,0,0,1}}=w^{p^s}\mu\in \{a^2 w+b|a,b\in\GF(q_0)\}$. We have not been able to find any instances where this condition holds.
}%
\end{remark}

\begin{proposition}
Let $(J(q+1,2),\mathcal{P})$ be a $G$-primitive decomposition with
$G=M(s,q)$ and for $P\in\P$ we have that $G_P=N_G(\PSL(2,q_0))$ where
$q=q_0^r$ for some odd prime $r$. Then  $\P$ is obtained by
Construction $\ref{con:design}$ or $\ref{con:Msqsubfield}$.
\end{proposition}
\proof
Let $q=p^f$ with $p$ a prime and $f$ an even integer.
As seen in the discussion before Lemma \ref{lem:subfieldedge}, $H:=G_P=\la \PSL(2,q_0),\phi^st_{\mu,0,0,1}\ra$ where $\mu$ is a primitive element of $\GF(q_0)$. Let $X=\{\infty\}\cup\GF(q)$. Then one orbit of $H$ on $X$ is
$\{\infty\}\cup\GF(q_0)$. Since $H$ is maximal in $G$, $H$ is exactly the stabiliser in $G$ of $\{\infty\}\cup\GF(q_0)$.

Suppose that $H$ contains $G_e$ for some edge
$e=\{\{v,w_1\},\{v,w_2\}\}$. 

Then by Lemma \ref{lem:Gedge}, $H$ contains
an element of $\PSL(2,q)$, and hence of $\PSL(2,q_0)$, which fixes $v$
and interchanges $w_1$ and $w_2$. Since, by Lemma
\ref{lem:subfieldorbits}, $\PSL(2,q_0)$ acts semiregularly on
$\GF(q)\setminus\GF(q_0)$, it follows that
$v\in\{\infty\}\cup\GF(q_0)$. Without loss of generality we may suppose that
$v=\infty$.  By Lemma \ref{lem:subfieldedge}, this means that both  $w_1+w_2$ and $(w_2-w_1)^{p^{2s}-1}$ lie in $\GF(q_0)$.
This is of course satisfied if $w_1,w_2\in\GF(q_0)$, and then we get Construction \ref{con:design} using $\mathcal{B}=(\{\infty\}\cup\GF(q_0))^{\PGammaL(2,q)}$, as $G$ is transitive on $\mathcal{B}$.
Now assume $w_1,w_2\notin \GF(q_0)$. Then by Lemma \ref{lem:subfieldedge}, $\gcd(\frac{q-1}{q_0-1},p^{2s}-1)\neq 1$. Moreover, $P=e^H$ and $\P=P^G$ are as obtained in Construction \ref{con:Msqsubfield}.
\Qed


\begin{thebibliography}{1}

\bibitem{BPP}
J.~Bamberg, G.~Pearce and C.~E.~Praeger,
Transitive decompositions of graph products: rank 3 product action
type, J. Group Theory, to appear.

\bibitem{bosak}
J.~Bos\'ak, Decompositions of Graphs, Kluwer Academic Publishers, Dordrecht, The Netherlands, 1990.


\bibitem{magma}
W.~Bosma, J.~Cannon and C.~Playoust, The Magma algebra system
  I: The user language, J. Symb. Comp. 24 3/4 (1997),
235--265. Also see the {\sc Magma} home page at
http://www.maths.usyd.edu.au:8000/u/magma/. 

\bibitem{BCN}
A.~E.~Brouwer, A.~M.~Cohen and A.~Neumaier, Distance-regular
graphs. Ergebnisse der Mathematik und ihrer Grenzgebiete (3) [Results
  in Mathematics and Related Areas (3)], 18. Springer-Verlag, Berlin,
1989.

\bibitem{PSL}
F.~Buekenhout, P.~Cara, K.~Vanmeerbeek,
Geometries of the group $\PSL(2,11)$. Geom. Dedicata 83 (2000),
169--206.

\bibitem{Cameron}
P.~J.~Cameron, Permutation Groups, London Math. Soc, Student Texts 45, 1999.

\bibitem{CamKorch}
P.~J.~Cameron and G.~Korchm\'aros, One-factorizations of complete
graphs with a doubly transitive automorphism group.  
Bull. London Math. Soc. 25 (1993), 1--6.

\bibitem{PGLsubs}
P.~J.~Cameron, G.~R.~Omidi and B.~Tayeh-Rezaie, $3$-Designs from
$\PGL(2,q)$,  Electron. J. Combin. 13 (2006), \#R50.

\bibitem{handbook}
C.~J.~Colbourn and J.~H.~Dinitz (Editors), The CRC Handbook of
Combinatorial Designs, CRC Press Series on Discrete Mathematics and
its Applications, CRC Press, Boca Raton, FL, 1996.


\bibitem{atlas}
J.~H. Conway, R.~T. Curtis, S.~P. Norton, R.~A. Parker, and R.~A. Wilson,
 Atlas of finite groups, Clarendon Press, Oxford, 1985.


\bibitem{Cuaresma}
M.~C.~Cuaresma, Homogeneous Factorisations of Johnson Graphs, PhD Thesis, University of the Philippines, 2004.

\bibitem{CGP}
M.~C.~Cuaresma, M.~Giudici and C.~E.~Praeger, Homogeneous
factorisations of Johnson graphs, submitted.

\bibitem{DGLP}
A.~Devillers, M.~Giudici, C.H.~Li and C.~E.~Praeger,
A remarkable Mathieu graph tower, submitted.

\bibitem{dempwolff}
U.~Dempwolff, On the second chomology of $\GL(n,2)$,  J. Austral. Math. Soc.  16  (1973), 207--209.

\bibitem{dickson}
L.~E. Dickson, Linear groups: {W}ith an exposition of the
  {G}alois field theory, Dover Publications Inc., New York, 1958.

\bibitem{DM}
J.~D.~Dixon and B.~Mortimer, Permutation Groups, Graduate
Texts in Mathematics 163, Springer, New York, 1996.

\bibitem{FP}
X.~G.~Fang and C.~E.~Praeger, Finite two-arc-transitive graphs
  admitting a Suzuki simple group, Comm. Alg. 27 (1999),
3727--3754.


\bibitem{GPZ}
A.~Gardiner, C.~E.~Praeger  and S.~Zhou, Cross ratio graphs.
J. London Math. Soc. (2)  64  (2001), 257--272.

\bibitem{PSLmax}
Michael Giudici, Maximal subgroups of almost simple groups with
  socle $\PSL(2,q)$, preprint: math.GR/0703685.

\bibitem{GLPP1}
M.~Giudici, C.H.~Li, P.~Poto\v cnik and C.~E.~Praeger,
Homogeneous factorisations of graphs and digraphs. 
European J. Combin. 27 (2006), no. 1, 11--37.

\bibitem{genpaper}
M.~Giudici, C.~H.~Li and Cheryl E.~Praeger, Symmetrical covers, decompositions and factorisations of graphs, in preparation.

\bibitem{iso1}
F.~Harary, R.~W.~Robinson and N.~C.~Wormald, Isomorphic
factorisations. I. Complete graphs, Trans. Amer. Math. Soc. 242
(1978), 243--260. 

\bibitem{iso10}
F.~Harary and R.~W.~ Robinson, Isomorphic factorisations X: unsolved
problems, J. Graph Theory 9 (1985), 67--86. 

\bibitem{heinrich}
K.~Heinrich, Graph decompositions and designs, in: The CRC Handbook of
Combinatorial Designs, Charles J. Colbourn and Jeffrey H. Dinitz,
(Editors), CRC Press Series on Discrete Mathematics and its
Applications, CRC Press, Boca Raton, FL, 1996, 361--366.

\bibitem{huppert1}
B.~Huppert, Endliche {G}ruppen. {I}, in \emph{Die Grundlehren der
  Mathematischen Wissenschaften, Band 134} (Springer-Verlag, Berlin, 1967).

\bibitem{sashasbook}
A.~A.~Ivanov, Geometry of sporadic groups I: Petersen and tilde
geometries, Cambridge University Press, Cambridge, 1999.

\bibitem{KK}
L.~A.~Kalu\v znin and M.~H.~Klin, On some maximal subgroups of
symmetric and alternating groups, Mat. Sbornik 87 (1972), 91--121 (in
Russian).

\bibitem{Kantor}
W.~M.~Kantor, $k$-homogeneous groups, Math. Z. 124 (1972), 261--265.

\bibitem{LP03}
C.~H.~Li and C.~E.~Praeger, On partitioning the orbitals of a
transitive
permutation group, Trans. Amer. Math. Soc. 355 (2003),
 637--653.

\bibitem{LPS}
M.~W.~Liebeck, C.~E.~Praeger and J.~Saxl, 
A classification of the maximal subgroups of the finite alternating and symmetric groups.
J. Algebra 111 (1987), no. 2, 365--383.

\bibitem{LiWa}
D.~Livingstone and A.~Wagner, Transitivity of finite permutation groups on unordered sets, Math. Z.
 90 (1965), 393--403.

\bibitem{orrisonetal}
D.~K.~Maslen, M.~E.~Orrison and D.~N.~Rockmore,
Computing isotypic projections with the Lanczos iteration,
SIAM J. Matrix Anal. Appl. 25 (2004), 784--803.

\bibitem{graphatlas}
R.~C.~Read and R.~J.~Wilson, 
An atlas of graphs.
Oxford Science Publications. The Clarendon Press, Oxford University
Press, New York, 1998.

\bibitem{wielandt}
H.~Wielandt, Finite permutation groups.
Academic Press, New York-London, 1964. 

\bibitem{sibley}
T.~Sibley, On classifying finite edge colored graphs with two
transitive automorphism groups. J. Combin. Theory Ser. B 90 (2004),
121--138. 

\bibitem{seitz}
G.~M.~Seitz, Flag-transitive subgroups of Chevalley groups.  Ann. of Math. (2)  97  (1973), 27--56.

\end{thebibliography}
\end{document}